\documentclass[11pt]{amsart}

\usepackage{amsmath,amssymb,latexsym,amsthm,newlfont,enumerate}%,cmap}
\usepackage[
top=2cm,
bottom=2cm,
left=2.5cm,
right=2.5cm,
headheight=10pt, % as per the warning by fancyhdr
includehead,includefoot,
heightrounded, % to avoid spurious underfull messages
]{geometry}
\usepackage{graphicx}
\usepackage{adjustbox}
\usepackage[all]{xy}
\usepackage{mathtools}
\usepackage{hyperref}
\hypersetup{colorlinks=true, linkcolor=red}
\usepackage{multicol}

\theoremstyle{plain}

\newtheorem{thm}{Theorem}[section]

\newtheorem{pro}[thm]{Proposition}

\newtheorem{lem}[thm]{Lemma}

\newtheorem{conj}{Conjecture}

\theoremstyle{definition}

\newtheorem{dfn}[thm]{Definition}

\newtheorem{que}[thm]{Question}

\newtheorem{rem}[thm]{Remark}

\newtheorem{exa}[thm]{Example}

\DeclareMathOperator{\Hom}{Hom}

\DeclareMathOperator{\Ext}{Ext}

\newcommand{\E}{\mathcal{E}}

\newcommand{\Z}{\mathbb{Z}}
\newcommand{\PP}{\mathbb{P}}
\newcommand{\D}{\mathrm{D}}

\newcommand{\A}{\mathcal{A}}

\newcommand{\cone}{\mathrm{Cone}}

\newcommand{\Sym}{\mathrm{Sym}}

\newcommand{\OO}{\mathcal{O}}

\newcommand{\Pic}{\mathrm{Pic}}

\title[Classification of full exceptional collections]{Classification of full exceptional collections on smooth toric Fano varieties with Picard rank two}
\date{\today}
\author{Dae-Won Lee}
\address{Department of Mathematics, Yonsei University, Seoul 120-749, Republic of Korea}
\email{daewonlee@yonsei.ac.kr}

\begin{document}

\begin{abstract} 
	In this paper, we give the complete classification of full exceptional collections on smooth toric Fano threefolds and fourfolds with Picard rank two. To be precise, we give a partial answer to the conjecture in \cite{Kuz} and \cite{LYY}: we provide all the exceptional collections of maximal length and prove that they are in fact full exceptional collections.
\end{abstract}

\maketitle

\section{Introduction}
 
 In this paper, the varieties are over algebraically closed field $k$ of characteristic zero.
 
 In recent years, the derived category of coherent sheaves has become one of the most important invariants in the study of algebraic geometry. Let $X$ be a smooth projective variety and we denote by $\D(X)$ the bounded derived category of coherent sheaves on $X$. Then $\D(X)$ contains geometric information on $X$. For example, the variety $X$ can be reconstructed from $\D(X)$ if $K_X$ is ample or antiample \cite{BO}. Many other deep connections between $X$ and $\D(X)$ are collected in the book \cite{Huy}.
 In the derived category of coherent sheaves of a variety, the notion of full exceptional collection and semiorthogonal decomposition is one of the most imperative objects. 
 
 There are well-known examples of semiorthogonal decompositions. The simplest one is that of projective spaces:
 
 \begin{thm} \cite{Bei}
 	Let $X$ be the projective space $\PP^n$. Then there is a semiorthogonal decomposition;
 	\begin{align}
 	\D(X)=\langle\OO_X,\OO_X(1),\cdots,\OO_X(n)\rangle. \nonumber
 	\end{align}
 \end{thm}
 
 This is the earliest example that is discovered by Beilinson. 
 We also recall that the semiorthogonal decompositions for Grasmannian and quadric are known \cite{Kap}.
 
 Note that an exceptional collection $\{E_1,\cdots,E_n\}$ on $\D(X)$ gives a semiorthogonal decomposition for $\D(X)$ of the form $\D(X)=\langle \mathcal{A},E_1,\cdots,E_n\rangle$, where $\mathcal{A}:=\langle E_1,\cdots,E_n\rangle^\perp$. Here, $\mathcal{A}$ is trivial if and only if the exceptional collection $\{E_1,\cdots,E_n\}$ is full. Thus, a natural and important question arises: \textit{When is an exceptional collection full?}
 
 One can na\"{i}vely expect that an exceptional collection of maximal length is full. However, recent studies show that phantom or quasi-phantom categories give a counterexample to such expectation \cite{GKMS}, \cite{GS}. Nonetheless, one may still hope the following conjecture is true, which is an interesting problem itself:
 
 \begin{conj} \cite{Kuz} \label{Kuzc}
 	Let $X$ be a smooth projective variety. If $D(X)$ has a full exceptional collection of line bundles of length $n$, then any exceptional collection of line bundles of length $n$ is full.
 \end{conj}
 
  The conjecture above is known to be true for $\PP^n$ \cite{Bei}, del Pezzo surfaces \cite{KO}, weak del Pezzo surfaces \cite{Kul}, Hirzebruch surfaces \cite{Hil} and smooth projective toric surface $X$ of Picard rank $3$ or $4$ \cite{HI}. It is also true for the following blow-ups of $\PP^3$; blow-up of $\PP^3$ at a point, a line or a cubic curve \cite{LYY}. In their work, they give a conjecture:
 
 \begin{conj} \cite{LYY} \label{Wanc}
 	All the exceptional collections of line bundles of maximal length on smooth toric Fano 3-folds and 4-folds are full.
 \end{conj}
  
 The following is our main result of this paper. It gives a partial positive answer to the above conjecture.
 
 \begin{thm} \label{main}
 	Let $X$ be any smooth toric Fano threefold or a toric Fano fourfold with Picard rank two. Then all the exceptional collections of maximal length consisting of line bundles on $X$ are full.
 \end{thm}
 
 We prove Theorem \ref{main} by finding the cohomologically zero line bundles and classifying all the exceptional collections of maximal length. After that, we show that they are in fact full. Our classification gives the following eleven cases; two of them are threefolds and the other nine of them are fourfolds. See Table \ref{tab:table14}. 
 
  	\begin{table}[h!] 
 	\begin{center}
 		\caption{Smooth toric Fano $3$-folds and $4$-folds with Picard rank $2$}
 		\label{tab:table14}
 		\begin{adjustbox}{max width=\textwidth}
 			\begin{tabular}{c|c|c}
 				Variety & Dimension & Maximal length \\
 				\hline
 				$\PP_{\PP^1}(\OO\oplus\OO\oplus\OO(1))$ &	$3$ & $6$\\
 				\hline
 				$\PP_{\PP^2}(\OO\oplus\OO(1))$ & $3$  & $6$ \\
 				\hline
 				$\PP_{\PP^2}(\OO\oplus\OO(2))$ &	$3$ & $6$\\
 				\hline
 				$\PP^2 \times \PP^1$ & $3$	 & $6$ \\
 				\hline
 				$\PP_{\PP^3}(\OO\oplus\OO(3))$ &	$4$ & $8$ \\
 				\hline
 				$\PP^1 \times \PP^3$ & $4$ & $8$  \\
 				\hline
 				$\PP_{\PP^1}(\OO\oplus\OO\oplus\OO\oplus\OO(1))$ & $4$ & $8$ \\
 				\hline
 				$\PP_{\PP^2}(\OO\oplus\OO\oplus\OO(2))$ & $4$ & $9$ \\
 				\hline
 				$\PP_{\PP^2}(\OO\oplus\OO\oplus\OO(1))$ & $4$ & $9$ \\
 				\hline
 				$\PP_{\PP^2}(\OO\oplus\OO(1)\oplus\OO(1))$ & $4$ & $9$ \\
 				\hline
 				$\PP^2 \times \PP^2$ & $4$ & $9$ \\
 				\hline
 				
 			\end{tabular}
 		\end{adjustbox}
 	\end{center}
 \end{table}
 
 There are originally $4$ smooth toric Fano threefolds with Picard rank two. However, \cite{LYY} considered two of them; blowup of $\PP^3$ at a point and at a line, which are the uppermost two varieties in Table \ref{tab:table14}. Hence for smooth toric Fano threefold, we only need to consider $\PP_{\PP^2}(\OO\oplus\OO(2))$ and $\PP^2 \times \PP^1$. For smooth toric Fano fourfolds, we consider $\PP_{\PP^3}(\OO\oplus\OO(3))$, $\PP_{\PP^3}(\OO\oplus\OO(2))$, $\PP_{\PP^3}(\OO\oplus\OO(1))$, $\PP^1 \times \PP^3$, $\PP_{\PP^1}(\OO\oplus\OO\oplus\OO\oplus\OO(1))$, $\PP_{\PP^2}(\OO\oplus\OO\oplus\OO(2))$, $\PP_{\PP^2}(\OO\oplus\OO\oplus\OO(1))$, $\PP_{\PP^2}(\OO\oplus\OO(1)\oplus\OO(1))$ and $\PP^2 \times \PP^2$. We classify all the exceptional collections of maximal length for the varieties above and show that they are full exceptional collections. 
 
 In general, it is hard to tell whether a given smooth projective variety has full exceptional collections or not and even if it is, classifying all the full exceptional collections is still more difficult. Nontheless, by proving Theorem \ref{main}, we are able to give explicit classification of the full exceptional collections of all smooth toric Fano threefolds and fourfolds with Picard rank two.
 
 After the first version of this paper was announced, Professor Alastair Craw introduced to the author of the paper \cite{Na} which constructs tilting bundles on every smooth toric Fano fourfolds by using the strong full exceptional collections of line bundles.
 
 This paper is organized as follows. In section $2$, we give a brief introduction and preliminary results on the derived category especially on the notion of semiorthogonal decomposition, full exceptional collection and mutation. In section $3$, we prove that all the exceptional collections of maximal length of line bundles on smooth toric Fano threefolds of Picard rank two are full. In section $4$, we prove the same statements for smooth toric Fano fourfolds of Picard rank two.
 
 \section*{Acknowledgments}
 The author would like to express his deep gratitute to his advisor Sung Rak Choi for introducing this topic and many valuable suggestions. He also thanks Professor Alastair Craw for introducing the paper \cite{Na}.
 
\section{Preliminary}

In this section, we first give definitions and basic facts on derived category.

Let $\mathcal{T}$ be a triangulated category. By axioms of triangulated category, for any morphism $f:A\rightarrow B$ in $\mathcal{T}$, there exists a distinguished triangle
 \begin{align}
 A\rightarrow B\rightarrow C\rightarrow A[1], \nonumber
 \end{align}
 where $C\cong \cone(f)$ is the $\emph{cone}$ of the morphism $f:A\rightarrow B$.

\begin{dfn}
	A \emph{semiorthogonal decomposition} of $\mathcal{T}$ is a collection of full triangulated subcategories $\{\A_1,\A_2,\cdots,\A_n\}$ of $\mathcal{T}$ such that
	\begin{itemize}
		\item [(1)] $\Hom(A_i,A_j)=0$ for all $1\leq j<i\leq n $ and any objects $A_i\in \A_i$, $A_j\in \A_j$ and
		\item [(2)] the smallest triangulated subcategory containing the collection $\{\A_1,\A_2,\cdots,\A_n\}$ coincides with $\mathcal{T}$.
	\end{itemize}
\end{dfn}

We will frequently use the following Theorem:

\begin{thm} (Orlov's projective bundle formula) \cite{Orl} \label{proj}
	Let $f:X:=\PP(\E)\rightarrow Y$ be the projective bundle of $\E$, where $\E$ is a vector bundle of rank $r+1$ on a smooth projective variety $Y$. If $\D(Y)$ has a full exceptional collection, then there is a semiorthogonal decomposition 
	\begin{align}
		\D(X)=\langle f^*\D(Y),f^*\D(Y)\otimes\OO_X(1),\cdots,f^*\D(Y)\otimes \OO_X(r)\rangle, \nonumber
	\end{align}
	where $\OO_X(1)$ is the Grothendieck bundle of $X$.
\end{thm}
 
 Next we recall the definition of exceptional collection and what it means to be full.
 
 \begin{dfn}
 	
 	\begin{itemize}
 		\item [(1)] An object $E\in \mathcal{T}$ is called an \emph{exceptional object} if 
 		\begin{align}
 		\Hom(E,E[m])=
 		\begin{cases}
 		k & \text{if } m=0,\\
 		0 & \text{if } m\neq 0.
 		\end{cases} \nonumber
 		\end{align}
 		\item [(2)] A sequence of exceptional objects $\{E_1,E_2,\cdots,E_n\}$ is called an \emph{exceptional collection} if 
 		\begin{align}
 		\Hom(E_i,E_j[m])=0 \text{ for any } 1\leq j<i\leq n \text{ and for all } m\in \Z.  \nonumber		
 		\end{align}
 		\item [(3)] An exceptional collection $\{E_1,E_2,\cdots,E_n\}$ is said to be \emph{full} if they generates the whole category $\mathcal{T}$.
 	\end{itemize}
 \end{dfn}

 There is a weaker notion of exceptional collection, which we call numerically exceptional collection.

 \begin{dfn}
 	\begin{itemize}
 		\item [(1)] An object $E\in \D(X)$ is said to be \emph{numerically exceptional} if $\chi(E,E)=1$, where $\chi(\mathcal{E},\mathcal{F})=\sum_s (-1)^s \mathrm{dimExt}^s(\mathcal{E},\mathcal{F})$.
 		\item [(2)] Objects $E_1,E_2,\cdots,E_n$ in $\D(X)$ is said to be \emph{numerically exceptional collection} if 
 		\begin{itemize}
 			\item [(a)] every $E_i$ is numerically exceptional and
 			\item [(b)] $\chi(E_i,E_j)=0$ for all $1\leq j<i\leq n$.
 		\end{itemize}
 	\end{itemize}
 \end{dfn}
 
 \begin{dfn}
 	A (numerically) exceptional collection $E_1,E_2,\cdots,E_n$ in $\D(X)$ is of \emph{maximal length} if $n=\mathrm{rank}(K_0/\mathrm{ker}\chi)$.
 \end{dfn}
 
 Note that any full exceptional collection is of maximal length.  The following Lemma is useful to the proof of our main result.
 
 \begin{lem}\cite{LYY}\label{mt1}
 	Let $X$ be a smooth projective variety and $\D(X)$ the bounded derived category of coherent sheaves on $X$. If $\D(X)$ has two exceptional collections of the same length with only one compoent different, then both are full or neither is.
 \end{lem}	

 \begin{dfn}
 	Let $X$ be a smooth projective variety. A \emph{cohomologically zero line bundle} is a line bundle $L\in \Pic(X)$ such that $H^i(X,L)=0$ for all $i\in \Z$.
 \end{dfn}

 \begin{exa}
 	\begin{itemize}
 		\item [1.] $\OO_{\PP^2 \times \PP^1}(a,b)$ is cohomologically zero if and only if $a=-1,-2$ or $b=-1$;
 		\item [2.] $\OO_{\PP^3 \times \PP^1}(a,b)$ is cohomologically zero if and only if $a=-1,-2,-3$ or $b=-1$;
 		\item [3.] $\OO_{\PP^2 \times \PP^2}(a,b)$ is cohomologically zero if and only if $a=-1,-2$ or $b=-1,-2$.
 	\end{itemize}
 \end{exa}
 
 \begin{rem} \label{LYY}
 	Let $X$ be a smooth projective variety. For any two line bundles $L_1, L_2$ on $X$, there are isomorphisms $\Hom(L_1,L_2[i])\cong \Ext^i(L_1,L_2)\cong H^i(X,L_2\otimes L_1^{-1})$. Hence, for line bundles $L_1,L_2,\cdots,L_n$ on $X$, they form an exceptional collection if and only if $L_j\otimes L_i^{-1}$ are cohomologically zero line bundles for all $1\leq j<i\leq n$.
 \end{rem}

 Therefore, finding cohomologically zero line bundles is essential in order to give a classification of exceptional collections consisting of line bundles of maximal length.
 
 \begin{lem}\cite{LYY} \label{nor}
 	Let $X$ be a smooth projective variety. Then the sequences $\{L_1,L_2,\cdots,L_n\}$ and $\{\OO_X, L_2\otimes L_1^{-1},\cdots,L_n\otimes L_1^{-1}\}$ are both (full) exceptional collections or neither is.
 \end{lem}
 
 The process in Lemma \ref{nor}, making $\{L_1,L_2,\cdots,L_n\}$ into $\{\OO_X, L_2\otimes L_1^{-1},\cdots,L_n\otimes L_1^{-1}\}$  will be often called as normalization.

 Now, we recall the definition of the mutation and its property.
 
 	Let $\mathcal{T}$ be a triangulated category and $E\in \mathcal{T}$ be an exceptional object. Then for any object $F\in \mathcal{T}$, consider two triangles:
 	\begin{align}
 	\mathrm{RHom}(E,F)\otimes E\rightarrow F\rightarrow L_E(F) ~\text{and} \nonumber\\
 	R_E(F)\rightarrow F\rightarrow \mathrm{RHom}(F,E)^*\otimes E,\nonumber
 	\end{align}
 	which are the projections of $F$ onto the left and the right orthogonals to $E$.
 	\begin{dfn}
 	We call 
 	\begin{align}
 	L_E(F)\cong \cone(\mathrm{RHom}(E,F)\otimes E\rightarrow F)  ~and~ R_E(F)\cong \cone(F\rightarrow \mathrm{RHom}(F,E)^*\otimes E)[-1],\nonumber
 	\end{align}
 	the \emph{left} and the \emph{right mutation functors}, respectively. 
   \end{dfn}

  Suppose $\mathcal{T}=\langle E_1, E_2,\cdots, E_n\rangle$ is a semiorthogonal decomposition. Then we have semiorthogonal decompositions:
  \begin{align}
  \mathcal{T}=\langle E_1,E_2,\cdots,E_{i-1},L_{E_i}(E_{i+1}),E_i,E_{i+1},\cdots,E_n\rangle ~\text{for}~ 1\leq i\leq n-1\nonumber
  \end{align}
  and
  \begin{align}
  \mathcal{T}=\langle E_1,E_2,\cdots,E_{i-2},E_i,R_{E_i}(E_{i-1}),E_{i+1},\cdots,E_n\rangle ~\text{for}~ 2\leq i\leq n. \nonumber
  \end{align}
 
 In particular, one has the following lemma:
 \begin{lem} \cite{Bo90}, \cite{Har}\label{mt}
 	\begin{itemize}
 		\item [(1)] Let $E,F\in \mathcal{T}$ be an exceptional pair, i.e., $\{E,F\}$ is an exceptional collection. If $\{F,E\}$ is also an exceptional pair, then $L_E(F)\cong F$ and $R_F(E)\cong E$.
 		\item [(2)] Let $\D(X)=\langle\A_1,\A_2,\cdots,A_n\rangle$ be an full exceptional collection in a derived category of smooth projective variety $X$. Then the following two collections 
 		\begin{align}
 		\langle\A_2,\cdots,\A_{n-1},\A_n,\A_1\otimes \omega_X^{-1}\rangle ~and~ \langle\A_n\otimes \omega_X,\A_1,\A_2,\cdots,\A_{n-1}\rangle\nonumber
 		\end{align}
 		are also full exceptional collections in $\D(X)$.
 	\end{itemize}
 \end{lem}

 Below is the Hirzebruch-Riemann-Roch Theorem for smooth projective fourfolds.

 \begin{lem}\label{HRR}
 	Let $X$ be a smooth projective variety with $dim X=4$ and $D$ be a divisor on $X$. Then
 	\begin{align}
 	\chi(X,\OO(D))&=\frac{1}{720}(-c_1^4+4c_1^2c_2+3c_2^2+c_1c_3-c_4)+\frac{1}{24}Dc_1c_2\nonumber\\
 	&\qquad +\frac{1}{24}D^2(c_1^2+c_2)+\frac{1}{12}D^3c_1+\frac{1}{24}D^4. \nonumber
 	\end{align}
 	\begin{proof}
 	Note that $\chi(X,E)=\int_X \mathrm{ch}(E)\cdot \mathrm{td}(T_X)$ by Hirzebruch-Riemann-Roch Theorem. If $E=\OO_X(D)$ is a line bundle, we have
 	\begin{align}
 	\mathrm{ch}(\OO_X(D))=1+D+\frac{1}{2}D^2+\frac{1}{6}D^3+\frac{1}{24}D^4. \nonumber
 	\end{align}
 	Also, since 
 	\begin{align}
 	\mathrm{td}(T_X)&=1+\frac{1}{2}c_1+\frac{1}{12}(c_1^2+c_2)+\frac{1}{24}c_1c_2 \nonumber\\
 	&\qquad +\frac{1}{720}(-c_1^4+4c_1^2c_2+3c_2^2+c_1c_3-c_4),\nonumber
 	\end{align}
 	we have the required result.
 	\end{proof}
 \end{lem}
 
 Next we list the varieties that we will consider:
 
  Smooth toric Fano threefolds with Picard rank two \cite{WW}:
  	\begin{itemize}
  		\item [1.] $X\cong \PP_{\PP^2}(\OO\oplus\OO(2))$
  		\item [2.] $X\cong \PP^2 \times \PP^1$
  	\end{itemize}
 
 Smooth toric Fano fourfolds with Picard rank two \cite{Bat}:
 \begin{multicols}{2}
 	 \begin{itemize}
 		\item [3.] $X\cong \PP_{\PP^3}(\OO\oplus\OO(3))$
 		\item [4.] $X\cong \PP_{\PP^3}(\OO\oplus\OO(2))$
 		\item [5.] $X\cong \PP_{\PP^3}(\OO\oplus\OO(1))$
 		\item [6.] $X\cong \PP^1 \times \PP^3$
 		\item [7.] $X\cong \PP_{\PP^1}(\OO\oplus\OO\oplus\OO\oplus\OO(1))$
 		\item [8.] $X\cong \PP_{\PP^2}(\OO\oplus\OO\oplus\OO(2))$
 		\item [9.] $X\cong \PP_{\PP^2}(\OO\oplus\OO\oplus\OO(1))$
 		\item [10.] $X\cong \PP_{\PP^2}(\OO\oplus\OO(1)\oplus\OO(1))$
 		\item [11.] $X\cong \PP^2 \times \PP^2$
 	\end{itemize}
 \end{multicols}

 Finally we briefly explain our strategy: we first find all the cohomologically zero line bundles on each variety. By the armument in Remark \ref{LYY}, we classify all the exceptional collections of maximal length consisting of line bundles. After that, by using Lemma \ref{proj} and Lemma \ref{mt}, we show that they are full exceptional collections of line bundles.

\section{Smooth toric Fano threefolds with Picard rank 2}
\subsection{When $X\cong \PP_{\PP^2}(\OO\oplus\OO(2))$.} 
Note that $\PP_{\PP^2}(\OO\oplus\OO(2))\cong \PP_{\PP^2}(\OO(-1)\oplus \OO(1))$. Let $f:\PP_{\PP^2}(\OO(-1)\oplus \OO(1))\rightarrow \PP^2$ be the projective bundle. Then its canonical divisor is $K_X=-3H-2D$ where $H$ is the pullback of hyperplane class in $\PP^3$ and $D$ is the class of $c_1(\OO_X(1))$. The Picard group of $X$ is 
\begin{align}
\Pic(X)\cong \Pic(\PP^2)\oplus \mathbb{Z}[D]=\mathbb{Z}[H]\oplus\mathbb{Z}[D], \nonumber
\end{align}
with intersection numbers
\begin{align}
H^3=0, ~H^2D=1, ~HD^2=0, ~HD^3=1. \nonumber
\end{align}
\subsubsection{Cohomologically zero line bundles}
\begin{lem}\label{CZ1}
	$H^0(X,aH+bD)=0$ if and only if $b<0$ or $a+b<0, ~b\geq 0$. Also, $H^3(X,aH+bD)=0$ if and only if $b>-2$ or $a+b>-5, ~b\leq -2$.
	\begin{proof}
		Note that we have the following isomorphisms:
		\begin{align}
		f_*(\OO_X(b))&\cong \Sym^b(\OO_{\PP^2}(-1)\oplus\OO_{\PP^2}(1)) \nonumber \\ 
		&\cong \bigoplus_{i=0}^b \Sym^i(\OO_{\PP^2}(-1))\otimes \Sym^{b-i}(\OO_{\PP^2}(1)) \nonumber \\
		&\cong \bigoplus_{i=0}^b \OO_{\PP^2}(b-2i) ~\text{if}~ b\geq 0, \nonumber
		\end{align} and $f_*(\OO_X(b))\cong 0$ if $b<0$.
		
		By projection formula, we have the following isomorphisms:
		\begin{align}
		H^0(X, f^*(\OO_{\PP^2}(a))\otimes \OO_X(b))&\cong H^0(\PP^2, \OO_{\PP^2}(a)\otimes  \bigoplus_{i=0}^b \OO_{\PP^2}(b-2i)) \nonumber \\
		&\cong \bigoplus_{i=0}^b H^0(\PP^2,\OO_{\PP^2}(a+b-2i)) ~\text{if}~ b\geq 0,\nonumber
		\end{align}
		and if $b<0$, then $H^0(X,aH+bD)=0$.
		Hence, $H^0(X,aH+bD)=0$ if and only if $b<0$ or $a+b<0, ~b\geq 0$.
		
		By Serre duality, one obtains the condition for $H^3(X,aH+bD)=0$.
	\end{proof}
\end{lem}
\begin{lem}
	$h^1(X,aH+bD)h^2(X,aH+bD)=0$ for all $a,b \in \mathbb{Z}$.
	\begin{proof}
		By Leray's spectral sequence and projection formula, if $b>-2$ then we have 
		\begin{align}
		H^1(X,f^*(\OO_{\PP^2}(a))\otimes \OO_X(b))\cong H^1(\PP^2,\OO_{\PP^2}(a)\otimes f_*\OO_X(b)). \nonumber
		\end{align}
		Since $f_*(\OO_X(b))\cong \Sym^b(\OO_{\PP^2}(-1)\oplus\OO_{\PP^2}(1))\cong \bigoplus_{i=0}^b \OO_{\PP^2}(b-2i)$ if $b\geq 0$, and $0$ otherwise, we have $H^1(X,aH+bD)=0$ for all $a\in \mathbb{Z}, b\geq-1$.
		
		By Serre duality, we have $H^2(X,aH+bD)=0$ for all $a\in \mathbb{Z}, b\leq -1$. Thus, either one of $H^1(X,aH+bD)$ or $H^2(X,aH+bD)$ is zero and this completes the proof.
	\end{proof}
\end{lem}
\begin{pro} \label{pro1}
	A line bundle $\OO_X(aH+bD)$ is cohomologically zero if one of the following holds:
	\begin{itemize}
		\item [(1)] $b=-1$,
		\item [(2)] $a=-2, b=-2$,
		\item [(3)] $a=-2, b=0$,
		\item [(4)] $a=-1, b=-2$,
		\item [(5)] $a=-1, b=0$.
	\end{itemize}
	\begin{proof}
		Since $c(X)=(1+3H+3H^2)(1+2D)$, we have $c_1=3H+2D$ and $c_2=6HD+3H^2$.
		Thus, by Hirzebruch-Riemann-Roch formula, we have
		\begin{align}
		\chi(\OO_X(aH+bD))&=\frac{1}{6}(6+9a+8b+9ab+3a^2+3b^2+3a^2b+b^3) \nonumber\\
		&=\frac{1}{6}(b+1)(3a^2+9a+b^2+2b+6). \nonumber
		\end{align}
		Therefore, by Lemma \ref{CZ1}, $\OO_X(aH+bD)$ is cohomologically zero if the following conditions hold:
		\begin{itemize}
			\item [(i)] $b<0$ or $a+b<0, ~b\geq 0$,
			\item [(ii)] $b>-2$ or $a+b>-5, ~b\leq -2$,
			\item [(iii)] $(b+1)(3a^2+9a+b^2+2b+6)=0$.
		\end{itemize}
		Note that the equation $3a^2+9a+b^2+2b+6=0$ has $4$ integer solutions; $(a,b)=(-2,-2), (-2,0),$ $(-1,-2), (-1,0)$. This completes the proof.
	\end{proof}
\end{pro}

\subsubsection{Classification of full exceptional collections}
\begin{thm} \label{cl1}
	Let $X$ be the projective bundle $f:\PP_{\PP^2}(\OO(-1)\oplus \OO(1))\rightarrow \PP^2$. Then the normalized sequence
	\begin{align}
	\{\OO_X,\OO_X(D_1),\cdots,\OO_X(D_5)\} \label{1} \tag{A}
	\end{align}
	is an exceptional collection of line bundles if and only if the ordered set of divisors $\{D_1,D_2,\cdots,D_5\}$ is one of the following types:
	\begin{itemize}
		\item [(1)] $\{H,2H,aH+D,(a+1)H+D,(a+2)H+D\}$,
		\item [(2)]
		$\{H,aH+D,(a+1)H+D,(a+2)H+D,2H+2D\}$
		\item [(3)] $\{aH+D,(a+1)H+D,(a+2)H+D,H+2D,2H+2D\}$,
	\end{itemize}
	where $a\in \mathbb{Z}$. Moreover, by mutations and normalizations, they are related as:
	\begin{align}
	(1)\rightarrow (2)\rightarrow (3)\rightarrow (1). \nonumber
	\end{align}
	\begin{proof}
		Write $D_0=0$. The sequence (\ref{1}) is an exceptional collection if and only if for any integers $0\leq i<j\leq 5$, the line bundles $\OO_X(D_i-D_j)$ are cohomologically zero.
		
		To make the sequence (\ref{1}) an exceptional collection, $\OO_X(D_i)$ must be one of the line bundles in Proposition \ref{pro1}. To find out all the possible exceptional collections, we need Table \ref{tab:table1}.
		
		From Table \ref{tab:table1}, we observe that $\OO_X(D_1)$ may be one of $B_0$ and $B_4$. To finish the proof, we shall discuss case by case:
		\begin{itemize}
			\item [Case1.] Suppose $\OO_X(D_1)=B_0$. In this case, the only possible combination is $\{aH+D,(a+1)H+D,(a+2)H+D,H+2D,2H+2D\}$, which is $(3)$ in Theorem \ref{cl1}.
			\item [Case2.] Suppose $\OO_X(D_1)=B_4$. In this case, there are two possible choices for $D_2$:
			\begin{itemize}
				\item [(a)]
				If $\OO_X(D_2)=B_0$, we have $\{H,aH+D,(a+1)H+D,(a+2)H+D,2H+2D\},$ which is $(1)$ in Theorem \ref{cl1},
				\item [(b)] If $\OO_X(D_2)=B_2$, we have $\{H,2H,aH+D,(a+1)H+D,(a+2)H+D\},$ which is $(2)$ in Theorem \ref{cl1}.
			\end{itemize}
		\end{itemize}
		
		\begin{table}[h!] 
			\begin{center}
				\caption{Exceptional pairs.}
				\label{tab:table1}
				\begin{adjustbox}{max width=\textwidth}
				\begin{tabular}{c|c|c|c|c|c}
					& $B_0'$ & $B_1$ & $B_2$ & $B_3$ & $B_4$\\
					\hline
					$B_0$ &	$a'=a+1, a+2$ & \checkmark  &  &  \checkmark  &  \\
					\hline
					$B_1$ &	 &  &  &  & \\
					\hline
					$B_2$ &	\checkmark &  &  &  & \\
					\hline
					$B_3$ &  & \checkmark &  &  &  \\
					\hline
					$B_4$ & \checkmark & \checkmark & \checkmark &  & \\
					\hline
				\end{tabular}
			
			\end{adjustbox}
			\vspace{3mm}
			
			where $B_0=aH+D, B_1=2H+2D, B_2=2H, B_3=H+2D, B_4=H$.
			\end{center}
		\end{table}
	\end{proof}
\end{thm}

\begin{thm}
	Let $X$ be the projective bundle $f:\PP_{\PP^2}(\OO(-1)\oplus\OO(1))\rightarrow \PP^2$. Then any exceptional collection of line bundles of length 6 on $X$ is full.
	\begin{proof}
		By Lemma \ref{nor}, we only need to show that any exceptional collection of line bundles of length $6$ in Theorem \ref{cl1} is full. By Lemma \ref{mt}, it suffices to show that the exceptional collection of type $(1)$ in Theorem \ref{cl1} is full. 
		
		By Beilinson's semiorthogonal decomposition, we have 
		\begin{align}
		\D(\PP^2)=\langle\OO_{\PP^2}, \OO_{\PP^2}(1), \OO_{\PP^2}(2)\rangle. \nonumber
		\end{align}
		By Orlov's projective bundle formula (Theorem \ref{proj}), we obtain the following semiorthogonal decomposition of $\D(X)$
		\begin{align}
		\D(X)&=\langle f^*\D(\PP^2),f^*\D(\PP^2)\otimes \OO_X(1)\rangle \nonumber\\
		&=\langle f^*\OO_{\PP^2},f^*\OO_{\PP^2}(1),f^*\OO_{\PP^2}(2),f^*\OO_{\PP^2}\otimes \OO_X(1),f^*\OO_{\PP^2}(1)\otimes \OO_X(1),f^*\OO_{\PP^2}(2)\otimes \OO_X(1)\rangle\nonumber\\
		&=\langle \OO_X,\OO_X(H),\OO_X(2H),\OO_X(aH+D),\OO_X((a+1)H+D),\OO_X((a+2)H+D)\rangle.\nonumber
		\end{align}
		Hence type $(1)$ is a full exceptional collection and thus, all the exceptional collections of length $6$ are full exceptional collections.
	\end{proof}
\end{thm}

\subsection{When $X \cong \PP^2 \times \PP^1$.}
Let $p_1: \PP^2 \times \PP^1\rightarrow \PP^2$ and $p_2: \PP^2 \times \PP^1\rightarrow \PP^1$ be the projections onto $\PP^2$ and $\PP^1$, respectively. Then the canonical divisor of $X$ is $K_X=-3H-2D$, where $H$ and $D$ be the pullback of the hyperplane class of $\PP^2$ and $\PP^1$, respectively.  

\subsubsection{Classification of full exceptional collections}
\begin{thm} \label{cl2}
	Let $X$ be $\PP^2 \times \PP^1$. Then the normalized sequence
	\begin{align}
	\{\OO_X,\OO_X(D_1),\cdots,\OO_X(D_5)\} \label{2} \tag{B}
	\end{align}
	is an exceptional collection of line bundles if and only if the ordered set of divisors $\{D_1,D_2,\cdots,D_5\}$ is one of the following types:
	\begin{itemize}
		\item [(1)] $\{aH+D,(a+1)H+D,(a+2)H+D,H+2D,2H+2D\}$,
		\item [(2)]
		$\{H,2H,aH+D,(a+1)H+D,(a+2)H+D\}$,
		\item [(3)] $\{H,aH+D,(a+1)H+D,(a+2)H+D,2H+2D\}$,
		\item [(4)] $\{H+bD, H+(b+1)D, 2H+cD, 2H+(c+1)D, 3H+D\}$,
		\item [(5)] $\{D, H+bD, H+(b+1)D, 2H+cD, 2H+(c+1)D\}$,
		\item [(6)] $\{D, H+bD, 2H+bD, H+(b+1)D, 2H+(b+1)D\}$,
		\item [(7)] $\{H+bD, 2H+bD, H+(b+1)D, 2H+(b+1)D,3H+D\}$,
		\item [(8)] $\{H, D, H+D, 2H+bD, 2H+(b+1)D\}$,
		\item [(9)] $\{-H+D, D, H+bD, H+(b+1)D, 2H+2D\}$,
		\item [(10)] $\{H, 2H+bD, 2H+(b+1)D, 3H+D, 4H+D\}$,
		\item [(11)] $\{H+bD, H+(b+1)D, 2H+D, 3H+D, 2H+2D\}$,
		\item [(12)] $\{H+D, 2H+D, H+2D, 3H+D, 2H+2D\}$,
		\item [(13)] $\{H, D, 2H, H+D, 2H+D\}$,
		\item [(14)] $\{-H+D, H, D, H+D, 2H+2D\}$,
		\item [(15)] $\{2H-D, H, 2H, 3H+D, 4H+D\}$,
		\item [(16)] $\{-H+D, D, H+2D, 2H+2D, H+3D\}$,
		\item [(17)] $\{H, 2H+D, 3H+D, 2H+2D, 4H+D\}$,
	\end{itemize}
	where $a,b,c\in \mathbb{Z}$. Moreover, by mutations and normalizations, they are related as:
	\begin{align}
	&(1)\rightarrow (2)\rightarrow (3)\rightarrow (1), \nonumber\\
	&(4)\rightarrow (5)\rightarrow (4),\nonumber\\
	&(6)\rightarrow (7)\rightarrow (8)\rightarrow (9)\rightarrow (10)\rightarrow (11)\rightarrow (6), \nonumber\\
	&(12)\rightarrow (13)\rightarrow (14)\rightarrow (15)\rightarrow (16)\rightarrow (17)\rightarrow (12). \nonumber
	\end{align}
	\begin{proof}
		Note that $\OO_X(aH+bD)$ is cohomologically zero line bundle if $a=-1, 2$ or $b=-1$. Based on Table \ref{tab:table2} and the rules below, we get the result.
		\begin{itemize}
			\item [(1)] $B_0$ can occur at most three times.
			\item [(2)] $B_1$ and $B_2$ can occur at most two times.
			\item [(3)] $B_1$ cannot occur after two $B_2$.
		\end{itemize} 
	\begin{table}[h!]
		\begin{center}
			\caption{Exceptional pairs.}
			\label{tab:table2}
			\begin{tabular}{c|c|c|c} 
				& $B_0'$ & $B_1'$ & $B_2'$ \\
				\hline
				$B_0$ &	$a'=a+1, a+2$ & $a=-1,0$ or $b'=2$  & $a=2$ or $b'=-1,0,1$  \\
				\hline
				$B_1$ &	$a'=2,3$ or $b=0$ & $b'=b+1$ & \checkmark \\
				\hline
				$B_2$ &	$a'=3,4$ or $c=0$ & $b'=c+1$ & $c'=c+1$ \\
				\hline
			\end{tabular}
			\vspace{3mm}
			
			where $B_0=aH+D, B_1=H+bD, B_2=2H+cD$.
		\end{center}
	\end{table}
	\end{proof}
\end{thm}

\begin{thm}
	Let $X$ be $\PP^2 \times \PP^1$. Then any exceptional collection of line bundles of length 6 on $X$ is full.
	\begin{proof}
		By Lemma \ref{nor}, we only need to show that any exceptional collection of line bundles of length $6$ in Theorem \ref{cl2} is full. By Lemma \ref{mt}, it suffices to show that the exceptional collections of type $(2), (5), (6)$ and $(12)$ in Theorem \ref{cl2} are full. 
		
		We first show that exceptional collection of type $(2)$ is full. Consider the projective bundle $f:X\cong \PP_{\PP^2}(\OO(-a)\oplus \OO(-a))\rightarrow \PP^2$, where $\OO_X(1)=aH+D$. Then by the projective bundle formula, we have the following semiorthogonal decomposition of $\D(X)$.
		\begin{align}
		\D(X)&=\langle f^*\D(\PP^2),f^*\D(\PP^2)\otimes \OO_X(1)\rangle \nonumber\\
		&=\langle f^*\OO_{\PP^2},f^*\OO_{\PP^2}(1),f^*\OO_{\PP^2}(2),f^*\OO_{\PP^3}\otimes \OO_X(1),f^*\OO_{\PP^2}(1)\otimes\OO_X(2),f^*\OO_{\PP^2}(2)\otimes \OO_X(1)\rangle\nonumber\\
		&=\langle \OO_X,\OO_X(H),\OO_X(2H),\OO_X(aH+D),\OO_X((a+1)H+D),\OO_X((a+2)H+D)\rangle.\nonumber
		\end{align}
		Hence type $(2)$ is a full exceptional collection.
		
		Now we show that exceptional collection of type $(5)$ is full. Consider the projective bundle $f:X\cong \PP_{\PP^1}(\OO(-b)\oplus \OO(-b)\oplus \OO(-b))\rightarrow \PP^1$, where $\OO_X(1)=H+bD$. Then by the projective bundle formula, we obtain the following semiorthogonal decomposition of $\D(X)$.
		\begin{align}
		\D(X)&=\langle f^*\D(\PP^1),f^*\D(\PP^1)\otimes \OO_X(1),f^*\D(\PP^1)\otimes \OO_X(2)\rangle \nonumber\\
		&=\langle f^*\OO_{\PP^1},f^*\OO_{\PP^1}(1),f^*\OO_{\PP^1}\otimes \OO_X(1),f^*\OO_{\PP^1}(1)\otimes\OO_X(1),f^*\OO_{\PP^1}\otimes \OO_X(2),f^*\OO_{\PP^1}(1)\otimes \OO_X(2)\rangle\nonumber\\
		&=\langle\OO_X,\OO_X(D),\OO_X(H+bD),\OO_X(H+(b+1)D),\OO_X(2H+2bD),\OO_X(2H+(2b+1)D)\rangle.\nonumber
		\end{align}
		Based on this semiorthogonal decomposition, we show that the exceptional collection of type $(5)$ is full by mathematical induction. For a given $b\in \mathbb{Z}$, assume that for $c=k$ the exceptional collection
		\begin{align}\label{mi1} \tag{E1}
		\{D, H+bD, H+(b+1)D, 2H+kD, 2H+(k+1)D\}
		\end{align}
		is full. Then for $c=k-1$, we have the exceptional collection
		\begin{align}\label{mi2} \tag{E2}
		\{D, H+bD, H+(b+1)D, 2H+(k-1)D, 2H+kD\}
		\end{align}
		and for $c=k+1$, we have the exceptional collection
		\begin{align}\label{mi3} \tag{E3}
		\{D, H+bD, H+(b+1)D, 2H+(k+1)D, 2H+(k+2)D\}.
		\end{align}
		Comparing the exceptional collections (\ref{mi2}) and (\ref{mi3}) with (\ref{mi1}) and by Lemma \ref{mt1}, the exceptional collections (\ref{mi2}) and (\ref{mi3}) are full. Hence the exceptional collection of type $(5)$ is also full.
		
		In order to show that exceptional collection of type $(6)$ is full, we set $c=b$ in the exceptional collection of type $(5)$ to get 
		\begin{align}
		\{D, H+bD, H+(b+1)D, 2H+bD, 2H+(b+1)D\}.\nonumber
		\end{align}
		Then by Lemma \ref{mt}, we see that 
		\begin{align}
		\{D, H+bD, 2H+bD, H+(b+1)D, 2H+(b+1)D\}\nonumber
		\end{align}
		is also a full exceptional collection.
		
		Setting $b=1$ in exceptional collection of type $(7)$ gives a full exceptional collection 
		\begin{align}
		\{H+D,2H+D,H+2D,2H+2D,3H+D\}.\nonumber
		\end{align}
		Mutating $2H+2D$ and $3H+D$ gives a full exceptional collection
		\begin{align}
		\{H+D,2H+D,H+2D,3H+D,2H+2D\},\nonumber
		\end{align}
		which is of type $(12)$. 
		
		Therefore, all the exceptional collections of length $6$ are full.
	\end{proof}
\end{thm}	

As a result, the conjecture \ref{Kuzc} holds true for smooth toric Fano threefolds of Picard rank two.

\section{Smooth toric Fano fourfolds with Picard rank 2}
\subsection{When $X=\PP_{\PP^3}(\OO\oplus\OO(3))$.}	
Let $f:\PP_{\PP^3}(\OO\oplus \OO(3))\rightarrow \PP^3$ be the projective bundle. Then its canonical divisor is $K_X=-7H-2D$ where $H$ is the pullback of hyperplane class in $\PP^3$ and $D$ is the class of $c_1(\OO_X(1))$. The Picard group of $X$ is 
\begin{align}
\Pic(X)\cong \Pic(\PP^3)\oplus \mathbb{Z}[D]=\mathbb{Z}[H]\oplus\mathbb{Z}[D], \nonumber
\end{align}
with intersection numbers
\begin{align}
H^4=0, ~H^3D=1, ~H^2D^2=-3, ~HD^3=9, ~D^4=-27. \nonumber
\end{align}
\subsubsection{Cohomologically zero line bundles}
\begin{lem}\label{CZ3}
	$H^0(X,aH+bD)=0$ if and only if $b<0$ or $a+3b<0, ~b\geq 0$. Also, $H^4(X,aH+bD)=0$ if and only if $b>-2$ or $a+3b>-13, ~b\leq -2$.
	\begin{proof}
		Note that we have the following isomorphisms:
		\begin{align}
		f_*(\OO_X(b))&\cong \Sym^b(\OO_{\PP^3}\oplus\OO_{\PP^3}(3)) \nonumber \\ 
		&\cong \bigoplus_{i=0}^b \Sym^i(\OO_{\PP^3})\otimes \Sym^{b-i}(\OO_{\PP^3}(3)) \nonumber \\
		&\cong \bigoplus_{i=0}^b \OO_{\PP^3}(3b-3i) ~\text{if}~ b\geq 0, \nonumber
		\end{align} and $f_*(\OO_X(b))\cong 0$ if $b<0$.
		
		By projection formula, we have the following isomorphisms:
		\begin{align}
		H^0(X, f^*(\OO_{\PP^3}(a))\otimes \OO_X(b))&\cong H^0(\PP^3, \OO_{\PP^3}(a)\otimes  \bigoplus_{i=0}^b \OO_{\PP^3}(3b-3i)) \nonumber \\
		&\cong \bigoplus_{i=0}^b H^0(\PP^3,\OO_{\PP^3}(a+3b-3i)) ~\text{if}~ b\geq 0,\nonumber
		\end{align}
		and if $b<0$, then $H^0(X,aH+bD)=0$.
		Hence, $H^0(X,aH+bD)=0$ if and only if $b<0$ or $a+3b<0, ~b\geq 0$.
		
		By Serre duality, one obtains the condition for $H^4(X,aH+bD)=0$.
    \end{proof}
\end{lem}
 \begin{pro} \label{pro3}
 	A line bundle $\OO_X(aH+bD)$ is cohomologically zero if one of the following holds:
 	\begin{itemize}
 		\item [(1)] $b=-1$,
 		\item [(2)] $a=-1, b=0$,
 		\item [(3)] $a=-2, b=0$,
 		\item [(4)] $a=-3, b=0$,
 		\item [(5)] $a=-4, b=-2$,
 		\item [(6)] $a=-5, b=-2$,
 		\item [(7)] $a=-6, b=-2$.
 	\end{itemize}
 \begin{proof}
 	Since $c(X)=(1+4H+6H^2+4H^3)\{(1+D)^2+2H(1+D)\}$, we have $c_1=7H+2D$, $c_2=18H^2+6HD$, $c_3=22H^3+12H^2D$ and $c_4=12H^4+8H^3D=8$.
 	Thus, by Lemma \ref{HRR}, we have
 	\begin{align}
 	\chi(\OO_X(aH+bD))&=\frac{1}{24}(24+44a-6b+24a^2-10ab+15b^2+4a^3+6a^2b-18ab^2 \nonumber\\
 	&\qquad +18b^3+4a^3b-18a^2b^2+36ab^3-27b^4) \nonumber\\
 	&=\frac{1}{24}(b+1)(2a-3b+4)(2a^2-6ab+8a+9b^2-3b+6). \nonumber
 	\end{align}
 	Therefore, by Lemma \ref{CZ3}, $\OO_X(aH+bD)$ is cohomologically zero if the following conditions hold:
 	\begin{itemize}
 		\item [$\mathrm{(i)}$] $b<0$ or $a+3b<0, ~b\geq 0$,
 		\item [$\mathrm{(ii)}$] $b>-2$ or $a+3b>-13, ~b\leq -2$,
 		\item [$\mathrm{(iii)}$] $(b+1)(2a-3b+4)(2a^2-6ab+8a+9b^2-3b+6)=0$.
 	\end{itemize}
 Note that the equation $(2a^2-6ab+8a+9b^2-3b+6)=0$ has $4$ integer solutions; $(a,b)=(-6,-2), (-4,-2), (-3,0), (-1,0)$. Also, the equation $2a-3b+4=0$ has only two integer solutions $(a,b)=(-5,-2), (-2,0)$ if it is constrained on the conditions $\mathrm{(i)}$ and $\mathrm{(ii)}$. This completes the proof.
 \end{proof}
 \end{pro}

 \subsubsection{Classification of full exceptional collections}
 \begin{thm} \label{cl3}
 	Let $X$ be the projective bundle $f:X=\PP_{\PP^3}(\OO\oplus \OO(3))\rightarrow \PP^3$. Then the normalized sequence
 	\begin{align}
 	\{\OO_X,\OO_X(D_1),\cdots,\OO_X(D_7)\} \label{3} \tag{C}
 	\end{align}
 	is an exceptional collection of line bundles if and only if the ordered set of divisors $\{D_1,D_2,\cdots,D_7\}$ is one of the following types:
 	\begin{itemize}
 		\item [(1)] $\{H,2H,3H,aH+D,(a+1)H+D,(a+2)H+D,(a+3)H+D\}$,
 		\item [(2)] $\{H,2H,aH+D,(a+1)H+D,(a+2)H+D,(a+3)H+D,6H+2D\}$,
 		\item [(3)] $\{H,aH+D,(a+1)H+D,(a+2)H+D,(a+3)H+D,,6H+2D\}$,
 		\item [(4)] $\{aH+D,(a+1)H+D,(a+2)H+D,(a+3)H+D,4H+2D,5H+2D,6H+2D\}$,
 	\end{itemize}
 where $a\in \mathbb{Z}$. Moreover, by mutations and normalizations, they are related as:
 \begin{align}
 (1)\rightarrow (2)\rightarrow (3)\rightarrow (4). \nonumber
 \end{align}
 \begin{proof}
 	Write $D_0=0$. The sequence (\ref{3}) is an exceptional collection if and only if for any integers $0\leq i<j\leq 7$, the line bundles $\OO_X(D_i-D_j)$ are cohomologically zero.
 	
 	To make the sequence (\ref{3}) an exceptional collection, $\OO_X(D_i)$ must be one of the line bundles in Proposition \ref{pro3}. To find out all the possible exceptional collections, we need Table \ref{tab:table3}.
 	
 	From Table \ref{tab:table3}, we observe that $\OO_X(D_1)$ must be one of $B_0$ or $B_1$. To finish the proof, we shall discuss case by case:
 	\begin{itemize}
 		\item [Case1.] Suppose $\OO_X(D_1)=B_0$. In this case, the only possible combination is $\{aH+D,(a+1)H+D,(a+2)H+D,(a+3)H+D,4H+2D,5H+2D,6H+2D\}$, which is $(4)$ in Theorem \ref{cl3}.
 		\item [Case2.] Suppose $\OO_X(D_1)=B_1$. In this case, there are two possibilities for $D_2$:
 		\begin{itemize}
 			\item [(a)] When $\OO_X(D_2)=B_0$. In this case, we have only one type of exceptional collections:
 			$\{H,aH+D,(a+1)H+D,(a+2)H+D,(a+3)H+D,6H+2D,5H+2D\}$, which is $(3)$ in Theorem \ref{cl3}
 			\item [(b)] When $\OO_X(D_2)=B_2$. In this case, we have two types of exceptional collections:
 			$\{H,2H,3H,aH+D,(a+1)H+D,(a+2)H+D,(a+3)H+D\}$ and $\{H,2H,aH+D,(a+1)H+D,(a+2)H+D,(a+3)H+D,6H+2D\}$, which are $(1)$ and $(2)$ in Theorem \ref{cl3}, respectively.
 		\end{itemize}
 	\end{itemize}
  	\begin{table}[h!] 
 	\begin{center}
 		\caption{Exceptional pairs.}
 		\label{tab:table3}
 		\begin{tabular}{c|c|c|c|c|c|c|c}
 			& $B_0'$ & $B_1$ & $B_2$ & $B_3$ & $B_4$ & $B_5$ & $B_6$\\
 			\hline
 			$B_0$ &	$a'=a+1, a+2, a+3$ &  &  &  & \checkmark  & \checkmark & \checkmark \\
 			\hline
 			$B_1$ &	\checkmark &  & \checkmark & \checkmark &  & \checkmark & \checkmark \\
 			\hline
 			$B_2$ &	\checkmark &  &  & \checkmark &   &  & \checkmark \\
 			\hline
 			$B_3$ & \checkmark &  &  &  &   &  &  \\
 			\hline
 			$B_4$ &  &  &  &  &  & \checkmark & \checkmark \\
 			\hline
 			$B_5$ &  &  &  &  &  &  & \checkmark \\
 			\hline
 			$B_6$ &  &  &  &  &  &  &  \\
 			\hline
 		\end{tabular}
 		\vspace{3mm}
 		
 		where $B_0=aH+D, B_1=H, B_2=2H, B_3=3H,$\\ 
 		$B_4=4H+2D, B_5=5H+2D, B_6=6H+2D$.
 	\end{center}
 \end{table}
 \end{proof}
 \end{thm}

 \begin{thm}
 	Let $X$ be the projective bundle $f:X=\PP_{\PP^3}(\OO\oplus \OO(3))\rightarrow \PP^3$. Then any exceptional collection of line bundles of length 8 on $X$ is full.
 	\begin{proof}
 		By Lemma \ref{nor}, we only need to show that any exceptional collection of line bundles of length $6$ in Theorem \ref{cl3} is full. By Lemma \ref{mt}, it suffices to show that the exceptional collection of type $(1)$ in Theorem \ref{cl3} is full. 
 		
 		By Beilinson's semiorthogonal decomposition, we have 
 		\begin{align}
 		\D(\PP^3)=\langle\OO_{\PP^3}, \OO_{\PP^3}(1),\OO_{\PP^3}(2),\OO_{\PP^3}(3)\rangle. \nonumber
 		\end{align}
 		By Orlov's projective bundle formula (Theorem \ref{proj}), we obtain the following semiorthogonal decomposition of $\D(X)$
 		\begin{align}
 		\D(X)&=\langle f^*\D(\PP^3),f^*\D(\PP^3)\otimes \OO_X(1)\rangle \nonumber\\
 		&=\langle f^*\OO_{\PP^3},f^*\OO_{\PP^3}(1),f^*\OO_{\PP^3}(2),f^*\OO_{\PP^3}(3),f^*\OO_{\PP^3}\otimes \OO_X(1),\nonumber\\
 		&\qquad f^*\OO_{\PP^3}(1)\otimes\OO_X(1),f^*\OO_{\PP^3}(2)\otimes \OO_X(1),f^*\OO_{\PP^3}(3)\otimes \OO_X(1)\rangle\nonumber\\
 		&=\langle\OO_X, \OO_X(H),\OO_X(2H),\OO_X(3H),\OO_X(aH+D),\nonumber\\
 		&\qquad \OO_X((a+1)H+D),\OO_X((a+2)H+D),\OO_X((a+3)H+D)\rangle.\nonumber
  		\end{align}
 		Hence type $(1)$ is a full exceptional collection and so all the exceptional collections of length $8$ on $X$ are full exceptional collections.
 	\end{proof}
 \end{thm}	

\subsection{When $X=\PP_{\PP^3}(\OO\oplus\OO(2))$}	
Let $f:\PP_{\PP^3}(\OO\oplus \OO(2))\rightarrow \PP^3$ be the projective bundle. Then its canonical divisor is $K_X=-6H-2D$ where $H$ is the pullback of hyperplane class in $\PP^3$ and $D$ is the class of $c_1(\OO_X(1))$. The Picard group of $X$ is 
\begin{align}
\Pic(X)\cong \Pic(\PP^3)\oplus \mathbb{Z}[D]=\mathbb{Z}[H]\oplus\mathbb{Z}[D], \nonumber
\end{align}
with intersection numbers
\begin{align}
H^4=0, ~H^3D=1, ~H^2D^2=-2, ~HD^3=4, ~D^4=-8. \nonumber
\end{align}
\subsubsection{Cohomologically zero line bundles}
\begin{lem}\label{CZ4}
	$H^0(X,aH+bD)=0$ if and only if $b<0$ or $a+2b<0, ~b\geq 0$. Also, $H^4(X,aH+bD)=0$ if and only if $b>-2$ or $a+2b>-10, ~b\leq -2$.
	\begin{proof}
		Note that we have the following isomorphisms:
		\begin{align}
		f_*(\OO_X(b))&\cong \Sym^b(\OO_{\PP^3}\oplus\OO_{\PP^3}(2)) \nonumber \\ 
		&\cong \bigoplus_{i=0}^b \Sym^i(\OO_{\PP^3})\otimes \Sym^{b-i}(\OO_{\PP^3}(2)) \nonumber \\
		&\cong \bigoplus_{i=0}^b \OO_{\PP^3}(2b-2i) ~\text{if}~ b\geq 0, \nonumber
		\end{align} and $f_*(\OO_X(b))\cong 0$ if $b<0$.
		
		By projection formula, we have the following isomorphisms:
		\begin{align}
		H^0(X, f^*(\OO_{\PP^3}(a))\otimes \OO_X(b))&\cong H^0(\PP^3, \OO_{\PP^3}(a)\otimes  \bigoplus_{i=0}^b \OO_{\PP^3}(2b-2i)) \nonumber \\
		&\cong \bigoplus_{i=0}^b H^0(\PP^3,\OO_{\PP^3}(a+2b-2i)) ~\text{if}~ b\geq 0,\nonumber
		\end{align}
		and if $b<0$, then $H^0(X,aH+bD)=0$.
		Hence, $H^0(X,aH+bD)=0$ if and only if $b<0$ or $a+2b<0, ~b\geq 0$.
		
		By Serre duality, one obtains the condition for $H^4(X,aH+bD)=0$.
	\end{proof}
\end{lem}
\begin{pro} \label{pro4}
	A line bundle $\OO_X(aH+bD)$ is cohomologically zero if one of the following holds:
	\begin{itemize}
		\item [(1)] $b=-1$,
		\item [(2)] $a=-1, b=0$,
		\item [(3)] $a=-2, b=0$,
		\item [(4)] $a=-3, b=0$,
		\item [(5)] $a=-3, b=-2$,
		\item [(6)] $a=-4, b=-2$,
		\item [(7)] $a=-5, b=-2$.
	\end{itemize}
	\begin{proof}
		Since $c(X)=(1+4H+6H^2+4H^3)\{(1+D)^2+2H(1+D)\}$, we have $c_1=6H+2D$, $c_2=14H^2+8HD$, $c_3=16H^3+12H^2D$ and $c_4=8H^4+8H^3D=8$.
		Thus, by Lemma \ref{HRR}, we have
		\begin{align}
		\chi(\OO_X(aH+bD))&=\frac{1}{6}(6+11a-b+6a^2+ab-b^2 \nonumber\\
		&\qquad +a^3+3a^2b-6ab^2+4b^3+a^3b-3a^2b^2+4ab^3-2b^4) \nonumber\\
		&=\frac{1}{6}(b+1)(a-b+2)(a^2-2ab+4a+2b^2-2b+3). \nonumber
		\end{align}
		Therefore, by Lemma \ref{CZ4}, $\OO_X(aH+bD)$ is cohomologically zero if the following conditions hold:
		\begin{itemize}
			\item [$\mathrm{(i)}$] $b<0$ or $a+2b<0, ~b\geq 0$,
			\item [$\mathrm{(ii)}$] $b>-2$ or $a+2b>-10, ~b\leq -2$,
			\item [$\mathrm{(iii)}$] $(b+1)(a-b+2)(a^2-2ab+4a+2b^2-2b+3)=0$.
		\end{itemize}
		Note that the equation $(a^2-2ab+4a+2b^2-2b+3)=0$ has $4$ integer solutions; $(a,b)=(-5,-2), (-3,-2), (-3,0), (-1,0)$. Also, the equation $2a-3b+4=0$ has only two integer solutions $(a,b)=(-4,-2), (-2,0)$ if it is constrained on the conditions $\mathrm{(i)}$ and $\mathrm{(ii)}$. This completes the proof.
	\end{proof}
\end{pro}

\subsubsection{Classification of full exceptional collections}
\begin{thm} \label{cl4}
	Let $X$ be the projective bundle $f:X=\PP_{\PP^3}(\OO\oplus \OO(2))\rightarrow \PP^3$. Then the normalized sequence
	\begin{align}
	\{\OO_X,\OO_X(D_1),\cdots,\OO_X(D_7)\} \label{4} \tag{D}
	\end{align}
	is an exceptional collection of line bundles if and only if the ordered set of divisors $\{D_1,D_2,\cdots,D_7\}$ is one of the following types:
	\begin{itemize}
		\item [(1)] $\{H,2H,3H,aH+D,(a+1)H+D,(a+2)H+D,(a+3)H+D\}$,
		\item [(2)] $\{H,2H,aH+D,(a+1)H+D,(a+2)H+D,(a+3)H+D,5H+2D\}$,
		\item [(3)] $\{H,aH+D,(a+1)H+D,(a+2)H+D,(a+3)H+D,4H+2D,5H+2D\}$,
		\item [(4)] $\{aH+D,(a+1)H+D,(a+2)H+D,(a+3)H+D,3H+2D,4H+2D,5H+2D\}$,
	\end{itemize}
	where $a\in \mathbb{Z}$. Moreover, by mutations and normalizations, they are related as:
	\begin{align}
	(1)\rightarrow (2)\rightarrow (3)\rightarrow (4). \nonumber
	\end{align}
	\begin{proof}
		Write $D_0=0$. The sequence (\ref{4}) is an exceptional collection if and only if for any integers $0\leq i<j\leq 7$, the line bundles $\OO_X(D_i-D_j)$ are cohomologically zero.
		
		To make the sequence (\ref{4}) an exceptional collection, $\OO_X(D_i)$ must be one of the line bundles in Proposition \ref{pro4}. To find out all the possible exceptional collections, we need Table \ref{tab:table4}.
		
		From Table \ref{tab:table4}, we observe that $\OO_X(D_1)$ must be one of $B_0$ or $B_1$. To finish the proof, we shall discuss case by case:
		\begin{itemize}
			\item [Case1.] Suppose $\OO_X(D_1)=B_0$. In this case, the only possible combination is $\{aH+D,(a+1)H+D,(a+2)H+D,(a+3)H+D,3H+2D,4H+2D,5H+2D\}$, which is $(4)$ in Theorem \ref{cl4}.
			\item [Case2.] Suppose $\OO_X(D_1)=B_1$. In this case, there are two possibilities for $D_2$:
			\begin{itemize}
				\item [(a)] When $\OO_X(D_2)=B_0$. In this case, we have only one type of exceptional collections:
				$\{H,aH+D,(a+1)H+D,(a+2)H+D,(a+3)H+D,4H+2D,5H+2D\}$, which is $(3)$ in Theorem \ref{cl4}
				\item [(b)] When $\OO_X(D_2)=B_2$. In this case, we have two types of exceptional collections:
				$\{H,2H,3H,aH+D,(a+1)H+D,(a+2)H+D,(a+3)H+D\}$ and $\{H,2H,aH+D,(a+1)H+D,(a+2)H+D,(a+3)H+D,5H+2D\}$, which are $(1)$ and $(2)$ in Theorem \ref{cl4}, respectively.
			\end{itemize}
		\end{itemize}
\begin{table}[h!]
	\begin{center}
		\caption{Exceptional pairs.}
		\label{tab:table4}
		\begin{tabular}{c|c|c|c|c|c|c|c} 
			& $B_0'$ & $B_1$ & $B_2$ & $B_3$ & $B_4$ & $B_5$ & $B_6$\\
			\hline
			$B_0$ &	$a'=a+1, a+2, a+3$ &  &  &  & \checkmark  & \checkmark & \checkmark \\
			\hline
			$B_1$ &	\checkmark &  & \checkmark & \checkmark &  & \checkmark & \checkmark \\
			\hline
			$B_2$ &	\checkmark &  &  & \checkmark &   &  & \checkmark \\
			\hline
			$B_3$ & \checkmark &  &  &  &   &  &  \\
			\hline
			$B_4$ &  &  &  &  &  & \checkmark & \checkmark \\
			\hline
			$B_5$ &  &  &  &  &  &  & \checkmark \\
			\hline
			$B_6$ &  &  &  &  &  &  &  \\
			\hline
		\end{tabular}
	\vspace{3mm}
	
	where $B_0=aH+D, B_1=H, B_2=2H, B_3=3H,$\\ 
	$B_4=3H+2D, B_5=4H+2D, B_6=5H+2D$.
	\end{center}
\end{table}
\end{proof}
\end{thm}

\begin{thm}
	Let $X$ be the projective bundle $f:X=\PP_{\PP^3}(\OO\oplus \OO(2))\rightarrow \PP^3$. Then any exceptional collection of line bundles of length 8 on $X$ is full.
	\begin{proof}
		By Lemma \ref{nor}, we only need to show that any exceptional collection of line bundles of length $6$ in Theorem \ref{cl4} is full. By Lemma \ref{mt}, it suffices to show that the exceptional collection of type $(1)$ in Theorem \ref{cl4} is full. 
		
		By Beilinson's semiorthogonal decomposition, we have 
		\begin{align}
		\D(\PP^3)=\langle\OO_{\PP^3}, \OO_{\PP^3}(1),\OO_{\PP^3}(2),\OO_{\PP^3}(3)\rangle. \nonumber
		\end{align}
		By Theorem \ref{proj}, we obtain the following semiorthogonal decomposition of $\D(X)$
		\begin{align}
		\D(X)&=\langle f^*\D(\PP^3),f^*\D(\PP^3)\otimes \OO_X(1)\rangle \nonumber\\
		&=\langle f^*\OO_{\PP^3},f^*\OO_{\PP^3}(1),f^*\OO_{\PP^3}(2),f^*\OO_{\PP^3}(3),f^*\OO_{\PP^3}\otimes \OO_X(1),\nonumber\\
        &\qquad f^*\OO_{\PP^3}(1)\otimes\OO_X(1),f^*\OO_{\PP^3}(2)\otimes \OO_X(1),f^*\OO_{\PP^3}(3)\otimes \OO_X(1)\rangle\nonumber\\
		&=\langle\OO_X,\OO_X(H),\OO_X(2H),\OO_X(3H),\nonumber\\
		&\qquad \OO_X(aH+D),\OO_X((a+1)H+D),\OO_X((a+2)H+D),\OO_X((a+3)H+D)\rangle.\nonumber
		\end{align}
		Hence type $(1)$ is a full exceptional collection and so all the exceptional collections of length $8$ on $X$ are full exceptional collections.
	\end{proof}
\end{thm}	

\subsection{When $X=\PP_{\PP^3}(\OO\oplus\OO(1))$}	
Let $f:\PP_{\PP^3}(\OO\oplus \OO(1))\rightarrow \PP^3$ be the projective bundle. Then its canonical divisor is $K_X=-5H-2D$ where $H$ is the pullback of hyperplane class in $\PP^3$ and $D$ is the class of $c_1(\OO_X(1))$. The Picard group of $X$ is 
\begin{align}
\Pic(X)\cong \Pic(\PP^3)\oplus \mathbb{Z}[D]=\mathbb{Z}[H]\oplus\mathbb{Z}[D], \nonumber
\end{align}
with intersection numbers
\begin{align}
H^4=0, ~H^3D=1, ~H^2D^2=-1, ~HD^3=1, ~D^4=-1. \nonumber
\end{align}
\subsubsection{Cohomologically zero line bundles}
\begin{lem}\label{CZ5}
	$H^0(X,aH+bD)=0$ if and only if $b<0$ or $a+b<0, ~b\geq 0$. Also, $H^4(X,aH+bD)=0$ if and only if $b>-2$ or $a+b>-7, ~b\leq -2$.
	\begin{proof}
		Note that we have the following isomorphisms:
		\begin{align}
		f_*(\OO_X(b))&\cong \Sym^b(\OO_{\PP^3}\oplus\OO_{\PP^3}(1)) \nonumber \\ 
		&\cong \bigoplus_{i=0}^b \Sym^i(\OO_{\PP^3})\otimes \Sym^{b-i}(\OO_{\PP^3}(1)) \nonumber \\
		&\cong \bigoplus_{i=0}^b \OO_{\PP^3}(b-i) ~\text{if}~ b\geq 0, \nonumber
		\end{align} and $f_*(\OO_X(b))\cong 0$ if $b<0$.
		
		By projection formula, we have the following isomorphisms:
		\begin{align}
		H^0(X, f^*(\OO_{\PP^3}(a))\otimes \OO_X(b))&\cong H^0(\PP^3, \OO_{\PP^3}(a)\otimes  \bigoplus_{i=0}^b \OO_{\PP^3}(b-i)) \nonumber \\
		&\cong \bigoplus_{i=0}^b H^0(\PP^3,\OO_{\PP^3}(a+b-i)) ~\text{if}~ b\geq 0,\nonumber
		\end{align}
		and if $b<0$, then $H^0(X,aH+bD)=0$.
		Hence, $H^0(X,aH+bD)=0$ if and only if $b<0$ or $a+b<0, ~b\geq 0$.
		
		By Serre duality, one obtains the condition for $H^4(X,aH+bD)=0$.
	\end{proof}
\end{lem}
\begin{pro} \label{pro5}
	A line bundle $\OO_X(aH+bD)$ is cohomologically zero if one of the following holds:
	\begin{itemize}
		\item [(1)] $b=-1$,
		\item [(2)] $a=-1, b=0$,
		\item [(3)] $a=-2, b=0$,
		\item [(4)] $a=-3, b=0$,
		\item [(5)] $a=-2, b=-2$,
		\item [(6)] $a=-3, b=-2$,
		\item [(7)] $a=-4, b=-2$,
		\item [(8)] $a=-2, b=1$,
		\item [(9)] $a=-3, b=-3$.
	\end{itemize}
	\begin{proof}
		Since $c(X)=(1+4H+6H^2+4H^3)\{(1+D)^2+H(1+D)\}$, we have $c_1=5H+2D$, $c_2=10H^2+8HD$, $c_3=10H^3+12H^2D$ and $c_4=4H^4+8H^3D=8$.
		Thus, by Lemma \ref{HRR}, we have
		\begin{align}
		\chi(\OO_X(aH+bD))&=\frac{1}{24}(24+44a+6b+24a^2+22ab-11b^2 \nonumber\\
		&\qquad +4a^3+18a^2b-18ab^2+6b^3+4a^3b-6a^2b^2+4ab^3-b^4) \nonumber\\
		&=\frac{1}{24}(b+1)(2a-b+4)(2a^2-2ab+8a+b^2-2b+6). \nonumber
		\end{align}
		Therefore, by Lemma \ref{CZ5}, $\OO_X(aH+bD)$ is cohomologically zero if the following conditions hold:
		\begin{itemize}
			\item [($\mathrm{i}$)] $b<0$ or $a+b<0, ~b\geq 0$,
			\item [($\mathrm{ii}$)] $b>-2$ or $a+b>-7, ~b\leq -2$,
			\item [($\mathrm{iii}$)] $(b+1)(2a-b+4)(2a^2-2ab+8a+b^2-2b+6)=0$.
		\end{itemize}
		Note that the equation $2a^2-2ab+8a+b^2-2b+6=0$ has $6$ integer solutions; $(a,b)=(-4,-2), (-2,-2),$ $(-3,0), (-1,0), (-3,-3), (-2,1)$ under the conditions $\mathrm{(i)}$ and $\mathrm{(ii)}$. Also, the equation $2a-b+4=0$ has only two integer solutions $(a,b)=(-3,-2), (-2,0)$ if it is constrained on the conditions $\mathrm{(i)}$ and $\mathrm{(ii)}$. This completes the proof.
	\end{proof}
\end{pro}

\subsubsection{Classification of full exceptional collections}
\begin{thm} \label{cl5}
	Let $X$ be the projective bundle $f:X=\PP_{\PP^3}(\OO\oplus \OO(1))\rightarrow \PP^3$. Then the normalized sequence
	\begin{align}
	\{\OO_X,\OO_X(D_1),\cdots,\OO_X(D_7)\} \label{5} \tag{E}
	\end{align}
	is an exceptional collection of line bundles if and only if the ordered set of divisors $\{D_1,D_2,\cdots,D_7\}$ is one of the following types:
	\begin{itemize}
		\item [(1)] $\{H,2H,3H,aH+D,(a+1)H+D,(a+2)H+D,(a+3)H+D\}$,
		\item [(2)] $\{H,2H,aH+D,(a+1)H+D,(a+2)H+D,(a+3)H+D,4H+2D\}$,
		\item [(3)] $\{H,aH+D,(a+1)H+D,(a+2)H+D,(a+3)H+D,3H+2D,4H+2D\}$,
		\item [(4)] $\{aH+D,(a+1)H+D,(a+2)H+D,(a+3)H+D,2H+2D,3H+2D,4H+2D\}$,
		\item [(5)]  $\{-H+D,D,H+D,2H+2D,3H+2D,4H+2D,3H+3D\}$,
		\item [(6)]  $\{H,2H,3H+D,4H+D,5H+D,4H+2D,6H+D\}$,
		\item [(7)]  $\{H,2H+D,3H+D,4H+D,3H+2D,5H+D,4H+2D\}$,
		\item [(8)]  $\{H+D,2H+D,3H+D,2H+2D,4H+D,3H+2D,4H+2D\}$,
		\item [(9)]  $\{H,2H,H+D,3H,2H+D,3H+D,4H+D\}$,
		\item [(10)]  $\{H,D,2H,H+D,2H+D,3H+D,4H+2D\}$,
		\item [(11)]  $\{-H+D,H,D,H+D,2H+D,3H+2D,4H+2D\}$,
		\item [(12)]  $\{2H-D,H,2H,3H,4H+D,5H+D,6H+D\}$,
	\end{itemize}
	where $a\in \mathbb{Z}$. Moreover, by mutations and normalizations, they are related as:
	\begin{align}
	(1)\rightarrow (2)\rightarrow (3)\rightarrow (4), \nonumber
	\end{align}
	and 
	\begin{align}
	(5)\rightarrow (6)\rightarrow (7)\rightarrow (8)\rightarrow (9)\rightarrow (10)\rightarrow  (11)\rightarrow (12), \nonumber
	\end{align}
	\begin{proof}
		Write $D_0=0$. The sequence (\ref{5}) is an exceptional collection if and only if for any integers $0\leq i<j\leq 7$, the line bundles $\OO_X(D_i-D_j)$ are cohomologically zero.
		
		To make the sequence (\ref{5}) in Theorem \ref{cl5} an exceptional collection, $\OO_X(D_i)$ must be one of the line bundles in Proposition \ref{pro5}. To find out all the possible exceptional collections, we need Table \ref{tab:table5}.
		
		From Table \ref{tab:table5}, we observe that $\OO_X(D_1)$ may be one of $B_0$, $B_1$ and $B_7$. To finish the proof, we shall discuss case by case. However, since the proof is exactly same as the previous ones, we simply list the exceptional collections of line bundles of length $8$ from now on.
		\begin{itemize}
			\item [Case1.] Suppose $\OO_X(D_1)=B_0$. In this case, possible $\OO_X(D_2)$ may be $B_0$ or $B_1$. 
			\begin{itemize}
				\item [(a)] If $\OO_X(D_2)=B_0$, then we have the following exceptional collections:
				\begin{itemize}
					\item [(i)] $\{aH+D,(a+1)H+D,(a+2)H+D,(a+3)H+D,2H+2D,3H+2D,4H+2D\}$, which is $(4)$ in Theorem \ref{cl5},
					\item [(ii)] $\{-H+D,D,H+D,2H+2D,3H+2D,4H+2D,3H+3D\}$, which is $(5)$ in Theorem \ref{cl5},
					\item [(iii)] $\{H+D,2H+D,3H+D,2H+2D,4H+D,3H+2D,4H+2D\}$, which is $(8)$ in Theroem \ref{cl5}.
				\end{itemize}
				\item [(b)] If $\OO_X(D_2)=B_1$, then we only have one exceptional collection:
				$\{-H+D,H,D,H+D,2H+D,3H+2D,4H+2D\}$, which is $(11)$ in Theorem \ref{cl5}.
			\end{itemize} 
			\item [Case2.] Suppose $\OO_X(D_1)=B_1$. In this case, there are two possibilities for $D_2$:
			\begin{itemize}
				\item [(a)] When $\OO_X(D_2)=B_0$. In this case, we have three types of exceptional collections:
				\begin{itemize}
					\item [(i)] $\{H,aH+D,(a+1)H+D,(a+2)H+D,(a+3)H+D,4H+2D,5H+2D\}$, which is $(3)$ in Theorem \ref{cl5},
					\item [(ii)] $\{H,2H+D,3H+D,4H+D,3H+2D,5H+D,4H+2D\}$, which is $(7)$ in Theorem \ref{cl5},
					\item [(iii)] $\{H,D,2H,H+D,2H+D,3H+D,4H+2D\}$, which is $(10)$ in Theorem \ref{cl5}.
				\end{itemize}
				\item [(b)] When $\OO_X(D_2)=B_2$. In this case, we have four types of exceptional collections:
				\begin{itemize}
					\item [(i)] $\{H,2H,3H,aH+D,(a+1)H+D,(a+2)H+D,(a+3)H+D\}$, which is $(1)$ in Theorem \ref{cl5},
					\item [(ii)] $\{H,2H,aH+D,(a+1)H+D,(a+2)H+D,(a+3)H+D,5H+2D\}$, which is $(2)$ in Theorem \ref{cl5},
					\item [(iii)] $\{H,2H,3H+D,4H+D,5H+D,4H+2D,6H+D\}$, which is $(6)$ in Theorem \ref{cl5},
					\item [(iv)] $\{H,2H,H+D,3H,2H+D,3H+D,4H+D\}$, which is $(9)$ in Theorem \ref{cl5}.
				\end{itemize}
			\end{itemize}
		\item [Case3.] When $\OO_X(D_1)=B_7$. In this case we have only one type of exceptional collection:
		$\{2H-D,H,2H,3H,4H+D,5H+D,6H+D\}$, which is $(12)$ in Theorem \ref{cl5}.
		\end{itemize}

\begin{table}[t]
	\begin{center}
		\caption{Exceptional pairs.}
		\label{tab:table5}
		\begin{tabular}{c|c|c|c|c|c|c|c|c|c}
			& $B_0'$ & $B_1$ & $B_2$ & $B_3$ & $B_4$ & $B_5$ & $B_6$ & $B_7$ & $B_8$\\
			\hline
			$B_0$ &	$a'=a+1, a+2, a+3$ & $a=-1$  & $a=0$ & $a=1$ & \checkmark  & \checkmark & \checkmark &  & $a=1,0,-1$\\
			\hline
			$B_1$ &	\checkmark &  & \checkmark & \checkmark &  & \checkmark & \checkmark &  & \\
			\hline
			$B_2$ &	\checkmark &  &  & \checkmark &   &   & \checkmark &  & \\
			\hline
			$B_3$ & \checkmark &  &  &  &   &  &  &  & \\
			\hline
			$B_4$ & $a=4$ &  &  &  &  & \checkmark & \checkmark &  & \checkmark\\
			\hline
			$B_5$ & $a=5$ &  &  &  &  &  & \checkmark &  & \checkmark \\
			\hline
			$B_6$ & $a=6$ &  &  &  &  &  &  &  & \checkmark \\
			\hline
			$B_7$ & $a=4,5,6$ & \checkmark & \checkmark & \checkmark &  &  &  &  & \\
			\hline
			$B_8$ &  &  &  &  &  &  &  & \\
			\hline
		\end{tabular}
	\vspace{3mm}

	where $B_0=aH+D, B_1=H, B_2=2H, B_3=3H,B_4=2H+2D,$\\
	 $B_5=3H+2D, B_6=4H+2D, B_7=2H-D, B_8=3H+3D$.
	\end{center}
\end{table}
	\end{proof}
\end{thm}

\begin{thm}
	Let $X$ be the projective bundle $f:X=\PP_{\PP^3}(\OO\oplus \OO(1))\rightarrow \PP^3$. Then any exceptional collection of line bundles of length 8 on $X$ is full.
	\begin{proof}
		By Lemma \ref{nor}, we only need to show that any exceptional collection of line bundles of length $8$ in Theorem \ref{cl5} is full. By Lemma \ref{mt}, it suffices to show that the exceptional collections of type $(1)$ and $(5)$ in Theorem \ref{cl5} are full. 
		
		\begin{itemize}
			\item [(i)] To show the fullness of type $(1)$, we use Beilinson's semiorthogonal decomposition, 
			\begin{align}
			\D(\PP^3)=\langle\OO_{\PP^3}, \OO_{\PP^3}(1),\OO_{\PP^3}(2),\OO_{\PP^3}(3)\rangle. \nonumber
			\end{align}
			By Orlov's projective bundle formula (Theorem \ref{proj}), we obtain the following semiorthogonal decomposition of $\D(X)$
			\begin{align}
			\D(X)&=\langle f^*\D(\PP^3),f^*\D(\PP^3)\otimes \OO_X(1)\rangle \nonumber\\
			&=\langle f^*\OO_{\PP^3},f^*\OO_{\PP^3}(1),f^*\OO_{\PP^3}(2),f^*\OO_{\PP^3}(3),f^*\OO_{\PP^3}\otimes \OO_X(1),\nonumber\\
			&\qquad f^*\OO_{\PP^3}(1)\otimes\OO_X(1),f^*\OO_{\PP^3}(2)\otimes \OO_X(1),f^*\OO_{\PP^3}(3)\otimes \OO_X(1)\rangle\nonumber\\
			&=\langle \OO_X,\OO_X(H),\OO_X(2H),\OO_X(3H),\OO_X(aH+D),\OO_X((a+1)H+D),\nonumber\\
			&\qquad \OO_X((a+2)H+D),\OO_X((a+3)H+D)\rangle.\nonumber
			\end{align}
			Hence type $(1)$ is a full exceptional collection and so the exceptional collections of length $8$ of type $(2), (3)$ and $(4)$ are all full exceptional collections.
			\item [(ii)] To prove the fullness of $(5)$, we use the mutation. Since $a$ is arbitrary for the types $(1)-(4)$, we may choose any integer. Let $a=1$ in $(1)$. Then we have a full exceptional collection 
			\begin{align}
			\langle\OO_X, \OO_X(H), \OO_X(2H), \OO_X(3H),\OO_X(H+D),\OO_X(2H+D),\OO_X(3H+D),\OO_X(4H+D)\rangle. \nonumber
			\end{align}
			Here, we mutate $\OO_X(3H)$ and $\OO_X(H+D)$. Since $-2H+D$ is a cohomologically zero and by Lemma \ref{mt} (1), $L_{\OO_X(3H)}(H+D)\cong 3H$. Therefore, we have
			\begin{align}
			 \langle\OO_X, \OO_X(H), \OO_X(2H), \OO_X(H+D),\OO_X(3H),\OO_X(2H+D),\OO_X(3H+D),\OO_X(4H+D)\rangle \nonumber
			\end{align}  and this is also a full exceptional collection.
			Hence all the exceptional collections of length $8$ on $X$ are full.
		\end{itemize} 
	\end{proof}
\end{thm}

\subsection{When $X\cong \PP^1 \times \PP^3$.}	
Let $p_1: \PP^1 \times \PP^3\rightarrow \PP^1$ and $p_2: \PP^1 \times \PP^3\rightarrow \PP^3$ be the projections onto $\PP^1$ and $\PP^3$, respectively. Then the canonical divisor of $X$ is $K_X=-2H-4D$, where $H$ and $D$ be the pullback of the hyperplane class of $\PP^1$ and $\PP^3$, respectively.

\subsubsection{Classification of full exceptional collections}
\begin{thm} \label{cl6}
	Let $X$ be $\PP^1 \times \PP^3$. Then the normalized sequence
	\begin{align}
	\{\OO_X,\OO_X(D_1),\cdots,\OO_X(D_7)\} \label{6} \tag{F}
	\end{align}
	is an exceptional collection of line bundles if and only if the ordered set of divisors $\{D_1,D_2,\cdots,D_7\}$ is one of the following types:
	\begin{itemize}
		\item [(1)] $\{H+aD,H+(a+1)D,H+(a+2)D,H+(a+3)D,2H+D,2H+2D,2H+3D\}$,
		\item [(2)] $\{D,2D,3D,H+aD,H+(a+1)D,H+(a+2)D,H+(a+3)D\}$,
		\item [(3)] $\{D,2D,H+aD,H+(a+1)D,H+(a+2)D,H+(a+3)D,2H+3D\}$,
		\item [(4)] $\{D,H+aD,H+(a+1)D,H+(a+2)D,H+(a+3)D,2H+2D,2H+3D\}$,
		\item [(5)]  $\{H-2D, H-D, H, bH+D, (b+1)H+D, 2H+2D, 2H+3D\}$,
		\item [(6)]  $\{D, 2D, dH+3D, (d+1)H+3D, H+4D, H+5D, H+6D\}$,
		\item [(7)]  $\{D, cH+2D, (c+1)H+2D, H+3D, H+4D, H+5D, 2H+3D\}$,
		\item [(8)]  $\{bH+D, (b+1)H+D, H+2D, H+3D, H+4D, 2H+2D, 2H+3D\}$,
		\item [(9)]  $\{H, bH+D, bH+2D, bH+3D, (b+1)H+D, (b+1)H+2D, (b+1)H+3D\}$,
		\item [(10)] $\{bH+D, bH+2D, bH+3D, (b+1)H+D, (b+1)H+2D, (b+1)H+3D, H+4D\}$,
		\item [(11)] $\{D, 2D, H, H+D, H+2D, dH+3D, (d+1)H+3D\}$,
		\item [(12)] $\{D, H-D, H, H+D, cH+2D, (c+1)H+2D, 2H+3D\}$.
		\item [(13)] $\{bH+D, bH+2D, (b+1)H+D, (b+1)H+2D, dH+3D, (d+1)H+3D, H+4D\}$,
		\item [(14)] $\{D, H, H+D, bH+2D, (b+1)H+2D, dH+3D, (d+1)H+3D\}$,
		\item [(15)] $\{H-D, H, bH+D, (b+1)H+D, cH+2D, (c+1)H+2D, 2H+3D\}$,
		\item [(16)] $\{D, bH+2D, (b+1)H+2D, dH+3D, (d+1)H+3D, H+4D, H+5D\}$,
		\item [(17)] $\{bH+D, (b+1)H+D, dH+2D, (d+1)H+2D, H+3D, H+4D, 2H+3D\}$,
		\item [(18)] $\{H, bH+D, (b+1)H+D, cH+2D, cH+3D, (c+1)H+2D, (c+1)H+3D\}$,
		\item [(19)] $\{bH+D, (b+1)H+D, cH+2D, cH+3D, (c+1)H+2D, (c+1)H+3D, H+4D\}$,
		\item [(20)] $\{H, bH+D, bH+2D, (b+1)H+D, (b+1)H+2D, dH+3D, (d+1)H+3D\}$,
		\item [(21)] $\{bH+D, (b+1)H+D, cH+2D, (c+1)H+2D, dH+3D, (d+1)H+3D, H+4D\}$,
		\item [(22)] $\{H, bH+D, (b+1)H+D, cH+2D, (c+1)H+2D, dH+3D, (d+1)H+3D\}$,
		\item [(23)] $\{D, cH+2D, cH+3D, (c+1)H+2D, (c+1)H+3D, H+4D, H+5D\}$,
		\item [(24)] $\{cH+D, cH+2D, (c+1)H+D, (c+1)H+2D, H+3D, H+4D, 2H+3D\}$,
		\item [(25)] $\{D, H, H+D, cH+2D, cH+3D, (c+1)H+2D, (c+1)H+3D\}$,
		\item [(26)] $\{H-D, H, cH+D, cH+2D, (c+1)H+D, (c+1)H+2D, 2H+3D\}$,
		\item [(27)] $\{bH+D, bH+2D, (b+1)H+D, bH+3D, (b+1)H+2D, (b+1)H+3D, H+4D\}$,
		\item [(28)] $\{D, H, 2D, H+D, H+2D, dH+3D, (d+1)H+3D\}$,
		\item [(29)] $\{H-D, D, H, H+D, cH+2D, (c+1)H+2D, 2H+3D\}$,
		\item [(30)] $\{-H+2D, D, 2D, dH+3D, (d+1)H+3D, H+4D, H+5D\}$,
		\item [(31)] $\{H-D, H, bH+D, (b+1)H+D, 2H+2D, 2H+3D, H+5D\}$
		\item [(32)] $\{D, cH+2D, (c+1)H+2D, H+3D, H+4D, 2H+3D, H+5D\}$,
		\item [(33)] $\{bH+D, (b+1)H+D, H+2D, H+3D, 2H+2D, H+4D, 2H+3D\}$,
		\item [(34)] $\{H, bH+D, bH+2D, (b+1)H+D, bH+3D, (b+1)H+2D, (b+1)H+3D\}$,
		\item [(35)] $\{H+D, H+2D, H+3D, 2H+D, H+4D, 2H+2D, 2H+3D\}$,
		\item [(36)] $\{D, 2D, H, 3D, H+D, H+2D, H+3D\}$,
		\item [(37)] $\{D, H-D, 2D, H, H+D, H+2D, 2H+3D\}$,
		\item [(38)] $\{H-2D, D, H-D, H, H+D, 2H+2D, 2H+3D\}$,
		\item [(39)] $\{-H+3D, D, 2D, 3D, H+4D, H+5D, H+6D\}$,
		\item [(40)] $\{H-2D, H-D, H, 2H+D, 2H+2D, 2H+3D, 3H+D\}$,
		\item [(41)] $\{D, 2D, H+3D, H+4D, H+5D, 2H+3D, H+6D\}$,
		\item [(42)] $\{D, H+2D, H+3D, H+4D, 2H+2D, H+5D, 2H+3D\}$,
		\item [(43)] $\{H+D, H+2D, H+3D, 2H+D, 2H+2D, H+4D, 2H+3D\}$,
		\item [(44)] $\{D, 2D, H, H+D, 3D, H+2D, H+3D\}$,
		\item [(45)] $\{D, H-D, H, 2D, H+D, H+2D, 2H+3D\}$,
		\item [(46)] $\{H-2D, H-D, D, H, H+D, 2H+2D, 2H+3D\}$,
		\item [(47)] $\{D, -H+3D, 2D, 3D, H+4D, H+5D, H+6D\}$,
		\item [(48)] $\{-H+2D, D, 2D, H+3D, H+4D, H+5D, 2H+3D\}$,
		\item [(49)] $\{H-D, H, 2H+D, 2H+2D, 2H+3D, 3H+D, 3H+2D\}$,
		\item [(50)] $\{D, H+2D, H+3D, H+4D, 2H+2D, 2H+3D, H+5D\}$,
		\item [(51)] $\{H-2D, H-D, H, 2H+D, 2H+2D, 3H+D, 2H+3D\}$,
		\item [(52)] $\{D, 2D, H+3D, H+4D, 2H+3D, H+5D, H+6D\}$,
		\item [(53)] $\{D, H+2D, H+3D, 2H+2D, H+4D, H+5D, 2H+3D\}$,
		\item [(54)] $\{H+D, H+2D, 2H+D, H+3D, H+4D, 2H+2D, 2H+3D\}$,
		\item [(55)] $\{D, H, 2D, 3D, H+D, H+2D, H+3D\}$,
		\item [(56)] $\{H-D, D, 2D, H, H+D, H+2D, 2H+3D\}$,
		\item [(57)] $\{-H+2D, -H+3D, D, 2D, 3D, H+4D, H+5D\}$,
		\item [(58)] $\{D, H-D, H, H+D, 2H+2D, 2H+3D, 3H+2D\}$,
		\item [(59)] $\{H+D, H+2D, 2H+D, H+3D, 2H+2D, H+4D, 2H+3D\}$,
		\item [(60)] $\{D, H, 2D, H+D, 3D, H+2D, H+3D\}$,
		\item [(61)] $\{H-D, D, H, 2D, H+D, H+2D, 2H+3D\}$,
		\item [(62)] $\{-H+2D, D, -H+3D, 2D, 3D, H+4D, H+5D\}$,
		\item [(63)] $\{H-D, D, H, H+D, 2H+2D, 2H+3D, 3H+2D\}$,
		\item [(64)] $\{-H+2D, D, 2D, H+3D, H+4D, 2H+3D, H+5D\}$,
		\item [(65)] $\{H-D, H, 2H+D, 2H+2D, 3H+D, 2H+3D, 3H+2D\}$,
		\item [(66)] $\{D, H+2D, H+3D, 2H+2D, H+4D, 2H+3D, H+5D\}$.
	\end{itemize}
	where $a,b,c,d\in \mathbb{Z}$. Moreover, by mutations and normalizations, they are related as:
	\begin{align}
	&(1)\rightarrow (2)\rightarrow (3)\rightarrow (4)\rightarrow (1), \nonumber\\
	&(5)\rightarrow (6)\rightarrow (7)\rightarrow (8)\rightarrow (9)\rightarrow (10)\rightarrow  (11)\rightarrow (12)\rightarrow (5), \nonumber\\
	&(13)\rightarrow (14)\rightarrow (15)\rightarrow (16)\rightarrow (17)\rightarrow (18)\rightarrow  (19)\rightarrow (20)\rightarrow (13), \nonumber\\
	&(21)\rightarrow (22)\rightarrow (21),\nonumber\\
	&(22)\rightarrow (24)\rightarrow (25)\rightarrow (26)\rightarrow (22),\nonumber\\
	&(27)\rightarrow (28)\rightarrow (29)\rightarrow (30)\rightarrow (31)\rightarrow (32)\rightarrow  (33)\rightarrow (34)\rightarrow (27), \nonumber\\
	&(35)\rightarrow (36)\rightarrow (37)\rightarrow (38)\rightarrow (39)\rightarrow (40)\rightarrow  (41)\rightarrow (42)\rightarrow (35),\nonumber\\
	&(43)\rightarrow (44)\rightarrow (45)\rightarrow (46)\rightarrow (47)\rightarrow (48)\rightarrow  (49)\rightarrow (50)\rightarrow (43),\nonumber\\
	&(51)\rightarrow (52)\rightarrow (53)\rightarrow (54)\rightarrow (55)\rightarrow (56)\rightarrow  (57)\rightarrow (58)\rightarrow (51),\nonumber\\
	&(59)\rightarrow (60)\rightarrow (61)\rightarrow (62)\rightarrow (63)\rightarrow (64)\rightarrow  (65)\rightarrow (66)\rightarrow (59).\nonumber
	\end{align}
	\begin{proof}
		Note that $\OO_X(aH+bD)$ is cohomologically zero line bundle if $a=-1$ or $b=-1,-2,-3$. Based on Table \ref{tab:table6} and the rules below, we get the result.
		\begin{itemize}
			\item [(1)] $B_0$ can occur at four times.
			\item [(2)] $B_1, B_2$ and $B_3$ can occur at most two times.
			\item [(3)] $B_1$ cannot occur twice after $B_2$ or $B_3$.
			\item [(4)] $B_2$ cannot occur twice after $B_3$.
			\item [(5)] $B_1$ cannot occur after two $B_2$ or $B_3$.
			\item [(6)] $B_2$ cannot occur after two $B_3$.
		\end{itemize}
	\begin{table}[h!]
		\begin{center}
			\caption{Exceptional pairs.}
			\label{tab:table6}
			\begin{adjustbox}{max width=\textwidth}
				\begin{tabular}{c|c|c|c|c} 
					& $B_0'$ & $B_1'$ & $B_2'$ & $B_3'$\\
					\hline
					$B_0$ &	$a'=a+1, a+2, a+3$ & $a=0,-1,-2$ or $b'=2$  & $a=-1,0,1$ or $c'=2$  & $a=0,1,2$ or $d'=2$  \\
					\hline
					$B_1$ &	$a'=2,3,4$ or $b=0$  & $b'=b+1$ & \checkmark & \checkmark\\
					\hline
					$B_2$ &	$a'=3,4,5$ or $c=0$  & $b'=c+1$ & $c'=c+1$ & \checkmark\\
					\hline
					$B_3$ & $a'=4,5,6$ or $d=0$  & $b'=d+1$ & $c'=d+1$ & $d'=d+1$\\
					\hline
				\end{tabular}
			\end{adjustbox}
			\vspace{3mm}
			
			where $B_0=H+aD, B_1=bH+D, B_2=cH+2D, B_3=dH+3D$.
		\end{center}
	\end{table}
	\end{proof}
\end{thm}

\begin{thm} \label{full1}
	Let $X$ be $\PP^1 \times \PP^3$. Then any exceptional collection of line bundles of length 6 on $X$ is full.
	\begin{proof}
		By Lemma \ref{nor}, we only need to show that any exceptional collection of line bundles of length $8$ in Theorem \ref{cl6} is full. By Lemma \ref{mt}, it suffices to show that the exceptional collections of type $(2), (22), (13)$, $(27)$, $(23)$, $(12)$, $(35)$, $(43)$, $(59)$ and $(54)$ in Theorem \ref{cl2} are full. 
		
		We first show that exceptional collection of type $(2)$ is full. Consider the projective bundle $f:X\cong \PP_{\PP^3}(\OO(-a)\oplus \OO(-a))\rightarrow \PP^3$, where $\OO_X(1)=H+aD$. Then by the projective bundle formula, we have the following semiorthogonal decomposition of $\D(X)$.
		\begin{align}
		\D(X)&=\langle f^*\D(\PP^3),f^*\D(\PP^3)\otimes \OO_X(1)\rangle \nonumber\\
		&=\langle f^*\OO_{\PP^3},f^*\OO_{\PP^3}(1),f^*\OO_{\PP^3}(2),f^*\OO_{\PP^3}(3),\nonumber\\
		&\qquad f^*\OO_{\PP^3}\otimes \OO_X(1),f^*\OO_{\PP^3}(1)\otimes\OO_X(1),f^*\OO_{\PP^3}(2)\otimes \OO_X(1),f^*\OO_{\PP^3}(3)\otimes \OO_X(1)\rangle\nonumber\\
		&=\langle \OO_X,\OO_X(D),\OO_X(2D),\OO_X(3D),\nonumber\\
		&\qquad \OO_X(H+aD),\OO_X(H+(a+1)D,\OO_X(H+(a+2)D),\OO_X(H+(a+3)D))\rangle.\nonumber
		\end{align}
		Hence type $(2)$ is a full exceptional collection.
		
		Now we show that exceptional collection of type $(22)$ is full. Consider the projective bundle $f:X\cong \PP_{\PP^1}(\OO(-b)\oplus \OO(-b)\oplus \OO(-b)\oplus \OO(-b))\rightarrow \PP^1$, where $\OO_X(1)=bH+D$. Then by the projective bundle formula, we obtain the following semiorthogonal decomposition of $\D(X)$.
		\begin{align}
		\D(X)&=\langle f^*\D(\PP^1),f^*\D(\PP^1)\otimes \OO_X(1),f^*\D(\PP^1)\otimes \OO_X(2),f^*\D(\PP^1)\otimes \OO_X(3)\rangle \nonumber\\
		&=\langle f^*\OO_{\PP^1},f^*\OO_{\PP^1}(1),f^*\OO_{\PP^1}\otimes \OO_X(1),f^*\OO_{\PP^1}(1)\otimes\OO_X(1),\nonumber\\
		&\qquad f^*\OO_{\PP^1}\otimes \OO_X(2),f^*\OO_{\PP^1}(1)\otimes \OO_X(2),f^*\OO_{\PP^1}\otimes \OO_X(3),f^*\OO_{\PP^1}(1)\otimes \OO_X(3)\rangle\nonumber\\
		&=\langle \OO_X,\OO_X(D),\OO_X(H+bD),\OO_X(H+(b+1)D),\nonumber\\
		&\qquad \OO_X(2H+2bD),\OO_X(2H+(2b+1)D),\OO_X(3H+3bD),\OO_X(3H+(3b+1)D)\rangle.\nonumber
		\end{align}
		Based on this semiorthogonal decomposition, we show that the exceptional collection of type $(5)$ is full by mathematical induction. For a given $b\in \mathbb{Z}$, assume that for $d=k$ the exceptional collection
		\begin{align}\label{mi4}
		\{H, bH+D, (b+1)H+D, 2bH+2D, (2b+1)H+2D, kH+3D,(k+1)H+3D\} \tag{E4}
		\end{align}
		is full. Then for $d=k-1$, we have the exceptional collection
		\begin{align}\label{mi5} \tag{E5}
		\{H, bH+D, (b+1)H+D, 2bH+2D, (2b+1)H+2D, (k-1)H+3D,kH+3D\}
		\end{align}
		and for $d=k+1$, we have the exceptional collection
		\begin{align}\label{mi6} \tag{E6}
		\{H, bH+D, (b+1)H+D, 2bH+2D, (2b+1)H+2D (k+1)H+3D,(k+2)H+3D\}.
		\end{align}
		Comparing the exceptional collections (\ref{mi5}) and (\ref{mi6}) with (\ref{mi4}) and by Lemma \ref{mt1}, the exceptional collections (\ref{mi5}) and (\ref{mi6}) are full. Similarly for fixed $b,d\in \mathbb{Z}$, we are able to show that the exceptional collection of type $(22)$ is full. Hence the exceptional collection of type $(22)$ is also full.
		
		For the other exceptional collections, we summarize it in the Table \ref{tab:table12}.
		
		\begin{table}[h!]
			\begin{center}
				\caption{Showing the fullness.}
				\label{tab:table12}
					\begin{tabular}{c|c} 
						Type & How to show \\
						\hline
						$(13)$ & Set $c=b$ in $(21)$ and mutate $(b+1)H+D$ and $bH+2D$  \\
						\hline
						$(27)$ & Set $d=b$ in $(13)$ and mutate $(b+1)H+2D$ and $bH+3D$  \\
						\hline
						$(23)$ & Mutate $(b+1)H+2D$ and $bH+3D$ in $(16)$  \\
						\hline
						$(12)$ & Mutate $H-D$ and $H$ in $(29)$  \\
						\hline
						$(35)$ & Set $a=1$ in $(1)$ and mutate $H+4D$ and $2H+D$  \\
						\hline
						$(43)$ & Mutate $H+4D$ and $2H+2D$ in $(35)$  \\
						\hline
						$(59)$ & Mutate $H+3D$ and $2H+D$ in $(43)$  \\
						\hline
						$(54)$ & Mutate $2H+2D$ and $H+4D$ in $(59)$  \\
						\hline
					\end{tabular}
				\end{center}
			\end{table}
		\newpage
		Therefore, all the exceptional collections of length $8$ are full.
	\end{proof}
\end{thm}

\subsection{When $X\cong \PP_{\PP^1}(\OO\oplus\OO\oplus\OO\oplus\OO(1))$.}	
Let $f:\PP_{\PP^1}(\OO\oplus \OO\oplus \OO \oplus\OO(1))\rightarrow \PP^1$ be the projective bundle. Then its canonical divisor is $K_X=-3H-4D$ where $H$ is the pullback of hyperplane class in $\PP^3$ and $D$ is the class of $c_1(\OO_X(1))$. The Picard group of $X$ is 
\begin{align}
\Pic(X)\cong \Pic(\PP^3)\oplus \mathbb{Z}[D]=\mathbb{Z}[H]\oplus\mathbb{Z}[D], \nonumber
\end{align}
with intersection numbers
\begin{align}
H^4=0, ~H^3D=0, ~H^2D^2=0, ~HD^3=1, ~D^4=-1. \nonumber
\end{align}
\subsubsection{Cohomologically zero line bundles}
\begin{lem}\label{CZ7}
	$H^0(X,aH+bD)=0$ if and only if $b<0$ or $a+b<0, ~b\geq 0$. Also, $H^4(X,aH+bD)=0$ if and only if $b>-4$ or $a+b>-7, ~b\leq -4$.
	\begin{proof}
		Note that we have the following isomorphisms:
		\begin{align}
		f_*(\OO_X(b))&\cong \Sym^b(\OO_{\PP^1}\oplus\OO_{\PP^1}\oplus\OO_{\PP^1}\oplus\OO_{\PP^1}(1)) \nonumber \\ 
		&\cong \bigoplus_{i=0}^b \Sym^i(\OO_{\PP^1}\oplus\OO_{\PP^1}\oplus\OO_{\PP^1})\otimes \Sym^{b-i}(\OO_{\PP^1}(1)) \nonumber \\
		&\cong \bigoplus_{i=0}^b \OO_{\PP^1}(b-i) ~\text{if}~ b\geq 0, \nonumber
		\end{align} and $f_*(\OO_X(b))\cong 0$ if $b<0$.
		
		By projection formula, we have the following isomorphisms:
		\begin{align}
		H^0(X, f^*(\OO_{\PP^1}(a))\otimes \OO_X(b))&\cong H^0(\PP^1, \OO_{\PP^1}(a)\otimes  \bigoplus_{i=0}^b \OO_{\PP^1}(b-2)) \nonumber \\
		&\cong \bigoplus_{i=0}^b H^0(\PP^1,\OO_{\PP^1}(a+b-i)) ~\text{if}~ b\geq 0,\nonumber
		\end{align}
		and if $b<0$, then $H^0(X,aH+bD)=0$.
		Hence, $H^0(X,aH+bD)=0$ if and only if $b<0$ or $a+b<0, ~b\geq 0$.
		
		By Serre duality, one obtains the condition for $H^4(X,aH+bD)=0$.
	\end{proof}
\end{lem}
\begin{pro} \label{pro7}
	A line bundle $\OO_X(aH+bD)$ is cohomologically zero if one of the following holds:
	\begin{itemize}
		\item [(1)] $b=-1$,
		\item [(2)] $b=-2$,
		\item [(3)] $b=-3$,
		\item [(4)] $a=-1, b=0$,
		\item [(5)] $a=-2, b=-4$.
	\end{itemize}
	\begin{proof}
		Since $c(X)=(1+2H)\{(1+D)^4+H(1+D)^3\}$, we have $c_1=3H+4D$, $c_2=2H^2+11HD+6D^2$, $c_3=6H^2D+15HD^2+4D^3$ and $c_4=6H^2D^2+9HD^3+D^4=8$.
		Thus, by Lemma \ref{HRR}, we have
		\begin{align}
		\chi(\OO_X(aH+bD))&=\frac{1}{24}(24+24a-24b+44ab+13b^2+24ab^2-2b^3+4ab^3-b^4) \nonumber\\
		&=\frac{1}{24}(b+1)(b+2)(b+3)(4a-b+4). \nonumber
		\end{align}
		Therefore, by Lemma \ref{CZ7}, $\OO_X(aH+bD)$ is cohomologically zero if the following conditions hold:
		\begin{itemize}
			\item [$\mathrm{(i)}$] $b<0$ or $a+b<0, ~b\geq 0$,
			\item [$\mathrm{(ii)}$] $b>-4$ or $a+b>-7, ~b\leq -4$,
			\item [$\mathrm{(iii)}$] $(b+1)(b+2)(b+3)(4a-b+4)=0$.
		\end{itemize}
		Note that the equation $4a-b+4=0$ has $2$ integer solutions; $(a,b)=(-1,0), (-2,-4)$ under conditions $\mathrm{(i)}$ and $\mathrm{(ii)}$. This completes the proof.
	\end{proof}
\end{pro}

\subsubsection{Classification of full exceptional collections}
\begin{thm} \label{cl7}
	Let $X$ be the projective bundle $f:\PP_{\PP^1}(\OO\oplus \OO\oplus \OO \oplus\OO(1))\rightarrow \PP^1$. Then the normalized sequence
	\begin{align}
	\{\OO_X,\OO_X(D_1),\cdots,\OO_X(D_7)\} \label{7} \tag{G}
	\end{align}
	is an exceptional collection of line bundles if and only if the ordered set of divisors $\{D_1,D_2,\cdots,D_7\}$ is one of the following types:
	\begin{itemize}
		\item [(1)] $\{H,aH+D,(a+1)H+D,bH+2D,(b+1)H+2D,cH+3D,(c+1)H+3D\}$,
		\item [(2)] $\{aH+D,(a+1)H+D,bH+2D,(b+1)H+2D,cH+3D,(c+1)H+3D,2H+4D\}$,
	\end{itemize}
	where $a,b,c\in \mathbb{Z}$. Moreover, by mutations and normalizations, they are related as:
	\begin{align}
	(1)\rightarrow (2)\rightarrow (1). \nonumber
	\end{align}
	\begin{proof}
		Write $D_0=0$. The sequence (\ref{7}) is an exceptional collection if and only if for any integers $0\leq i<j\leq 7$, the line bundles $\OO_X(D_i-D_j)$ are cohomologically zero.
		
		To make the sequence (\ref{7}) an exceptional collection, $\OO_X(D_i)$ must be one of the line bundles in Proposition \ref{pro7}. To find out all the possible exceptional collections, we need Table \ref{tab:table7}.
		
		From Table \ref{tab:table7}, we observe that $\OO_X(D_1)$ must be one of $B_0$ or $B_3$. To finish the proof, we shall discuss case by case:
		\begin{itemize}
			\item [Case1.] Suppose $\OO_X(D_1)=B_0$. In this case, the only possible combination is $\{aH+D,(a+1)H+D,bH+2D,(b+1)H+2D,cH+3D,(c+1)H+3D,4H+2D\}$, which is $(2)$ in Theorem \ref{cl7}.
			\item [Case2.] Suppose $\OO_X(D_1)=B_3$. In this case, the only possible choice for $D_2$ is $B_0$:
			$\{H,aH+D,(a+1)H+D,bH+2D,(b+1)H+2D,cH+3D,(c+1)H+3D\},$
			which is $(1)$ in Theorem \ref{cl7}.
		\end{itemize}
		\begin{table}[h!]
			\begin{center}
				\caption{Exceptional pairs.}
				\label{tab:table7}
				\begin{tabular}{c|c|c|c|c|c}
					& $B_0'$ & $B_1'$ & $B_2'$ & $B_3$ & $B_4$\\
					\hline
					$B_0$ & $a'=a+1$ & \checkmark & \checkmark &  & \checkmark\\
					\hline
					$B_1$ &	 & $b'=b+1$ & \checkmark &  & \checkmark\\
					\hline
					$B_2$ &	 &  & $c'=c+1$ &  & \checkmark\\
					\hline
					$B_3$ & \checkmark & \checkmark & \checkmark &  &\\
					\hline
					$B_4$ &  &  &  &  &\\
					\hline
				\end{tabular}
			\vspace{3mm}
			
			where $B_0=aH+D, B_1=bH+2D, B_2=cH+3D, B_3=H, B_4=2H+4D$.
			\end{center}
		\end{table}
	\end{proof}
\end{thm}

\begin{thm}
	Let $X$ be the projective bundle $f:\PP_{\PP^1}(\OO\oplus \OO\oplus \OO \oplus\OO(1))\rightarrow \PP^1$. Then any exceptional collection of line bundles of length 8 on $X$ is full.
	\begin{proof}
		By Lemma \ref{nor}, we only need to show that any exceptional collection of line bundles of length $6$ in Theorem \ref{cl7} is full. By Lemma \ref{mt}, it suffices to show that the exceptional collection of type $(1)$ in Theorem \ref{cl7} is full. 
		
		By Beilinson's semiorthogonal decomposition, we have 
		\begin{align}
		\D(\PP^1)=\langle \OO_{\PP^1}, \OO_{\PP^3}(1)\rangle. \nonumber
		\end{align}
		By Orlov's projective bundle formula (Theorem \ref{proj}), we obtain the following semiorthogonal decomposition of $\D(X)$
		\begin{align}
		\D(X)&=\langle f^*\D(\PP^1),f^*\D(\PP^1)\otimes \OO_X(1),f^*\D(\PP^1)\otimes \OO_X(2),f^*\D(\PP^1)\otimes \OO_X(3)\rangle \nonumber\\
		&=\langle f^*\OO_{\PP^1},f^*\OO_{\PP^1}(1),f^*\OO_{\PP^1}\otimes \OO_X(1),f^*\OO_{\PP^1}(1)\otimes \OO_X(1),\nonumber\\
		&\qquad f^*\OO_{\PP^1}\otimes\OO_X(2),f^*\OO_{\PP^1}(1)\otimes \OO_X(2),f^*\OO_{\PP^1}\otimes \OO_X(3),f^*\OO_{\PP^1}(1)\otimes \OO_X(3)\rangle\nonumber\\
		&=\langle \OO_X,\OO_X(H),\OO_X(aH+D),\OO_X((a+1)H+D),\nonumber\\
		&\qquad \OO_X(2aH+2D),\OO_X((2a+1)H+2D),\OO_X(3aH+3D),\OO_X((3a+1)H+3D)\rangle.\nonumber
		\end{align}
		Based on this semiorthogonal decomposition, we may prove that the exceptional collection of type $(1)$ is full by mathematical induction  exactly the same argument as in the proof of Theorem \ref{full1}.
		Hence type $(1)$ is a full exceptional collection and so all the exceptional collections of length $8$ are full exceptional collections.
	\end{proof}
\end{thm}

\subsection{When $X\cong \PP_{\PP^2}(\OO\oplus\OO\oplus\OO(2))$.}	
Let $f:\PP_{\PP^2}(\OO\oplus \OO \oplus\OO(2))\rightarrow \PP^2$ be the projective bundle. Then its canonical divisor is $K_X=-5H-3D$ where $H$ is the pullback of hyperplane class in $\PP^3$ and $D$ is the class of $c_1(\OO_X(1))$. The Picard group of $X$ is 
\begin{align}
\Pic(X)\cong \Pic(\PP^2)\oplus \mathbb{Z}[D]=\mathbb{Z}[H]\oplus\mathbb{Z}[D], \nonumber
\end{align}
with intersection numbers
\begin{align}
H^4=0, ~H^3D=0, ~H^2D^2=1, ~HD^3=-2, ~D^4=4. \nonumber
\end{align}
\subsubsection{Cohomologically zero line bundles}
\begin{lem}\label{CZ8}
	$H^0(X,aH+bD)=0$ if and only if $b<0$ or $a+2b<0, ~b\geq 0$. Also, $H^4(X,aH+bD)=0$ if and only if $b>-2$ or $a+b>-11, ~b\leq -2$.
	\begin{proof}
		Note that we have the following isomorphisms:
		\begin{align}
		f_*(\OO_X(b))&\cong \Sym^b(\OO_{\PP^2}\oplus \OO_{\PP^2} \oplus\OO_{\PP^2}(2)) \nonumber \\ 
		&\cong \bigoplus_{i=0}^b \Sym^i(\OO_{\PP^2}\oplus \OO_{\PP^2})\otimes \Sym^{b-i}(\OO_{\PP^2}(2)) \nonumber \\
		&\cong \bigoplus_{i=0}^b \OO_{\PP^2}(2b-2i) ~\text{if}~ b\geq 0, \nonumber
		\end{align} and $f_*(\OO_X(b))\cong 0$ if $b<0$.
		
		By projection formula, we have the following isomorphisms:
		\begin{align}
		H^0(X, f^*(\OO_{\PP^2}(a))\otimes \OO_X(b))&\cong H^0(\PP^2, \OO_{\PP^2}(a)\otimes  \bigoplus_{i=0}^b \OO_{\PP^2}(2b-2i)) \nonumber \\
		&\cong \bigoplus_{i=0}^b H^0(\PP^2,\OO_{\PP^2}(a+2b-2i)) ~\text{if}~ b\geq 0,\nonumber
		\end{align}
		and if $b<0$, then $H^0(X,aH+bD)=0$.
		Hence, $H^0(X,aH+bD)=0$ if and only if $b<0$ or $a+2b<0, ~b\geq 0$.
		
		By Serre duality, one obtains the condition for $H^4(X,aH+bD)=0$.
	\end{proof}
\end{lem}
\begin{pro} \label{pro8}
	A line bundle $\OO_X(aH+bD)$ is cohomologically zero if one of the following holds:
	\begin{itemize}
		\item [(1)] $b=-1$,
		\item [(2)] $b=-2$,
		\item [(3)] $a=-1, b=0$,
		\item [(4)] $a=-2, b=0$,
		\item [(5)] $a=-3, b=-3$,
		\item [(6)] $a=-4, b=-3$.
	\end{itemize}
	\begin{proof}
		Since $c(X)=(1+3H+3H^2)\{(1+D)^3+2H(1+D)^2\}$, we have $c_1=5H+3D$, $c_2=9H^2+13HD+3D^2$, $c_3=6H^3+21H^2D+11HD^2+D^3$ and $c_4=12H^3D+15H^2D^2+3HD^3=9$.
		Thus, by Lemma \ref{HRR}, we have
		\begin{align}
		\chi(\OO_X(aH+bD))&=\frac{1}{12}(12+18a+10b+6a^2+19ab-2b^2 \nonumber\\
		&\qquad +9a^2b-3ab^2+2b^3+3a^2b^2-4ab^3+2b^4) \nonumber\\
		&=\frac{1}{12}(b+1)(b+2)(3a^2-4ab+9a+2b^2-4b+6). \nonumber
		\end{align}
		Therefore, by Lemma \ref{CZ8}, $\OO_X(aH+bD)$ is cohomologically zero if the following conditions hold:
		\begin{itemize}
			\item [$\mathrm{(i)}$] $b<0$ or $a+2b<0, ~b\geq 0$,
			\item [$\mathrm{(ii)}$] $b>-3$ or $a+2b>-11, ~b\leq -3$,
			\item [$\mathrm{(iii)}$] $(b+1)(b+2)(3a^2-4ab+9a+2b^2-4b+6)=0$.
		\end{itemize}
		Note that the equation $3a^2-4ab+9a+2b^2-4b+6=0$ has $4$ integer solutions; $(a,b)=(-1,0),(-2,0),(-3,-3),(-4,-3)$ under conditions $\mathrm{(i)}$ and $\mathrm{(ii)}$. This completes the proof.
	\end{proof}
\end{pro}

\subsubsection{Classification of full exceptional collections}
\begin{thm} \label{cl8}
	Let $X$ be the projective bundle $f:\PP_{\PP^2}(\OO\oplus \OO \oplus\OO(2))\rightarrow \PP^2$. Then the normalized sequence
	\begin{align}
	\{\OO_X,\OO_X(D_1),\cdots,\OO_X(D_8)\} \label{8} \tag{H}
	\end{align}
	is an exceptional collection of line bundles if and only if the ordered set of divisors $\{D_1,D_2,\cdots,D_8\}$ is one of the following types:
	\begin{itemize}
		\item [(1)] $\{H,2H,aH+D,(a+1)H+D,(a+2)H+D,bH+2D,(b+1)H+2D,(b+2)H+2D\}$,
		\item [(2)]
		$\{H,aH+D,(a+1)H+D,(a+2)H+D,bH+2D,(b+1)H+2D,(b+2)H+2D,4H+3D\}$
		\item [(3)] $\{aH+D,(a+1)H+D,(a+2)H+D,bH+2D,(b+1)H+2D,(b+2)H+2D,3H+3D,4H+3D\}$,
	\end{itemize}
	where $a,b\in \mathbb{Z}$. Moreover, by mutations and normalizations, they are related as:
	\begin{align}
	(1)\rightarrow (2)\rightarrow (3)\rightarrow (1). \nonumber
	\end{align}
	\begin{proof}
		Write $D_0=0$. The sequence (\ref{8}) is an exceptional collection if and only if for any integers $0\leq i<j\leq 8$, the line bundles $\OO_X(D_i-D_j)$ are cohomologically zero.
		
		To make the sequence (\ref{8}) an exceptional collection, $\OO_X(D_i)$ must be one of the line bundles in Proposition \ref{pro8}. To find out all the possible exceptional collections, we need Table \ref{tab:table8}.
		
		From Table \ref{tab:table8}, we observe that $\OO_X(D_1)$ may be one of $B_0$ and $B_3$. To finish the proof, we shall discuss case by case:
		\begin{itemize}
			\item [Case1.] Suppose $\OO_X(D_1)=B_0$. In this case, the only possible combination is $\{aH+D,(a+1)H+D,(a+2)H+D,bH+2D,(b+1)H+2D,(b+2)H+2D,3H+3D,4H+3D\}$, which is $(3)$ in Theorem \ref{cl8}.
			\item [Case2.] Suppose $\OO_X(D_1)=B_3$. In this case, there are two possible choices for $D_2$:
			\begin{itemize}
			\item [(a)]
			If $\OO_X(D_2)=B_4$, we have $\{H,2H,aH+D,(a+1)H+D,(a+2)H+D,bH+2D,(b+1)H+2D,(b+2)H+2D\},$ which is $(1)$ in Theorem \ref{cl8},
			\item [(b)] If $\OO_X(D_2)=B_0$, we have $\{H,aH+D,(a+1)H+D,(a+2)H+D,bH+2D,(b+1)H+2D,(b+2)H+2D,4H+3D\},$ which is $(2)$ in Theorem \ref{cl8}.
			\end{itemize}
		\end{itemize}
		
		\begin{table}[h!]
			\begin{center}
				\caption{Exceptional pairs.}
				\label{tab:table8}
				\begin{tabular}{c|c|c|c|c|c|c} 
					& $B_0'$ & $B_1'$ & $B_2$ & $B_3$ & $B_4$ & $B_5$\\
					\hline
					$B_0$ & $a'=a+1, a+2$ & \checkmark &  &  & \checkmark & \checkmark\\
					\hline
					$B_1$ &	 & $b'=b+1, a+2$ &  &  & \checkmark &  \checkmark\\
					\hline
					$B_2$ &	\checkmark & \checkmark &  & \checkmark &  & \checkmark\\
					\hline
					$B_3$ & \checkmark & \checkmark &  &  &  &  \\
					\hline
					$B_4$ &  &  &  &  &  & \checkmark\\
					\hline
					$B_5$ &  &  &  &  &  & \\
					\hline
				\end{tabular}
			\vspace{3mm}
			
			where $B_0=aH+D, B_1=bH+2D, B_2=H, B_3=2H, B_4=3H+3D, B_5=4H+3D$.
			\end{center}
	\end{table}
	\end{proof}
\end{thm}

\begin{thm} \label{full2}
	Let $X$ be the projective bundle $f:\PP_{\PP^2}(\OO\oplus \OO \oplus\OO(2))\rightarrow \PP^2$. Then any exceptional collection of line bundles of length 9 on $X$ is full.
	\begin{proof}
		By Lemma \ref{nor}, we only need to show that any exceptional collection of line bundles of length $9$ in Theorem \ref{cl8} is full. By Lemma \ref{mt}, it suffices to show that the exceptional collection of type $(1)$ in Theorem \ref{cl8} is full. 
		
		By Beilinson's semiorthogonal decomposition, we have 
		\begin{align}
		\D(\PP^2)=\langle \OO_{\PP^2}, \OO_{\PP^2}(1), \OO_{\PP^2}(2)\rangle. \nonumber
		\end{align}
		By Orlov's projective bundle formula (Theorem \ref{proj}), we obtain the following semiorthogonal decomposition of $\D(X)$
		\begin{align}
		\D(X)&=\langle f^*\D(\PP^2),f^*\D(\PP^2)\otimes \OO_X(1),f^*\D(\PP^2)\otimes \OO_X(2)\rangle \nonumber\\
		&=\langle f^*\OO_{\PP^2},f^*\OO_{\PP^2}(1),f^*\OO_{\PP^2}(2),f^*\OO_{\PP^2}\otimes \OO_X(1),f^*\OO_{\PP^2}(1)\otimes \OO_X(1),f^*\OO_{\PP^2}(2)\otimes \OO_X(1),\nonumber\\
		&\qquad f^*\OO_{\PP^2}\otimes \OO_X(2),f^*\OO_{\PP^2}(1)\otimes \OO_X(2),f^*\OO_{\PP^2}(2)\otimes \OO_X(2)\rangle\nonumber\\
		&=\langle \OO_X,\OO_X(H),\OO_X(2H),\OO_X(aH+D),\OO_X((a+1)H+D),\OO_X((a+2)H+D),\nonumber\\
		&\qquad \OO_X(2aH+2D),\OO_X((2a+1)H+2D),\OO_X((2a+2)H+2D)\rangle.\nonumber
		\end{align}
		Based on this semiorthogonal decomposition, we show that the exceptional collection of type $(1)$ is full by mathematical induction. For given $b,c\in \mathbb{Z}$, assume that for $b=k$ the exceptional collection
		\begin{align}\label{mi13}
		\{H, 2H, aH+D, (a+1)H+D, (a+2)H+D, kH+2D, (k+1)H+2D,(k+2)H+2D\} \tag{E7}
		\end{align}
		is full. Then for $b=k-1$, we have the exceptional collection
		\begin{align}\label{mi14} \tag{E8}
		\{H, 2H, aH+D, (a+1)H+D, (a+2)H+D, (k-1)H+2D, kH+2D,(k+1)H+2D\}
		\end{align}
		and for $b=k+1$, we have the exceptional collection
		\begin{align}\label{mi15} \tag{E9}
		\{H, 2H, aH+D, (a+1)H+D, (a+2)H+D, (k+1)H+2D, (k+2)H+2D,(k+3)H+2D\}.
		\end{align}
		Comparing the exceptional collections (\ref{mi13}) and (\ref{mi14}) with (\ref{mi15}) and by Lemma \ref{mt1}, the exceptional collections (\ref{mi14}) and (\ref{mi15}) are full. 
		Hence type $(1)$ is a full exceptional collection and so all the exceptional collections of length $9$ are full exceptional collections.
	\end{proof}
\end{thm}

\subsection{When $X\cong \PP_{\PP^2}(\OO\oplus\OO\oplus\OO(1))$.}	 
Let $f:\PP_{\PP^2}(\OO(-1)\oplus \OO \oplus\OO)\rightarrow \PP^2$ be the projective bundle. Then its canonical divisor is $K_X=-4H-3D$ where $H$ is the pullback of hyperplane class in $\PP^3$ and $D$ is the class of $c_1(\OO_X(1))$. The Picard group of $X$ is 
\begin{align}
\Pic(X)\cong \Pic(\PP^2)\oplus \mathbb{Z}[D]=\mathbb{Z}[H]\oplus\mathbb{Z}[D], \nonumber
\end{align}
with intersection numbers
\begin{align}
H^4=0, ~H^3D=0, ~H^2D^2=1, ~HD^3=-1, ~D^4=1. \nonumber
\end{align}
\subsubsection{Cohomologically zero line bundles}
\begin{lem}\label{CZ9}
	$H^0(X,aH+bD)=0$ if and only if $b<0$ or $a+b<0, ~b\geq 0$. Also, $H^4(X,aH+bD)=0$ if and only if $b>-3$ or $a+b>-7, ~b\leq -3$.
	\begin{proof}
		Note that we have the following isomorphisms:
		\begin{align}
		f_*(\OO_X(b))&\cong \Sym^b(\OO_{\PP^2}\oplus \OO_{\PP^2} \oplus\OO_{\PP^2}(1)) \nonumber \\ 
		&\cong \bigoplus_{i=0}^b \Sym^i(\OO_{\PP^2}\oplus \OO_{\PP^2})\otimes \Sym^{b-i}(\OO_{\PP^2}(1)) \nonumber \\
		&\cong \bigoplus_{i=0}^b \OO_{\PP^2}(b-i) ~\text{if}~ b\geq 0, \nonumber
		\end{align} and $f_*(\OO_X(b))\cong 0$ if $b<0$.
		
		By projection formula, we have the following isomorphisms:
		\begin{align}
		H^0(X, f^*(\OO_{\PP^2}(a))\otimes \OO_X(b))&\cong H^0(\PP^2, \OO_{\PP^2}(a)\otimes  \bigoplus_{i=0}^b \OO_{\PP^2}(2b-2i)) \nonumber \\
		&\cong \bigoplus_{i=0}^b H^0(\PP^2,\OO_{\PP^2}(a+b-i)) ~\text{if}~ b\geq 0,\nonumber
		\end{align}
		and if $b<0$, then $H^0(X,aH+bD)=0$.
		Hence, $H^0(X,aH+bD)=0$ if and only if $b<0$ or $a+b<0, ~b\geq 0$.
		
		By Serre duality, one obtains the condition for $H^4(X,aH+bD)=0$.
	\end{proof}
\end{lem}
\begin{pro} \label{pro9}
	A line bundle $\OO_X(aH+bD)$ is cohomologically zero if one of the following holds:
	\begin{itemize}
		\item [(1)] $b=-1$,
		\item [(2)] $b=-2$,
		\item [(3)] $a=-1, b=0$,
		\item [(4)] $a=-2, b=0$,
		\item [(5)] $a=-2, b=-3$,
		\item [(6)] $a=-3, b=-3$.
	\end{itemize}
	\begin{proof}
		Since $c(X)=(1+3H+3H^2)\{(1+D)^3+H(1+D)^2\}$, we have $c_1=4H+3D$, $c_2=6H^2+11HD+3D^2$, $c_3=3H^3+15H^2D+10HD^2+D^3$ and $c_4=6H^3D+12H^2D^2+3HD^3=9$.
		Thus, by Lemma \ref{HRR}, we have
		\begin{align}
		\chi(\OO_X(aH+bD))&=\frac{1}{24}(24+36a+26b+12a^2+46ab-b^2 \nonumber\\
		&\qquad +18a^2b+6ab^2-2b^3+6a^2b^2-4ab^3+b^4) \nonumber\\
		&=\frac{1}{24}(b+1)(b+2)(6a^2-4ab+18a+b^2-5b+12). \nonumber
		\end{align}
		Therefore, by Lemma \ref{CZ9}, $\OO_X(aH+bD)$ is cohomologically zero if the following conditions hold:
		\begin{itemize}
			\item [$\mathrm{(i)}$] $b<0$ or $a+b<0, ~b\geq 0$,
			\item [$\mathrm{(ii)}$] $b>-3$ or $a+b>-7, ~b\leq -3$,
			\item [$\mathrm{(iii)}$] $(b+1)(b+2)(6a^2-4ab+18a+b^2-5b+12)=0$.
		\end{itemize}
		Note that the equation $6a^2-4ab+18a+b^2-5b+12=0$ has $4$ integer solutions; $(a,b)=(-1,0),(-2,0),(-2,-3),(-3,-3)$ under conditions $\mathrm{(i)}$ and $\mathrm{(ii)}$. This completes the proof.
	\end{proof}
\end{pro}

\subsubsection{Classification of full exceptional collections}
\begin{thm} \label{cl9}
	Let $X$ be the projective bundle $f:\PP_{\PP^2}(\OO\oplus \OO \oplus\OO(1))\rightarrow \PP^2$. Then the normalized sequence
	\begin{align}
	\{\OO_X,\OO_X(D_1),\cdots,\OO_X(D_8)\} \label{9} \tag{I}
	\end{align}
	is an exceptional collection of line bundles if and only if the ordered set of divisors $\{D_1,D_2,\cdots,D_8\}$ is one of the following types:
	\begin{itemize}
		\item [(1)] $\{H,2H,aH+D,(a+1)H+D,(a+2)H+D,bH+2D,(b+1)H+2D,(b+2)H+2D\}$,
		\item [(2)]
		$\{H,aH+D,(a+1)H+D,(a+2)H+D,bH+2D,(b+1)H+2D,(b+2)H+2D,3H+3D\}$
		\item [(3)] $\{aH+D,(a+1)H+D,(a+2)H+D,bH+2D,(b+1)H+2D,(b+2)H+2D,2H+3D,3H+3D\}$,
	\end{itemize}
	where $a,b\in \mathbb{Z}$. Moreover, by mutations and normalizations, they are related as:
	\begin{align}
	(1)\rightarrow (2)\rightarrow (3)\rightarrow (1). \nonumber
	\end{align}
	\begin{proof}
		Write $D_0=0$. The sequence (\ref{9}) is an exceptional collection if and only if for any integers $0\leq i<j\leq 8$, the line bundles $\OO_X(D_i-D_j)$ are cohomologically zero.
		
		To make the sequence (\ref{9}) an exceptional collection, $\OO_X(D_i)$ must be one of the line bundles in Proposition \ref{pro9}. To find out all the possible exceptional collections, we need Table \ref{tab:table9}.
		
		From Table \ref{tab:table9}, we observe that $\OO_X(D_1)$ may be one of $B_0$ and $B_3$. To finish the proof, we shall discuss case by case:
		\begin{itemize}
			\item [Case1.] Suppose $\OO_X(D_1)=B_0$. In this case, the only possible combination is $\{aH+D,(a+1)H+D,(a+2)H+D,bH+2D,(b+1)H+2D,(b+2)H+2D,2H+3D,3H+3D\}$, which is $(3)$ in Theorem \ref{cl9}.
			\item [Case2.] Suppose $\OO_X(D_1)=B_3$. In this case, there are two possible choices for $D_2$:
			\begin{itemize}
				\item [(a)]
				If $\OO_X(D_2)=B_4$, we have $\{H,2H,aH+D,(a+1)H+D,(a+2)H+D,bH+2D,(b+1)H+2D,(b+2)H+2D\},$ which is $(1)$ in Theorem \ref{cl9},
				\item [(b)] If $\OO_X(D_2)=B_0$, we have $\{H,aH+D,(a+1)H+D,(a+2)H+D,bH+2D,(b+1)H+2D,(b+2)H+2D,3H+3D\},$ which is $(2)$ in Theorem \ref{cl9}.
			\end{itemize}
		\end{itemize}
		
		\begin{table}[t]
			\begin{center}
				\caption{Exceptional pairs.}
				\label{tab:table9}
				\begin{tabular}{c|c|c|c|c|c|c} 
					& $B_0'$ & $B_1'$ & $B_2$ & $B_3$ & $B_4$ & $B_5$\\
					\hline
					$B_0$ & $a'=a+1, a+2$ & \checkmark &  &  & \checkmark & \checkmark\\
					\hline
					$B_1$ &	 & $a'=a+1, a+2$ &  &  & \checkmark &  \checkmark\\
					\hline
					$B_2$ &	\checkmark & \checkmark &  & \checkmark &  & \checkmark \\
					\hline
					$B_3$ & \checkmark & \checkmark &  &  &  &  \\
					\hline
					$B_4$ &  &  &  &  &  & \checkmark\\
					\hline
					$B_5$ &  &  &  &  &  & \\
					\hline
				\end{tabular}
			\vspace{3mm}
			
			where $B_0=aH+D, B_1=bH+2D, B_2=H, B_3=2H, B_4=2H+3D, B_5=3H+3D$.
			\end{center}
		\end{table}
	\end{proof}
\end{thm}

\begin{thm}
	Let $X$ be the projective bundle $f:\PP_{\PP^2}(\OO\oplus \OO \oplus\OO(1))\rightarrow \PP^2$. Then any exceptional collection of line bundles of length 9 on $X$ is full.
	\begin{proof}
		By Lemma \ref{nor}, we only need to show that any exceptional collection of line bundles of length $9$ in Theorem \ref{cl9} is full. By Lemma \ref{mt}, it suffices to show that the exceptional collection of type $(1)$ in Theorem \ref{cl9} is full. 
		
		By Beilinson's semiorthogonal decomposition, we have 
		\begin{align}
		\D(\PP^2)=\langle \OO_{\PP^2}, \OO_{\PP^2}(1), \OO_{\PP^2}(2)\rangle. \nonumber
		\end{align}
		By Orlov's projective bundle formula (Theorem \ref{proj}), we obtain the following semiorthogonal decomposition of $\D(X)$
		\begin{align}
		\D(X)&=\langle f^*\D(\PP^2),f^*\D(\PP^2)\otimes \OO_X(1),f^*\D(\PP^2)\otimes \OO_X(2)\rangle \nonumber\\
		&=\langle f^*\OO_{\PP^2},f^*\OO_{\PP^2}(1),f^*\OO_{\PP^2}(2),f^*\OO_{\PP^2}\otimes \OO_X(1),f^*\OO_{\PP^2}(1)\otimes \OO_X(1),f^*\OO_{\PP^2}(2)\otimes \OO_X(1),\nonumber\\
		&\qquad f^*\OO_{\PP^2}\otimes \OO_X(2),f^*\OO_{\PP^2}(1)\otimes \OO_X(2),f^*\OO_{\PP^2}(2)\otimes \OO_X(2)\rangle \nonumber\\
		&=\langle \OO_X,\OO_X(H),\OO_X(2H),\OO_X(aH+D),\OO_X((a+1)H+D),\OO_X((a+2)H+D),\nonumber\\
		&\qquad \OO_X(2aH+2D),\OO_X((2a+1)H+2D),\OO_X((2a+2)H+2D)\rangle.\nonumber
		\end{align}
		As in the proof of Theorem \ref{full2}, we can conclude that the type $(1)$ is a full exceptional collection and so all the exceptional collections of length $9$ are full exceptional collections.
	\end{proof}
\end{thm}

\subsection{When $X\cong \PP_{\PP^2}(\OO\oplus\OO(1)\oplus\OO(1))$.}	
Let $f:\PP_{\PP^2}(\OO(-1)\oplus \OO \oplus\OO)\rightarrow \PP^2$ be the projective bundle. Then its canonical divisor is $K_X=-2H-3D$ where $H$ is the pullback of hyperplane class in $\PP^2$ and $D$ is the class of $c_1(\OO_X(1))$. The Picard group of $X$ is 
\begin{align}
\Pic(X)\cong \Pic(\PP^2)\oplus \mathbb{Z}[D]=\mathbb{Z}[H]\oplus\mathbb{Z}[D], \nonumber
\end{align}
with intersection numbers
\begin{align}
H^4=0, ~H^3D=0, ~H^2D^2=1, ~HD^3=1, ~D^4=1. \nonumber
\end{align}
\subsubsection{Cohomologically zero line bundles}
\begin{lem}\label{CZ10}
	$H^2(X,aH+bD)=0$ if either one of the following holds:
	\begin{itemize}
		\item [(1)] $-3<b<0$,
		\item [(2)] $a>b-3$ and $b\geq 0$,
		\item [(3)] $a<b+4$ and $b\leq -3$.
	\end{itemize}
	\begin{proof}
		Note that we have the following isomorphisms:
		\begin{align}
		f_*(\OO_X(b))&\cong \Sym^b(\OO_{\PP^2}(-1)\oplus \OO_{\PP^2} \oplus\OO_{\PP^2}) \nonumber \\ 
		&\cong \bigoplus_{i=0}^b \Sym^i(\OO_{\PP^2}(-1))\otimes \Sym^{b-i}(\OO_{\PP^2}\oplus \OO_{\PP^2}) \nonumber \\
		&\cong \bigoplus_{i=0}^b \OO_{\PP^2}(-i) ~\text{if}~ b\geq 0, \nonumber
		\end{align} and $f_*(\OO_X(b))\cong 0$ if $b<0$.
		
		By projection formula and Leray's spectral sequence, we have the following isomorphisms:
		\begin{align}
		H^2(X, f^*(\OO_{\PP^2}(a))\otimes \OO_X(b))&\cong H^0(\PP^2, \OO_{\PP^2}(a)\otimes R^2f_*(\OO_X(b))) \oplus H^2(\PP^2,\OO_{\PP^2}(a)\otimes f_*(\OO_X(b))) \nonumber \\
		&\cong H^2(\PP^2,\OO_{\PP^2}(a)\otimes f_*(\OO_X(b))) ~\text{if}~ b>-3 \nonumber\\
		&\cong \bigoplus_{i=0}^b H^2(\PP^2,\OO_{\PP^2}(a-i)) ~\text{if}~ b\geq 0\nonumber\\
		&\cong \bigoplus_{i=0}^b H^0(\PP^2,\OO_{\PP^2}(-3-a+i)) ~\text{by Serre duality}. \nonumber
		\end{align}
		Thus, $H^2(X,aH+bD)=0$ if $b\geq 0$, $a>b-3$ or $-3<b<0$ and by using Serre duality, we obtain $H^2(X,aH+bD)=0$ if $a<b+4$ and $b\leq -3$.
	\end{proof}
\end{lem}

\begin{pro} \label{pro10}
	A line bundle $\OO_X(aH+bD)$ is cohomologically zero if one of the following holds:
	\begin{itemize}
		\item [(1)] $b=-1$,
		\item [(2)] $b=-2$,
		\item [(3)] $a=-1, b=0$,
		\item [(4)] $a=-2, b=0$,
		\item [(5)] $a=0, b=-3$,
		\item [(6)] $a=-1, b=-3$.
	\end{itemize}
	\begin{proof}
		Since $c(X)=(1+3H+3H^2)\{(1+D)^3-H(1+D)^2\}$, we have $c_1=2H+3D$, $c_2=7HD+3D^2$, $c_3=-3H^3+3H^2D+8HD^2+D^3$ and $c_4=-6H^3D+6H^2D^2+3HD^3=9$.
		Thus, by Lemma \ref{HRR}, we have
		\begin{align}
		\chi(\OO_X(aH+bD))&=\frac{1}{24}(24+36a+50b+12a^2+62ab+35^2 \nonumber\\
		&\qquad +18a^2b+30ab^2+10b^3+6a^2b^2+4ab^3+b^4) \nonumber\\
		&=\frac{1}{24}(b+1)(b+2)(6a^2+4ab+18a+b^2+7b+12). \nonumber
		\end{align}
		Therefore, by Lemma \ref{CZ10}, $\OO_X(aH+bD)$ is cohomologically zero if the following conditions hold:
		\begin{itemize}
			\item [$\mathrm{(i)}$] $-3<b<0$ or $a>b-3$, $b\geq 0$ or $a<b+4$, $b\leq -3$,
			\item [$\mathrm{(ii)}$] $(b+1)(b+2)(6a^2+4ab+18a+b^2+7b+12)=0$.
		\end{itemize}
		Note that the equation $6a^2+4ab+18a+b^2+7b+12=0$ has $4$ integer solutions; $(a,b)=(-1,0),(-2,0),(0,-3),(-1,-3)$ under condition $\mathrm{(i)}$. This completes the proof.
	\end{proof}
\end{pro}

\subsubsection{Classification of full exceptional collections}
\begin{thm} \label{cl10}
	Let $X$ be the projective bundle $f:\PP_{\PP^2}(\OO\oplus \OO \oplus\OO(1))\rightarrow \PP^2$. Then the normalized sequence
	\begin{align}
	\{\OO_X,\OO_X(D_1),\cdots,\OO_X(D_8)\} \label{10} \tag{J}
	\end{align}
	is an exceptional collection of line bundles if and only if the ordered set of divisors $\{D_1,D_2,\cdots,D_8\}$ is one of the following types:
	\begin{itemize}
		\item [(1)] $\{H,2H,aH+D,(a+1)H+D,(a+2)H+D,bH+2D,(b+1)H+2D,(b+2)H+2D\}$,
		\item [(2)]
		$\{H,aH+D,(a+1)H+D,(a+2)H+D,bH+2D,(b+1)H+2D,(b+2)H+2D,H+3D\}$
		\item [(3)] $\{aH+D,(a+1)H+D,(a+2)H+D,bH+2D,(b+1)H+2D,(b+2)H+2D,3D,H+3D\}$,
	\end{itemize}
	where $a,b\in \mathbb{Z}$. Moreover, by mutations and normalizations, they are related as:
	\begin{align}
	(1)\rightarrow (2)\rightarrow (3)\rightarrow (1). \nonumber
	\end{align}
	\begin{proof}
		Write $D_0=0$. The sequence (\ref{10}) is an exceptional collection if and only if for any integers $0\leq i<j\leq 8$, the line bundles $\OO_X(D_i-D_j)$ are cohomologically zero.
		
		To make the sequence (\ref{10}) an exceptional collection, $\OO_X(D_i)$ must be one of the line bundles in Proposition \ref{pro10}. To find out all the possible exceptional collections, we need Table \ref{tab:table10}.
		
		From Table \ref{tab:table10}, we observe that $\OO_X(D_1)$ may be one of $B_0$ and $B_3$. To finish the proof, we shall discuss case by case:
		\begin{itemize}
			\item [Case1.] Suppose $\OO_X(D_1)=B_0$. In this case, the only possible combination is $\{aH+D,(a+1)H+D,(a+2)H+D,bH+2D,(b+1)H+2D,(b+2)H+2D,3D,H+3D\}$, which is $(3)$ in Theorem \ref{cl10}.
			\item [Case2.] Suppose $\OO_X(D_1)=B_2$. In this case, there are two possible choices for $D_2$:
			\begin{itemize}
				\item [(a)]
				If $\OO_X(D_2)=B_3$, we have $\{H,2H,aH+D,(a+1)H+D,(a+2)H+D,bH+2D,(b+1)H+2D,(b+2)H+2D\},$ which is $(1)$ in Theorem \ref{cl10},
				\item [(b)] If $\OO_X(D_2)=B_0$, we have $\{H,aH+D,(a+1)H+D,(a+2)H+D,bH+2D,(b+1)H+2D,(b+2)H+2D,H+3D\},$ which is $(2)$ in Theorem \ref{cl10}.
			\end{itemize}
		\end{itemize}
		
	\begin{table}[h!]
		\begin{center}
			\caption{Exceptional pairs.}
			\label{tab:table10}
			\begin{tabular}{c|c|c|c|c|c|c} 
				& $B_0'$ & $B_1'$ & $B_2$ & $B_3$ & $B_4$ & $B_5$\\
				\hline
				$B_0$ & $a'=a+1, a+2$ & \checkmark &  &  & \checkmark & \checkmark\\
				\hline
				$B_1$ &	 & $a'=a+1, a+2$ &  &  & \checkmark &  \checkmark\\
				\hline
				$B_2$ &	\checkmark & \checkmark &  & \checkmark &  & \checkmark \\
				\hline
				$B_3$ & \checkmark & \checkmark &  &  &  &  \\
				\hline
				$B_4$ &  &  &  &  &  & \checkmark\\
				\hline
				$B_5$ &  &  &  &  &  & \\
				\hline
			\end{tabular}
		\vspace{3mm}
		
		where $B_0=aH+D, B_1=bH+2D, B_2=H, B_3=2H, B_4=3D, B_5=H+3D$.
		\end{center}
	\end{table}
	\end{proof}
\end{thm}

\begin{thm}
	Let $X$ be the projective bundle $f:\PP_{\PP^2}(\OO(-1)\oplus \OO \oplus\OO)\rightarrow \PP^2$. Then any exceptional collection of line bundles of length 9 on $X$ is full.
	\begin{proof}
		By Lemma \ref{nor}, we only need to show that any exceptional collection of line bundles of length $9$ in Theorem \ref{cl10} is full. By Lemma \ref{mt}, it suffices to show that the exceptional collection of type $(1)$ in Theorem \ref{cl10} is full. 
		
		By Beilinson's semiorthogonal decomposition, we have 
		\begin{align}
		\D(\PP^2)=\langle \OO_{\PP^2}, \OO_{\PP^2}(1), \OO_{\PP^2}(2)\rangle. \nonumber
		\end{align}
		By Orlov's projective bundle formula (Theorem \ref{proj}), we obtain the following semiorthogonal decomposition of $\D(X)$
		\begin{align}
		\D(X)&=\langle f^*\D(\PP^2),f^*\D(\PP^2)\otimes \OO_X(1),f^*\D(\PP^2)\otimes \OO_X(2)\rangle \nonumber\\
		&=\langle f^*\OO_{\PP^2},f^*\OO_{\PP^2}(1),f^*\OO_{\PP^2}(2),f^*\OO_{\PP^2}\otimes \OO_X(1),f^*\OO_{\PP^2}(1)\otimes \OO_X(1),f^*\OO_{\PP^2}(2)\otimes \OO_X(1),\nonumber\\
		&\qquad f^*\OO_{\PP^2}\otimes \OO_X(2),f^*\OO_{\PP^2}(1)\otimes \OO_X(2),f^*\OO_{\PP^2}(2)\otimes \OO_X(2)\rangle\nonumber\\
		&=\langle \OO_X,\OO_X(H),\OO_X(2H),\OO_X(aH+D),\OO_X((a+1)H+D),\OO_X((a+2)H+D),\nonumber\\
		&\qquad \OO_X(2aH+2D),\OO_X((2a+1)H+2D),\OO_X((2a+2)H+2D)\rangle.\nonumber
		\end{align}
		As in the proof of Theorem \ref{full2}, we can conclude that the type $(1)$ is a full exceptional collection and so all the exceptional collections of length $9$ are full exceptional collections.
	\end{proof}
\end{thm}

\subsection{When $X\cong \PP^2 \times \PP^2$.}	
Let $p_1: \PP^2 \times \PP^2\rightarrow \PP^2$ and $p_2: \PP^2 \times \PP^2\rightarrow \PP^2$ be the projections onto the first and the second factor, respectively. 
Then the canonical divisor of $X$ is $K_X=-3H-3D$, where $H$ and $D$ be the pullback of the hyperplane class of $\PP^2$ of the first factor and $\PP^2$ of the second factor, respectively.

\subsubsection{Classification of full exceptional collections}
\begin{thm} \label{cl11}
	Let $X$ be $\PP^2 \times \PP^2$. Then the normalized sequence
	\begin{align}
	\{\OO_X,\OO_X(D_1),\cdots,\OO_X(D_8)\} \label{11} \tag{K}
	\end{align}
	is an exceptional collection of line bundles if and only if the ordered set of divisors $\{D_1,D_2,\cdots,D_8\}$ is one of the following types:
	\begin{itemize}
		\item [(1)] $\{H+aD,H+(a+1)D,H+(a+2)D,2H+bD,2H+(b+1)D,2H+(b+2)D,3H+D,3H+2D\}$,
		\item [(2)] $\{D,2D,H+aD,H+(a+1)D,H+(a+2)D,2H+bD,2H+(b+1)D,2H+(b+2)D\}$,
		\item [(3)] $\{D,H+aD,H+(a+1)D,H+(a+2)D,2H+bD,2H+(b+1)D,2H+(b+2)D,3H+2D\}$,
		
		\item [(4)] $\{H+aD,H+(a+1)D,H+(a+2)D,2H,2H+D,3H+D,2H+2D,3H+2D\}$,
		\item [(5)]  $\{D, 2D, H+aD, H+(a+1)D, 2H+(a+1)D, H+(a+2)D, 2H+(a+2)D,2H+(a+3)D\}$,
		\item [(6)]  $\{D, H+aD, H+(a+1)D, 2H+(a+1)D, H+(a+2)D, 2H+(a+2)D,2H+(a+3)D, 3H+2D\}$,
		\item [(7)]  $\{H+aD, H+(a+1)D, 2H+(a+1)D, H+(a+2)D, 2H+(a+2)D,2H+(a+3)D, 3H+D,3H+2D\}$,
		\item [(8)]  $\{D, H+D, 2D, H+2D, H+3D, 2H+bD, 2H+(b+1)D,2H+(b+2)D\}$,
		\item [(9)]  $\{H, D, H+D, H+2D, 2H+bD, 2H+(b+1)D,2H+(b+2)D,3H+2D\}$,
		\item [(10)] $\{-H+D, D, 2D, H+aD, H+(a+1)D, H+(a+2)D, 2H+2D,2H+3D\}$,
		\item [(11)] $\{H, H+D, 2H+bD, 2H+(b+1)D, 2H+(b+2)D, 3H+D, 3H+2D,4H+2D\}$,
		\item [(12)] $\{D, H+aD, H+(a+1)D, H+(a+2)D, 2H+D, 2H+2D, 3H+2D,2H+3D\}$,
		
		\item [(13)] $\{H+aD,H+(a+1)D,H+(a+2)D, 2H+D, 2H+2D, 3H+D, 2H+3D,3H+2D\}$,
		\item [(14)] $\{D, 2D, H+aD, H+(a+1)D, 2H+aD, H+(a+2)D, 2H+(a+1)D,2H+(a+2)D\}$,
		\item [(15)] $\{D, H+aD, H+(a+1)D, 2H+aD, H+(a+2)D, 2H+(a+1)D,2H+(a+2)D,3H+2D\}$,
		\item [(16)] $\{H+aD, H+(a+1)D, 2H+aD, H+(a+2)D, 2H+(a+1)D,2H+(a+2)D, 3H+D, 3H+2D\}$,
		\item [(17)] $\{D, H, 2D, H+D, H+2D, 2H+bD, 2H+(b+1)D,2H+(b+2)D\}$,
		\item [(18)] $\{H-D, D, H, H+D, 2H+bD, 2H+(b+1)D, 2H+(b+2)D,3H+2D\}$,
		\item [(19)] $\{-H+2D, D, 2D, H+aD, H+(a+1)D, H+(a+2)D, 2H+3D,2H+4D\}$,
		\item [(20)] $\{H-D, H, 2H+bD, 2H+(b+1)D, 2H+(b+2)D, 3H+D, 3H+2D,4H+D\}$,
		\item [(21)] $\{D, H+aD, H+(a+1)D, H+(a+2)D, 2H+2D, 2H+3D, 3H+2D,2H+4D\}$,
		
		\item [(22)] $\{H+aD,H+(a+1)D,H+(a+2)D, 2H+D, 3H+D, 2H+2D, 2H+3D,3H+2D\}$,
		\item [(23)] $\{D, 2D, H+aD, 2H+aD, H+(a+1)D, H+(a+2)D, 2H+(a+1)D,2H+(a+2)D\}$,
		\item [(24)] $\{D, H+aD, 2H+aD, H+(a+1)D, H+(a+2)D, 2H+(a+1)D,2H+(a+2)D, 3H+2D\}$,
		\item [(25)] $\{H+aD, 2H+aD, H+(a+1)D, H+(a+2)D, 2H+(a+1)D,2H+(a+2)D, 3H+D, 3H+2D\}$,
		\item [(26)] $\{H, D, 2D, H+D, H+2D, 2H+bD, 2H+(b+1)D,2H+(b+2)D\}$,
		\item [(27)] $\{-H+D, -H+2D, D, 2D, H+aD, H+(a+1)D, H+(a+2)D,2H+3D\}$,
		\item [(28)] $\{D, H, H+D, 2H+bD, 2H+(b+1)D, 2H+(b+2)D,3H+2D,4H+2D\}$,
		\item [(29)] $\{H-D, H, 2H+bD, 2H+(b+1)D, 2H+(b+2)D, 3H+D, 4H+D,3H+2D\}$,
		\item [(30)] $\{D, H+aD, H+(a+1)D, H+(a+2)D, 2H+2D, 3H+2D, 2H+3D,2H+4D\}$,
		
		\item [(31)] $\{H+aD,H+(a+1)D,H+(a+2)D, 2H+D, 3H+D, 2H+2D, 3H+2D,2H+3D\}$
		\item [(32)] $\{D, 2D, H+aD, 2H+aD, H+(a+1)D, 2H+(a+1)D, H+(a+2)D,2H+(a+2)D\}$,
		\item [(33)] $\{D, H+aD, 2H+aD, H+(a+1)D, 2H+(a+1)D, H+(a+2)D,2H+(a+2)D, 3H+2D\}$,
		\item [(34)] $\{H+aD, 2H+aD, H+(a+1)D, 2H+(a+1)D, H+(a+2)D,2H+(a+2)D, 3H+D, 3H+2D\}$,
		\item [(35)] $\{H, D, H+D, 2D, H+2D, 2H+bD, 2H+(b+1)D,2H+(b+2)D\}$,
		\item [(36)] $\{-H+D, D, -H+2D, 2D, H+aD, H+(a+1)D, H+(a+2)D,2H+3D\}$,
		\item [(37)] $\{H, D, H+D, 2H+bD, 2H+(b+1)D, 2H+(b+2)D, 3H+2D,4H+2D\}$,
		\item [(38)] $\{-H+D, D, H+aD, H+(a+1)D, H+(a+2)D, 2H+2D, 3H+2D,2H+3D\}$,
		\item [(39)] $\{H, 2H+bD, 2H+(b+1)D, 2H+(b+2)D, 3H+D, 4H+D, 3H+2D,4H+2D\}$,
		
		\item [(40)] $\{H+aD, H+(a+1)D, 2H+aD, 2H+(a+1)D, H+(a+2)D, 2H+(a+2)D, 3H+D,3H+2D\}$,
		\item [(41)] $\{D, H, H+D, 2D, H+2D, 2H+bD, 2H+(b+1)D,2H+(b+2)D\}$,
		\item [(42)] $\{H-D, H, D, H+D, 2H+bD, 2H+(b+1)D, 2H+(b+2)D,3H+2D\}$,
		\item [(43)] $\{D, -H+2D, 2D, H+aD, H+(a+1)D, H+(a+2)D, 2H+3D,2H+4D\}$,
		\item [(44)] $\{-H+D, D, H+aD, H+(a+1)D, H+(a+2)D, 2H+2D, 2H+3D,3H+2D\}$,
		\item [(45)] $\{H, 2H+bD, 2H+(b+1)D, 2H+(b+2)D, 3H+D, 3H+2D, 4H+D,4H+2D\}$,
		\item [(46)] $\{H+aD, H+(a+1)D, H+(a+2)D, 2H+D, 2H+2D, 3H+D, 3H+2D,2H+3D\}$,
		\item [(47)] $\{D, 2D, H+aD, H+(a+1)D, 2H+aD, 2H+(a+1)D, H+(a+2)D,2H+(b+2D)\}$,
		\item [(48)] $\{D, H+aD, H+(a+1)D, 2H+aD, 2H+(a+1)D, H+(a+2)D,2H+(b+2D), 3H+2D\}$,
		
		\item [(49)] $\{H-D, H, 2H-D, 2H, cH+D, (c+1)H+D, (c+2)H+D,3H+2D\}$,
		\item [(50)] $\{D, H, H+D, dH+2D, (d+1)H+2D, (d+2)H+2D, 2H+3D,2H+4D\}$,
		\item [(51)] $\{H-D, H, cH+D, (c+1)H+D, (c+2)H+D, 2H+2D, 2H+3D,3H+2D\}$,
		\item [(52)] $\{D, dH+2D, (d+1)H+2D, (d+2)H+2D, H+3D, H+4D, 2H+3D,2H+4D\}$,
		\item [(53)] $\{cH+D, (c+1)H+D, (c+2)H+D, H+2D, H+3D, 2H+2D, 2H+3D,3H+2D\}$,
		\item [(54)] $\{H, 2H, cH+D, cH+2D, (c+1)H+D, (c+1)H+2D, (c+2)H+D,(c+2)H+2D\}$,
		\item [(55)] $\{H, cH+D, cH+2D, (c+1)H+D, (c+1)H+2D, (c+2)H+D,(c+2)H+2D, 2H+3D\}$,
		\item [(56)] $\{cH+D, cH+2D, (c+1)H+D, (c+1)H+2D, (c+2)H+D,(c+2)H+2D, H+3D, 2H+3D\}$,
		\item [(57)] $\{D, H, H+D, 2H, 2H+D, dH+2D, (d+1)H+2D,(d+2)H+2D\}$,
		
		\item [(58)] $\{H, H+D, 2H, 2H+D, 3H+D, dH+2D, (d+1)H+2D,(d+2)H+2D\}$,
		\item [(59)] $\{D, H, H+D, 2H+D, dH+2D, (d+1)H+2D,(d+2)H+2D, 2H+3D\}$,
		\item [(60)] $\{H-D, H, 2H, cH+D, (c+1)H+D,(c+2)H+D, 2H+2D, 3H+2D, 2H+4D\}$,
		\item [(61)] $\{D, H+D, dH+2D, (d+1)H+2D,(d+2)H+2D, H+3D, 2H+3D,2H+4D\}$,
		\item [(62)] $\{H, cH+D, (c+1)H+D,(c+2)H+D, H+2D, 2H+2D, 2H+3D,3H+2D\}$,
		\item [(63)] $\{cH+D, (c+1)H+D,(c+2)H+D, 2D, H+2D, H+3D, 2H+2D,2H+3D\}$,
		\item [(64)] $\{H, 2H, cH+D, (c+1)H+D, (c+1)H+2D, (c+1)H+D, (c+2)H+2D,(c+3)H+2D\}$,
		\item [(65)] $\{H, cH+D, (c+1)H+D, (c+1)H+2D, (c+1)H+D, (c+2)H+2D,(c+3)H+2D, 2H+3D\}$,
		\item [(66)] $\{cH+D, (c+1)H+D, (c+1)H+2D, (c+1)H+D, (c+2)H+2D,(c+3)H+2D, H+3D, 2H+3D\}$,
		
		\item [(67)] $\{H-D, H, cH+D, (c+1)H+D,(c+2)H+D, 2H+2D, 3H+2D,2H+3D\}$,
		\item [(68)] $\{D, dH+2D, (d+1)H+2D,(d+2)H+2D, H+3D,2H+3D,H+4D,2H+4D\}$,
		\item [(69)] $\{cH+D, (c+1)H+D,(c+2)H+D, H+2D, 2H+2D, H+3D,2H+3D,3H+2D\}$,
		\item [(70)] $\{H, 2H, cH+D, (c+1)H+D,cH+2D, (c+1)H+2D, (c+2)H+D,(c+2)H+2D\}$,
		\item [(71)] $\{H, cH+D, (c+1)H+D,cH+2D, (c+1)H+2D, (c+2)H+D,(c+2)H+2D,2H+3D\}$,
		\item [(72)] $\{cH+D, (c+1)H+D,cH+2D, (c+1)H+2D, (c+2)H+D,(c+2)H+2D, H+3D,2H+3D\}$,
		\item [(73)] $\{H, D, H+D, 2H, 2H+D, dH+2D, (d+1)H+2D,(d+2)H+2D\}$,
		\item [(74)] $\{-H+D, D, H, H+D, dH+2D, (d+1)H+2D,(d+2)H+2D, 2H+3D\}$,
		\item [(75)] $\{H, 2H-D, 2H, cH+D, (c+1)H+D,(c+2)H+D, 3H+2D, 4H+2D\}$,
		
		\item [(76)] $\{H, 2H, cH+D, (c+1)H+D,(c+2)H+D, dH+2D, (d+1)H+2D,(d+2)H+2D\}$,
		\item [(77)] $\{H, cH+D, (c+1)H+D,(c+2)H+D, dH+2D, (d+1)H+2D,(d+2)H+2D,2H+3D\}$,
		\item [(78)] $\{cH+D, (c+1)H+D,(c+2)H+D, dH+2D, (d+1)H+2D,(d+2)H+2D,H+3D,2H+3D\}$,
		
		\item [(79)] $\{H, 2H, cH+D, (c+1)H+D,cH+2D, (c+2)H+D, (c+1)H+2D,(c+2)H+2D\}$,
		\item [(80)] $\{H, cH+D, (c+1)H+D,cH+2D, (c+2)H+D, (c+1)H+2D,(c+2)H+2D,2H+3D\}$,
		\item [(81)] $\{cH+D, (c+1)H+D,cH+2D, (c+2)H+D, (c+1)H+2D,(c+2)H+2D, H+3D,2H+3D\}$,
		\item [(82)] $\{H, D, 2H, H+D, 2H+D, dH+2D, (d+1)H+2D,(d+2)H+2D\}$,
		\item [(83)] $\{-H+D, H, D, H+D, dH+2D, (d+1)H+2D,(d+2)H+2D, 2H+3D\}$,
		\item [(84)] $\{2H-D, H, 2H, cH+D, (c+1)H+D,(c+2)H+D, 3H+2D, 4H+2D\}$,
		\item [(85)] $\{-H+D, D, dH+2D, (d+1)H+2D, (d+2)H+2D, H+3D,2H+3D,H+4D\}$,
		\item [(86)] $\{H, cH+D, (c+1)H+D,(c+2)H+D, 2H+2D,3H+2D,2H+3D,4H+2D\}$,
		\item [(87)] $\{cH+D, (c+1)H+D,(c+2)H+D, H+2D, 2H+2D, H+3D, 3H+2D, 2H+3D\}$,
		
		\item [(88)] $\{H, 2H, cH+D, cH+2D,(c+1)H+D, (c+2)H+D, (c+1)H+2D,(c+2)H+2D\}$,
		\item [(89)] $\{H, cH+D, cH+2D,(c+1)H+D, (c+2)H+D, (c+1)H+2D,(c+2)H+2D,2H+3D\}$,
		\item [(90)] $\{cH+D, cH+2D,(c+1)H+D, (c+2)H+D, (c+1)H+2D,(c+2)H+2D, H+3D,2H+3D\}$,
		\item [(91)] $\{D, H, 2H, H+D, 2H+D, dH+2D, (d+1)H+2D,(d+2)H+2D\}$,
		\item [(92)] $\{H-D, 2H-D, H, 2H, cH+D, (c+1)H+D,(c+2)H+D, 3H+2D\}$,
		\item [(93)] $\{H, D, H+D, dH+2D, (d+1)H+2D,(d+2)H+2D, 2H+3D, 2H+4D\}$,
		\item [(94)] $\{-H+D, D, dH+2D, (d+1)H+2D, (d+2)H+2D, H+3D,H+4D,2H+3D\}$,
		\item [(95)] $\{H, cH+D, (c+1)H+D,(c+2)H+D, 2H+2D,2H+3D,3H+2D,4H+2D\}$,
		\item [(96)] $\{cH+D, (c+1)H+D,(c+2)H+D, H+2D, H+3D, 2H+2D, 3H+2D, 2H+3D\}$,
		
		\item [(97)] $\{H, H+D, 2H+D, H+2D,2H+2D, 3H+D, 2H+3D,3H+2D\}$,
		\item [(98)] $\{D, H+D, 2D,H+2D, 2H+D, H+3D,2H+2D,2H+3D\}$,
		\item [(99)] $\{H, D,H+D, 2H, H+2D,2H+D, 2H+2D,3H+2D\}$,
		\item [(100)] $\{-H+D, D, H, 2D, H+D, H+2D, 2H+2D,2H+3D\}$,
		\item [(101)] $\{H, 2H-D, H+D, 2H, 2H+D, 3H+D,3H+2D, 4H+2D\}$,
		\item [(102)] $\{H-D, D, H, H+D, 2H+D,2H+2D, 3H+2D, 2H+3D\}$,
		\item [(103)] $\{-H+2D, D, 2D, H+2D, H+3D, 2H+3D,H+4D,2H+4D\}$,
		\item [(104)] $\{H-D, H, 2H,2H+D, 3H+D,2H+2D,3H+2D,4H+D\}$,
		\item [(105)] $\{D, H+D,H+2D, 2H+2D, H+3D, 2H+3D, 3H+2D, 2H+4D\}$,
		
		\item [(106)] $\{H, H+D, 2H, H+2D,2H+D, 3H+D, 2H+2D,3H+2D\}$,
		\item [(107)] $\{D, H, 2D,H+D, 2H+D, H+2D,2H+2D,2H+3D\}$,
		\item [(108)] $\{H-D, D,H, 2H, H+D,2H+D, 2H+2D,3H+2D\}$,
		\item [(109)] $\{-H+2D, D, H+D, 2D, H+2D, H+3D, 2H+3D,2H+4D\}$,
		\item [(110)] $\{H-D, 2H-D, H, 2H, 2H+D, 3H+D,3H+2D, 4H+D\}$,
		\item [(111)] $\{H, D, H+D, H+2D, 2H+2D,2H+3D, 3H+2D, 2H+4D\}$,
		\item [(112)] $\{-H+D, D, 2D, H+2D, H+3D, 2H+2D,H+4D,2H+3D\}$,
		\item [(113)] $\{H, H+D, 2H+D,2H+2D, 3H+D,2H+3D,3H+2D,4H+2D\}$,
		\item [(114)] $\{D, H+D,H+2D, 2H+D, H+3D, 2H+2D, 3H+2D, 2H+3D\}$,
		
		\item [(115)] $\{H+D, H+2D, 2H+D, H+3D,2H+2D, 3H+D, 2H+3D,3H+2D\}$,
		\item [(116)] $\{D, H, 2D,H+D, 2H, H+2D,2H+D,2H+2D\}$,
		\item [(117)] $\{H-D, D,H, 2H-D, H+D,2H, 2H+D,3H+2D\}$,
		\item [(118)] $\{-H+2D, D, H, 2D, H+D, H+2D, 2H+3D,2H+4D\}$,
		\item [(119)] $\{H-D, 2H-2D, H, 2H-D, 2H, 3H+D,3H+2D, 4H+D\}$,
		\item [(120)] $\{H-D, D, H, H+D, 2H+2D,2H+3D, 3H+2D, 2H+4D\}$,
		\item [(121)] $\{-H+2D, D, 2D, H+3D, H+4D, 2H+3D,H+5D,2H+4D\}$,
		\item [(122)] $\{H-D, H, 2H+D,2H+2D, 3H+D,2H+3D,3H+2D,4H+D\}$,
		\item [(123)] $\{D, H+2D,H+3D, 2H+2D, H+4D, 2H+3D, 3H+2D, 2H+4D\}$,
		
		\item [(124)] $\{H+D, H+2D, 2H+D, H+3D,2H+2D, 3H+D, 3H+2D,2H+3D\}$,
		\item [(125)] $\{D, H, 2D,H+D, 2H, 2H+D,H+2D,2H+2D\}$,
		\item [(126)] $\{H-D, D,H, 2H-D, 2H,H+D, 2H+D,3H+2D\}$,
		\item [(127)] $\{-H+2D, D, H, H+D, 2D, H+2D, 2H+3D,2H+4D\}$,
		\item [(128)] $\{H-D, 2H-2D, 2H-D, H, 2H, 3H+D,3H+2D, 4H+D\}$,
		\item [(129)] $\{H-D, H, D, H+D, 2H+2D,2H+3D, 3H+2D, 2H+4D\}$,
		\item [(130)] $\{D, -H+2D, 2D, H+3D, H+4D, 2H+3D,H+5D,2H+4D\}$,
		\item [(131)] $\{-H+D, H, H+2D,H+3D, 2H+2D,H+4D,2H+3D,3H+2D\}$,
		\item [(132)] $\{H, 2H+D,2H+2D, 3H+D, 2H+3D, 3H+2D, 4H+D, 4H+2D\}$,
		
		\item [(133)] $\{H, H+D, 2H+D, H+2D,3H+D, 2H+2D, 2H+3D,3H+2D\}$,
		\item [(134)] $\{D, H+D, 2D,2H+D, H+2D, H+3D,2H+2D,2H+3D\}$,
		\item [(135)] $\{H, D,2H, H+D, H+2D,2H+D, 2H+2D,3H+2D\}$,
		\item [(136)] $\{-H+D, H, D, 2D, H+D, H+2D, 2H+2D,2H+3D\}$,
		\item [(137)] $\{2H-D, H, H+D, 2H, 2H+D, 3H+D,3H+2D, 4H+2D\}$,
		\item [(138)] $\{-H+D, -H+2D, D, 2D, H+2D,H+3D, 2H+3D, H+4D\}$,
		\item [(139)] $\{D, H, H+D, 2H+D, 2H+2D, 3H+2D,2H+3D,4H+2D\}$,
		\item [(140)] $\{H-D, H, 2H,2H+D, 3H+D,2H+2D,4H+D,3H+2D\}$,
		\item [(141)] $\{D, H+D,H+2D, 2H+2D, H+3D, 3H+2D, 2H+3D, 2H+4D\}$,
		
		\item [(142)] $\{H+D, H+2D, 2H+D, H+3D,3H+D, 2H+2D, 2H+3D,3H+2D\}$,
		\item [(143)] $\{D, H, 2D,2H, H+D, H+2D,2H+D,2H+2D\}$,
		\item [(144)] $\{H-D, D,2H-D, H, H+D,2H, 2H+D,3H+2D\}$,
		\item [(145)] $\{-H+2D, H, D, 2D, H+D, H+2D, 2H+3D,2H+4D\}$,
		\item [(146)] $\{2H-2D, H-D, H, 2H-D, 2H, 3H+D,3H+2D, 4H+D\}$,
		\item [(147)] $\{-H+D, -H+2D, D, 2D, H+3D,H+4D, 2H+3D, H+5D\}$,
		\item [(148)] $\{D, H, H+D, 2H+2D, 2H+3D, 3H+2D,2H+4D,4H+2D\}$,
		\item [(149)] $\{H-D, H, 2H+D,2H+2D, 3H+D,2H+3D,4H+D,3H+2D\}$,
		\item [(150)] $\{D, H+2D,H+3D, 2H+2D, H+4D, 3H+2D, 2H+3D, 2H+4D\}$,
		
		\item [(151)] $\{H+D, H+2D, 2H+D, H+3D,3H+D, 2H+2D, 3H+2D,2H+3D\}$,
		\item [(152)] $\{D, H, 2D,2H, H+D, 2H+D,2H+D,2H+2D\}$,
		\item [(153)] $\{H-D, D,2H-D, H, 2H,H+D, 2H+D,3H+2D\}$,
		\item [(154)] $\{-H+2D, H, D, H+D, 2D, H+2D, 2H+3D,2H+4D\}$,
		\item [(155)] $\{2H-2D, H-D, 2H-D, H, 2H, 3H+D,3H+2D, 4H+D\}$,
		\item [(156)] $\{-H+D, D, -H+2D, 2D, H+3D,H+4D, 2H+3D, H+5D\}$,
		\item [(157)] $\{H, D, H+D, 2H+2D, 2H+3D, 3H+2D,2H+4D,4H+2D\}$,
		\item [(158)] $\{-H+D, D, H+2D,H+3D, 2H+2D,H+4D,3H+2D,2H+3D\}$,
		\item [(159)] $\{H, 2H+D,2H+2D, 3H+D, 2H+3D, 4H+D, 3H+2D, 4H+2D\}$,
		
		\item [(160)] $\{H, H+D, 2H, 2H+D,H+2D, 3H+D, 2H+2D,3H+2D\}$,
		\item [(161)] $\{D, H, H+D,2D, 2H+D, H+2D,2H+2D,2H+3D\}$,
		\item [(162)] $\{H-D, H,D, 2H, H+D,2H+D, 2H+2D,3H+2D\}$,
		\item [(163)] $\{D, -H+2D, H+D, 2D, H+2D, H+3D, 2H+3D,2H+4D\}$,
		\item [(164)] $\{-H+D, H, D, H+D, H+2D, 2H+2D,2H+3D, 3H+2D\}$,
		\item [(165)] $\{2H-D, H, 2H, 2H+D, 3H+D,3H+2D, 4H+D, 4H+2D\}$,
		\item [(166)] $\{-H+D, D, 2D, H+2D, H+3D, 2H+2D,2H+3D,4H+2D\}$,
		\item [(167)] $\{H, H+D, 2H+D,2H+2D, 3H+D,3H+2D,2H+3D,4H+2D\}$,
		\item [(168)] $\{D, H+D,H+2D, 2H+D, 2H+2D, H+3D, 3H+2D, 2H+3D\}$,
		
		\item [(169)] $\{H+D, H+2D, 2H+D, 2H+2D,3H+D, H+3D, 2H+3D,3H+2D\}$,
		\item [(170)] $\{D, H, H+D,2H, 2D, H+2D,2H+D,2H+2D\}$,
		\item [(171)] $\{H-D, H,2H-D, D, H+D,2H, 2H+D,3H+2D\}$,
		\item [(172)] $\{D, H, -H+2D, 2D, H+D, H+2D, 2H+3D,2H+4D\}$,
		\item [(173)] $\{H-D, -H+D, D, H, H+D, 2H+2D,2H+3D, 3H+2D\}$,
		\item [(174)] $\{-2H+2D, -H+2D, D, 2D, H+3D,H+4D, 2H+3D, 2H+4D\}$,
		\item [(175)] $\{H, 2H-D, 2H, 3H+D, 3H+2D, 4H+D,4H+2D,5H+D\}$,
		\item [(176)] $\{H-D, H, 2H+D,2H+2D, 3H+D,3H+2D,4H+D,2H+3D\}$,
		\item [(177)] $\{D, H+2D,H+3D, 2H+2D, 2H+3D, 3H+2D, H+4D, 2H+4D\}$,
		
		\item [(178)] $\{H+D, H+2D, 2H+D, 2H+2D,3H+D, H+3D, 3H+2D,2H+3D\}$,
		\item [(179)] $\{D, H, H+D,2H, 2D, 2H+D,H+2D,2H+2D\}$,
		\item [(180)] $\{H-D, H,2H-D, D, 2H,H+D, 2H+D,3H+2D\}$,
		\item [(181)] $\{D, H, -H+2D, H+D, 2D, H+2D, 2H+3D,2H+4D\}$,
		\item [(182)] $\{H-D, -H+D, H, D, H+D, 2H+2D,2H+3D, 3H+2D\}$,
		\item [(183)] $\{-2H+2D, D, -H+2D, 2D, H+3D,H+4D, 2H+3D, 2H+4D\}$,
		\item [(184)] $\{2H-D,H, 2H, 3H+D, 3H+2D, 4H+D,4H+2D,5H+D\}$,
		\item [(185)] $\{-H+D, D, H+2D,H+3D, 2H+2D,2H+3D,3H+2D,H+4D\}$,
		\item [(186)] $\{H, 2H+D,2H+2D, 3H+D, 3H+2D, 4H+D, 2H+3D, 4H+2D\}$,
		
		\item [(187)] $\{H-D, H, 2H+D, 2H+2D,3H+D, 4H+D, 2H+3D,3H+2D\}$,
		\item [(188)] $\{D, H+2D, H+3D,2H+2D, 3H+2D, H+4D,2H+3D,2H+4D\}$,
		\item [(189)] $\{H+D, H+2D,2H+D, 3H+D, H+3D,2H+2D, 2H+3D,3H+2D\}$,
		\item [(190)] $\{D, H, 2H, 2D, H+D, H+2D, 2H+D,2H+2D\}$,
		\item [(191)] $\{H-D, 2H-D, D, H, H+D, 2H,2H+D, 3H+2D\}$,
		\item [(192)] $\{H, -H+2D, D, 2D, H+D,H+2D, 2H+3D, 2H+4D\}$,
		\item [(193)] $\{-2H+2D,-H+D, -H+2D, D, 2D, H+3D,H+4D,2H+3D\}$,
		\item [(194)] $\{H-D, H, 2H-D,2H, 3H+D,3H+2D,4H+D,5H+D\}$,
		\item [(195)] $\{D, H,H+D, 2H+2D, 2H+3D, 3H+2D, 4H+2D, 2H+4D\}$,
		
		\item [(196)] $\{H-D, H, 2H+D, 2H+2D,3H+D, 4H+D, 3H+2D,2H+3D\}$,
		\item [(197)] $\{D, H+2D, H+3D,2H+2D, 3H+2D, 2H+3D,H+4D,2H+4D\}$,
		\item [(198)] $\{H+D, H+2D,2H+D, 3H+D, 2H+2D,H+3D, 2H+3D,3H+2D\}$,
		\item [(199)] $\{D, H, 2H, H+D, 2D, H+2D, 2H+D,2H+2D\}$,
		\item [(200)] $\{H-D, 2H-D, H, D, H+D, 2H,2H+D, 3H+2D\}$,
		\item [(201)] $\{H, D, -H+2D, 2D, H+D,H+2D, 2H+3D, 2H+4D\}$,
		\item [(202)] $\{-H+D,-2H+2D, -H+2D, D, 2D, H+3D,H+4D,2H+3D\}$,
		\item [(203)] $\{-H+D, D, H,H+D, 2H+2D,2H+3D,3H+2D,4H+2D\}$,
		\item [(204)] $\{H, 2H-D,2H, 3H+D, 3H+2D, 4H+D, 5H+D, 4H+2D\}$,
		
		\item [(205)] $\{H-D, H, 2H-D, 2H,3H+D, 4H+D, 3H+2D,5H+D\}$,
		\item [(206)] $\{D, H, H+D,2H+2D, 3H+2D, 2H+3D,4H+2D,2H+4D\}$,
		\item [(207)] $\{H-D, H,2H+D, 3H+D, 2H+2D,4H+D, 2H+3D,3H+2D\}$,
		\item [(208)] $\{D, H+2D, 2H+2D, H+3D, 3H+2D, H+4D, 2H+3D,2H+4D\}$,
		\item [(209)] $\{H+D, 2H+D, H+2D, 3H+D, H+3D, 2H+2D,2H+3D, 3H+2D\}$,
		\item [(210)] $\{H, D, 2H, 2D, H+D,H+2D, 2H+D, 2H+2D\}$,
		\item [(211)] $\{-H+D,H, -H+2D, D, 2D, H+D,H+2D,2H+3D\}$,
		\item [(212)] $\{2H-D, D, H,H+D, 2H,2H+D,3H+2D,4H+2D\}$,
		\item [(213)] $\{-2H+2D, -H+D,-H+2D, D, 2D, H+3D, 2H+3D, H+4D\}$,
		
		\item [(214)] $\{H-D, H, 2H+D, 2H+2D,3H+D, 3H+2D, 2H+3D,4H+D\}$,
		\item [(215)] $\{D, H+2D, H+3D,2H+2D, 2H+3D, H+4D,3H+2D,2H+4D\}$,
		\item [(216)] $\{H+D, H+2D,2H+D, 2H+2D, H+3D,3H+D, 2H+3D,3H+2D\}$,
		\item [(217)] $\{D, H, H+D, 2D, 2H, H+2D, 2H+D,2H+2D\}$,
		\item [(218)] $\{H-D, H, D, 2H-D, H+D, 2H,2H+D, 3H+2D\}$,
		\item [(219)] $\{D, -H+2D, H, 2D, H+D,H+2D, 2H+3D, 2H+4D\}$,
		\item [(220)] $\{-H+D,H-D, D, H, H+D, 2H+2D,2H+3D,3H+2D\}$,
		\item [(221)] $\{2H-2D, H, 2H-D,2H, 3H+D,3H+2D,4H+D,4H+2D\}$,
		\item [(222)] $\{-H+2D, D,2D, H+3D, H+4D, 2H+3D, 2H+4D, H+5D\}$,
		
		\item [(223)] $\{H+D, H+2D, 2H+D, 3H+D,H+3D, 2H+2D, 3H+2D,2H+3D\}$,
		\item [(224)] $\{D, H, 2H,2D, H+D, 2H+D,H+2D,2H+2D\}$,
		\item [(225)] $\{H-D, 2H-D,D, H, 2H,H+D, 2H+D,3H+2D\}$,
		\item [(226)] $\{H, -H+2D, D, H+D, 2D, H+2D, 2H+3D,2H+4D\}$,
		\item [(227)] $\{-2H+2D, -H+D, D, -H+2D, 2D, H+3D,H+4D, 2H+3D\}$,
		\item [(228)] $\{H-D, 2H-D, H, 2D, 3H+D,3H+2D, 4H+D, 5H+D\}$,
		\item [(229)] $\{H,D, H+D, 2H+2D, 2H+3D, 3H+2D,4H+2D,2H+4D\}$,
		\item [(230)] $\{-H+D, D, H+2D,H+3D, 2H+2D,3H+2D,H+4D,2H+3D\}$,
		\item [(231)] $\{H, 2H+D,2H+2D, 3H+D, 4H+D, 2H+3D, 3H+2D, 4H+2D\}$,
		
		\item [(232)] $\{H+D, H+2D, 2H+D, 3H+D,2H+2D, H+3D, 3H+2D,2H+3D\}$,
		\item [(233)] $\{D, H, 2H,H+D, 2D, 2H+D,H+2D,2H+2D\}$,
		\item [(234)] $\{H-D, 2H-D,H, D, 2H,H+D, 2H+D,3H+2D\}$,
		\item [(235)] $\{H, D, -H+2D, H+D, 2D, H+2D, 2H+3D,2H+4D\}$,
		\item [(236)] $\{-H+D, -2H+2D, D, -H+2D, 2D, H+3D,H+4D, 2H+3D\}$,
		\item [(237)] $\{-H+D, H, D, H+D, 2H+2D,2H+3D, 3H+2D, 4H+2D\}$,
		\item [(238)] $\{2H-D,H, 2H, 3H+D, 3H+2D, 4H+D,5H+D,4H+2D\}$,
		\item [(239)] $\{-H+D, D, H+2D,H+3D, 2H+2D,3H+2D,2H+3D,H+4D\}$,
		\item [(240)] $\{H, 2H+D,2H+2D, 3H+D, 4H+D, 3H+2D, 2H+3D, 4H+2D\}$,
		
		\item [(241)] $\{H-D, H, 2H+D, 3H+D,2H+2D, 2H+3D, 4H+D,3H+2D\}$,
		\item [(242)] $\{D, H+2D, 2H+2D,H+3D, H+4D, 3H+2D,2H+3D,2H+4D\}$,
		\item [(243)] $\{H+D, 2H+D,H+2D, H+3D, 3H+D,2H+2D, 2H+3D,3H+2D\}$,
		\item [(244)] $\{H, D, 2D, 2H, H+D, H+2D, 2H+D,2H+2D\}$,
		\item [(245)] $\{-H+D, -H+2D, H, D, 2D, H+D,H+2D, 2H+3D\}$,
		\item [(246)] $\{D, 2H-D, H, H+D, 2H,2H+D, 3H+2D, 4H+2D\}$,
		\item [(247)] $\{2H-2D,H-D, H, 2H-D, 2H, 3H+D,4H+D,3H+2D\}$,
		\item [(248)] $\{-H+D, -H+2D, D,2D, H+3D,2H+3D,H+4D,H+5D\}$,
		\item [(249)] $\{D, H,H+D, 2H+2D, 3H+2D, 2H+3D, 2H+4D, 4H+2D\}$,
		
		\item [(250)] $\{H-D, H, 2H+D, 3H+D,2H+2D, 2H+3D, 3H+2D,4H+D\}$,
		\item [(251)] $\{D, H+2D, 2H+2D,H+3D, H+4D, 2H+3D,3H+2D,2H+4D\}$,
		\item [(252)] $\{H+D, 2H+D,H+2D, H+3D, 2H+2D,3H+D, 2H+3D,3H+2D\}$,
		\item [(253)] $\{H, D, 2D, H+D, 2H, H+2D, 2H+D,2H+2D\}$,
		\item [(254)] $\{-H+D, -H+2D, D, H, 2D, H+D,H+2D, 2H+3D\}$,
		\item [(255)] $\{D, H, 2H-D, H+D, 2H,2H+D, 3H+2D, 4H+2D\}$,
		\item [(256)] $\{H-D,2H-2D, H, 2H-D, 2H, 3H+D,4H+D,3H+2D\}$,
		\item [(257)] $\{H-D, D, H,H+D, 2H+2D,3H+2D,2H+3D,2H+4D\}$,
		\item [(258)] $\{-H+2D, D,2D, H+3D, 2H+3D, H+4D, H+5D, 2H+4D\}$,
		
		\item [(259)] $\{H-D, H, 2H+D, 3H+D,2H+2D, 4H+D, 3H+2D,2H+3D\}$,
		\item [(260)] $\{D, H+2D, 2H+2D,H+3D, 3H+2D, 2H+3D,H+4D,2H+4D\}$,
		\item [(261)] $\{H+D, 2H+D,H+2D, 3H+D, 2H+2D,H+3D, 2H+3D,3H+2D\}$,
		\item [(262)] $\{H, D, 2H, H+D, 2D, H+2D, 2H+D,2H+2D\}$,
		\item [(263)] $\{-H+D, H, D, -H+2D, 2D, H+D,H+2D, 2H+3D\}$,
		\item [(264)] $\{2H-D, H, D, H+D, 2H,2H+D, 3H+2D, 4H+2D\}$,
		\item [(265)] $\{-H+D,-2H+2D, -H+2D, D, 2D, H+3D,2H+3D,H+4D\}$,
		\item [(266)] $\{-H+D, D, H,H+D, 2H+2D,3H+2D,2H+3D,4H+2D\}$,
		\item [(267)] $\{H, 2H-D,2H, 3H+D, 4H+D, 3H+2D, 5H+D, 4H+2D\}$,
		
		\item [(268)] $\{H-D, H, 2H+D, 3H+D,2H+2D, 3H+2D, 2H+3D,4H+D\}$,
		\item [(269)] $\{D, H+2D, 2H+2D,H+3D, 2H+3D, H+4D,3H+2D,2H+4D\}$,
		\item [(270)] $\{H+D, 2H+D,H+2D, 2H+2D, H+3D,3H+D, 2H+3D,3H+2D\}$,
		\item [(271)] $\{H, D, H+D, 2D, 2H, H+2D, 2H+D,2H+2D\}$,
		\item [(272)] $\{-H+D, D, -H+2D, H, 2D, H+D,H+2D, 2H+3D\}$,
		\item [(273)] $\{H, D, 2H-D, H+D, 2H,2H+D, 3H+2D, 4H+2D\}$,
		\item [(274)] $\{-H+D,H-D, D, H, H+D, 2H+2D,3H+2D,2H+3D\}$,
		\item [(275)] $\{2H-2D, H, 2H-D,2H, 3H+D,4H+D,3H+2D,4H+2D\}$,
		\item [(276)] $\{-H+2D, D,2D, H+3D, 2H+3D, H+4D, 2H+4D, H+5D\}$,
		
		\item [(277)] $\{H-D, H, 2H+D, 3H+D,2H+2D, 3H+2D, 4H+D,2H+3D\}$,
		\item [(278)] $\{D, H+2D, 2H+2D,H+3D, 2H+3D, 3H+2D,H+4D,2H+4D\}$,
		\item [(279)] $\{H+D, 2H+D,H+2D, 2H+2D, 3H+D,H+3D, 2H+3D,3H+2D\}$,
		\item [(280)] $\{H, D, H+D, 2H, 2D, H+2D, 2H+D,2H+2D\}$,
		\item [(281)] $\{-H+D, D, H, -H+2D, 2D, H+D,H+2D, 2H+3D\}$,
		\item [(282)] $\{H, 2H-D,D, H+D, 2H,2H+D, 3H+2D, 4H+2D\}$,
		\item [(283)] $\{H-D,-H+D, D, H, H+D, 2H+2D,3H+2D,2H+3D\}$,
		\item [(284)] $\{-2H+2D, -H+2D, D,2D, H+3D,2H+3D,H+4D,2H+4D\}$,
		\item [(285)] $\{H, 2H-D,2H, 3H+D, 4H+D, 3H+2D, 4H+2D, 5H+D\}$,
		
		\item [(286)] $\{H+D, H+2D, 2H+D, 2H+2D,H+3D, 3H+D, 3H+2D,2H+3D\}$,
		\item [(287)] $\{D, H, H+D,2D, 2H, 2H+D,H+2D,2H+2D\}$,
		\item [(288)] $\{H-D, H,D, 2H-D, 2H,H+D, 2H+D,3H+2D\}$,
		\item [(289)] $\{D, -H+2D, H, H+D, 2D, H+2D, 2H+D,3H+2D\}$,
		\item [(290)] $\{-H+D, H-D, H, D, H+D, 2H+2D,2H+3D, 3H+2D\}$,
		\item [(291)] $\{2H-2D, 2H-D,H, 2H, 3H+D,3H+2D, 4H+D, 4H+2D\}$,
		\item [(292)] $\{D,-H+2D, 2D, H+3D, H+4D, 2H+3D,3H+4D,H+5D\}$,
		\item [(293)] $\{-H+D, D, H+2D,H+3D, 2H+2D,2H+3D,H+4D,3H+2D\}$,
		\item [(294)] $\{H, 2H+D,2H+2D, 3H+D, 3H+2D, 2H+3D, 4H+D, 4H+2D\}$,
		
		\item [(295)] $\{H-D, H, D, H+D,2H+2D, 3H+2D, 2H+3D,2H+4D\}$,
		\item [(296)] $\{D, -H+2D, 2D,H+3D, 2H+3D, H+4D,H+5D,2H+4D\}$,
		\item [(297)] $\{-H+D, D,H+2D, 2H+2D, H+3D,H+4D, 2H+3D,3H+2D\}$,
		\item [(298)] $\{H, 2H+D, 3H+D, 2H+2D, 2H+3D, 3H+2D, 4H+D,4H+2D\}$,
		\item [(299)] $\{H+D, 2H+D, H+2D, H+3D, 2H+2D, 3H+D,3H+2D, 2H+3D\}$,
		\item [(300)] $\{H, D,2D, H+D, 2H,2H+D, H+2D, 2H+2D\}$,
		\item [(301)] $\{-H+D,-H+2D, D, H, H+D, 2D,H+2D,2H+3D\}$,
		\item [(302)] $\{D, H, 2H-D,2H, H+D,2H+D,3H+2D,4H+2D\}$,
		\item [(303)] $\{H-D, 2H-2D,2H-D, H, 2H, 3H+D, 4H+D, 3H+2D\}$,
		
		\item [(304)] $\{H+D, 2H+D, H+2D, H+3D,3H+D, 2H+2D, 3H+2D,2H+3D\}$,
		\item [(305)] $\{H, D, 2D,2H, H+D, 2H+D,H+2D,2H+2D\}$,
		\item [(306)] $\{-H+D, -H+2D,H, D, H+D,2D, H+2D,2H+3D\}$,
		\item [(307)] $\{D, 2H-D, H, 2H, H+D, 2H+D, 3H+2D,4H+2D\}$,
		\item [(308)] $\{2H-2D, H-D, 2H-D, H, 2H, 3H+D,4H+D, 3H+2D\}$,
		\item [(309)] $\{-H+D, D,-H+2D, 2D, H+3D,2H+3D, H+4D, H+5D\}$,
		\item [(310)] $\{H,D, H+D, 2H+2D, 3H+2D, 2H+3D,2H+4D,4H+2D\}$,
		\item [(311)] $\{-H+D, D, H+2D,2H+2D, H+3D,H+4D,3H+2D,2H+3D\}$,
		\item [(312)] $\{H, 2H+D,3H+D, 2H+2D, 2H+3D, 4H+D, 3H+2D, 4H+2D\}$,
		
		\item [(313)] $\{H+D, 2H+D, H+2D, 2H+2D,3H+D, H+3D, 3H+2D,2H+3D\}$,
		\item [(314)] $\{H, D, H+D,2H, 2D, 2H+D,H+2D,2H+2D\}$,
		\item [(315)] $\{-H+D, D,H, -H+2D, H+D,2D, H+2D,2H+3D\}$,
		\item [(316)] $\{H, 2H-D, D, 2H, H+D, 2H+D, 3H+2D,4H+2D\}$,
		\item [(317)] $\{H-D, -H+D, H, D, H+D, 2H+2D,3H+2D, 2H+3D\}$,
		\item [(318)] $\{-2H+2D, D,-H+2D, 2D, H+3D,2H+3D, H+4D, 2H+4D\}$,
		\item [(319)] $\{2H-D,H, 2H, 3H+D, 4H+D, 3H+2D,4H+2D,5H+D\}$,
		\item [(320)] $\{-H+D, D, H+2D,2H+2D, H+3D,2H+3D,3H+2D,H+4D\}$,
		\item [(321)] $\{H, 2H+D,3H+D, 2H+2D, 3H+2D, 4H+D, 2H+3D, 4H+2D\}$,
		
		\item [(322)] $\{H-D, 2H-D, H, 2H,3H+D, 4H+D, 3H+2D,5H+D\}$,
		\item [(323)] $\{H, D, H+D,2H+2D, 3H+2D, 2H+3D,4H+2D,2H+4D\}$,
		\item [(324)] $\{-H+D, D,H+2D, 2H+2D, H+3D,3H+2D, H+4D,2H+3D\}$,
		\item [(325)] $\{H, 2H+D, 3H+D, 2H+2D, 4H+D, 2H+3D, 3H+2D,4H+2D\}$,
		\item [(326)] $\{H+D, 2H+D, H+2D, 3H+D, H+3D, 2H+2D,3H+2D, 2H+3D\}$,
		\item [(327)] $\{H, D,2H, 2D, H+D,2H+D, H+2D, 2H+2D\}$,
		\item [(328)] $\{-H+D,H, -H+2D, D, H+D, 2D,H+2D,2H+3D\}$,
		\item [(329)] $\{2H-D, D, H,2H, H+D,2H+D,3H+2D,4H+2D\}$,
		\item [(330)] $\{-2H+2D, -H+D,D, -H+2D, 2D, H+3D, 2H+3D, H+4D\}$,
		
		\item [(331)] $\{H+D, 2H+D, H+2D, 3H+D,2H+2D, H+3D, 3H+2D,2H+3D\}$,
		\item [(332)] $\{H, D, 2H,H+D, 2D, 2H+D,H+2D,2H+2D\}$,
		\item [(333)] $\{-H+D, H,D, -H+2D, H+D,2D, H+2D,2H+3D\}$,
		\item [(334)] $\{2H-D, H, D, 2H, H+D, 2H+D, 3H+2D,4H+2D\}$,
		\item [(335)] $\{-H+D, -2H+2D, D, -H+2D, 2D, H+3D,2H+3D, 4H+2D\}$,
		\item [(336)] $\{-H+D, H,D, H+D, 2H+2D,3H+2D, 2H+3D, 4H+2D\}$,
		\item [(337)] $\{2H-D,H, 2H, 3H+D, 4H+D, 3H+2D,2H+3D,H+4D\}$,
		\item [(338)] $\{-H+D, D, H+2D,2H+2D, H+3D,3H+2D,2H+3D,H+4D\}$,
		\item [(339)] $\{H, 2H+D,3H+D, 2H+2D, 4H+D, 3H+2D, 2H+3D, 4H+2D\}$,
		
		\item [(340)] $\{H, 2H+D, 3H+D, 2H+2D,3H+2D, 2H+3D, 4H+D,4H+2D\}$,
		\item [(341)] $\{H+D, 2H+D, H+2D,2H+2D, H+3D, 3H+D,3H+2D,2H+3D\}$,
		\item [(342)] $\{H, D,H+D, 2D, 2H,2H+D, H+2D,2H+2D\}$,
		\item [(343)] $\{-H+D, D, -H+2D, H, H+D, 2D, H+2D,2H+3D\}$,
		\item [(344)] $\{H, D, 2H-D, 2H, H+D, 2H+D,3H+2D, 4H+2D\}$,
		\item [(345)] $\{-H+D, H-D,H, D, H+D,2H+2D, 3H+2D, 2H+3D\}$,
		\item [(346)] $\{2H-2D,2H-D, H, 2H, 3H+D, 4H+D,3H+2D,4H+2D\}$,
		\item [(347)] $\{D, -H+2D, 2D,H+3D, 2H+3D,H+4D,2H+4D,H+5D\}$,
		\item [(348)] $\{-H+D, D,H+2D, 2H+2D, H+3D, 2H+3D, H+4D, 3H+2D\}$,
	\end{itemize}
	where $a,b,c,d\in \mathbb{Z}$. Moreover, by mutations and normalizations, they are related as:
	
	\begin{align}
	&(1)\rightarrow (2)\rightarrow (3)\rightarrow (1), \nonumber\\
	&(4)\rightarrow (5)\rightarrow (6)\rightarrow (7)\rightarrow (8)\rightarrow (9)\rightarrow  (10)\rightarrow (11)\rightarrow (12)\rightarrow (4), \nonumber\\
	&(13)\rightarrow (14)\rightarrow (15)\rightarrow (16)\rightarrow (17)\rightarrow (18)\rightarrow  (19)\rightarrow (20)\rightarrow (21)\rightarrow (13), \nonumber\\
	&(22)\rightarrow (23)\rightarrow (24)\rightarrow (25)\rightarrow (26)\rightarrow (27)\rightarrow (28)\rightarrow (29)\rightarrow (30)\rightarrow (22),\nonumber\\
	&(31)\rightarrow (32)\rightarrow (33)\rightarrow (34)\rightarrow (35)\rightarrow (36)\rightarrow (37)\rightarrow (38)\rightarrow (39)\rightarrow (31), \nonumber\\
	&(40)\rightarrow (41)\rightarrow (42)\rightarrow (43)\rightarrow (44)\rightarrow (45)\rightarrow (46)\rightarrow (47)\rightarrow (48)\rightarrow (40),\nonumber\\
	&(49)\rightarrow (50)\rightarrow (51)\rightarrow (52)\rightarrow (53)\rightarrow (54)\rightarrow (55)\rightarrow (56)\rightarrow (57)\rightarrow (49),\nonumber\\
	&(58)\rightarrow (59)\rightarrow (60)\rightarrow (61)\rightarrow (62)\rightarrow (63)\rightarrow (64)\rightarrow (65)\rightarrow (66)\rightarrow (58),\nonumber\\
	&(67)\rightarrow (68)\rightarrow (69)\rightarrow (70)\rightarrow (71)\rightarrow (72)\rightarrow (73)\rightarrow (74)\rightarrow (75)\rightarrow (67),\nonumber\\
	&(76)\rightarrow (77)\rightarrow (78)\rightarrow (76),\nonumber\\
	&(79)\rightarrow (80)\rightarrow (81)\rightarrow (82)\rightarrow (83)\rightarrow (84)\rightarrow (85)\rightarrow (86)\rightarrow (87)\rightarrow (79),\nonumber\\
	&(88)\rightarrow (89)\rightarrow (90)\rightarrow (91)\rightarrow (92)\rightarrow (93)\rightarrow (94)\rightarrow (95)\rightarrow (96)\rightarrow (88),\nonumber\\
	&(97)\rightarrow (98)\rightarrow (90)\rightarrow (100)\rightarrow (101)\rightarrow (102)\rightarrow (103)\rightarrow (104)\rightarrow (105)\rightarrow (97),\nonumber\\
	&(106)\rightarrow (107)\rightarrow (108)\rightarrow (109)\rightarrow (110)\rightarrow (111)\rightarrow (112)\rightarrow (113)\rightarrow (114)\rightarrow (106),\nonumber\\
	&(115)\rightarrow (116)\rightarrow (117)\rightarrow (118)\rightarrow (119)\rightarrow (120)\rightarrow (121)\rightarrow (122)\rightarrow (123)\rightarrow (115),\nonumber\\
	&(124)\rightarrow (125)\rightarrow (126)\rightarrow (127)\rightarrow (128)\rightarrow (129)\rightarrow (130)\rightarrow (131)\rightarrow (132)\rightarrow (124),\nonumber\\
	&(133)\rightarrow (134)\rightarrow (135)\rightarrow (136)\rightarrow (137)\rightarrow (138)\rightarrow (139)\rightarrow (140)\rightarrow (141)\rightarrow (133),\nonumber
	\end{align}
	\begin{align}
	&(142)\rightarrow (143)\rightarrow (144)\rightarrow (145)\rightarrow (146)\rightarrow (147)\rightarrow (148)\rightarrow (149)\rightarrow (150)\rightarrow (142),\nonumber\\
	&(151)\rightarrow (152)\rightarrow (153)\rightarrow (154)\rightarrow (155)\rightarrow (156)\rightarrow (157)\rightarrow (158)\rightarrow (159)\rightarrow (151),\nonumber\\
	&(160)\rightarrow (161)\rightarrow (162)\rightarrow (163)\rightarrow (164)\rightarrow (165)\rightarrow (166)\rightarrow (167)\rightarrow (168)\rightarrow (160),\nonumber\\
	&(169)\rightarrow (170)\rightarrow (171)\rightarrow (172)\rightarrow (173)\rightarrow (174)\rightarrow (175)\rightarrow (176)\rightarrow (177)\rightarrow (169),\nonumber\\
	&(178)\rightarrow (179)\rightarrow (180)\rightarrow (181)\rightarrow (182)\rightarrow (183)\rightarrow (184)\rightarrow (185)\rightarrow (186)\rightarrow (178),\nonumber\\
	&(187)\rightarrow (188)\rightarrow (189)\rightarrow (190)\rightarrow (191)\rightarrow (192)\rightarrow (193)\rightarrow (194)\rightarrow (195)\rightarrow (187),\nonumber\\
	&(196)\rightarrow (197)\rightarrow (198)\rightarrow (199)\rightarrow (200)\rightarrow (201)\rightarrow (202)\rightarrow (203)\rightarrow (204)\rightarrow (196),\nonumber\\
	&(205)\rightarrow (206)\rightarrow (207)\rightarrow (208)\rightarrow (200)\rightarrow (210)\rightarrow (211)\rightarrow (212)\rightarrow (213)\rightarrow (205),\nonumber\\
	&(214)\rightarrow (215)\rightarrow (216)\rightarrow (217)\rightarrow (218)\rightarrow (219)\rightarrow (220)\rightarrow (221)\rightarrow (222)\rightarrow (214),\nonumber\\
	&(223)\rightarrow (224)\rightarrow (225)\rightarrow (226)\rightarrow (227)\rightarrow (228)\rightarrow (229)\rightarrow (230)\rightarrow (231)\rightarrow (223),\nonumber\\
	&(232)\rightarrow (233)\rightarrow (234)\rightarrow (235)\rightarrow (236)\rightarrow (237)\rightarrow (238)\rightarrow (239)\rightarrow (240)\rightarrow (232),\nonumber\\
	&(241)\rightarrow (242)\rightarrow (243)\rightarrow (244)\rightarrow (245)\rightarrow (246)\rightarrow (247)\rightarrow (248)\rightarrow (249)\rightarrow (241),\nonumber\\
	&(250)\rightarrow (251)\rightarrow (252)\rightarrow (253)\rightarrow (254)\rightarrow (255)\rightarrow (256)\rightarrow (257)\rightarrow (258)\rightarrow (250),\nonumber\\
	&(259)\rightarrow (260)\rightarrow (261)\rightarrow (262)\rightarrow (263)\rightarrow (264)\rightarrow (265)\rightarrow (266)\rightarrow (267)\rightarrow (259),\nonumber\\
	&(268)\rightarrow (269)\rightarrow (270)\rightarrow (271)\rightarrow (272)\rightarrow (273)\rightarrow (274)\rightarrow (275)\rightarrow (276)\rightarrow (268),\nonumber\\
	&(277)\rightarrow (278)\rightarrow (279)\rightarrow (280)\rightarrow (281)\rightarrow (282)\rightarrow (283)\rightarrow (284)\rightarrow (285)\rightarrow (277),\nonumber\\
	&(286)\rightarrow (287)\rightarrow (288)\rightarrow (289)\rightarrow (290)\rightarrow (291)\rightarrow (292)\rightarrow (293)\rightarrow (294)\rightarrow (286),\nonumber\\
	&(295)\rightarrow (296)\rightarrow (297)\rightarrow (298)\rightarrow (299)\rightarrow (300)\rightarrow (301)\rightarrow (302)\rightarrow (303)\rightarrow (295),\nonumber\\
	&(304)\rightarrow (305)\rightarrow (306)\rightarrow (307)\rightarrow (308)\rightarrow (309)\rightarrow (310)\rightarrow (311)\rightarrow (312)\rightarrow (304),\nonumber\\
	&(313)\rightarrow (314)\rightarrow (315)\rightarrow (316)\rightarrow (317)\rightarrow (318)\rightarrow (319)\rightarrow (320)\rightarrow (321)\rightarrow (313),\nonumber\\
	&(322)\rightarrow (323)\rightarrow (324)\rightarrow (325)\rightarrow (326)\rightarrow (327)\rightarrow (328)\rightarrow (329)\rightarrow (330)\rightarrow (322),\nonumber\\
	&(331)\rightarrow (332)\rightarrow (333)\rightarrow (334)\rightarrow (335)\rightarrow (336)\rightarrow (337)\rightarrow (338)\rightarrow (339)\rightarrow (331),\nonumber\\
	&(340)\rightarrow (341)\rightarrow (342)\rightarrow (343)\rightarrow (344)\rightarrow (345)\rightarrow (346)\rightarrow (347)\rightarrow (348)\rightarrow (340),\nonumber
	\end{align}
	\begin{proof}
		Note that $\OO_X(aH+bD)$ is cohomologically zero line bundle if $a=-1, 2$ or $b=-1,-2$. Based on Table \ref{tab:table11} and the rules below, we get the result.
		\begin{itemize}
			\item [(1)] Since $B_0, B_1$ and $B_2, B_3$ are symmetric, we only need to find exceptional collections such that $D_1$ is $B_0$ or $B_1$.
			\item [(2)] $B_0,B_1,B_2$ and $B_3$ can occur at most three times.
			\item [(3)] $B_0$ cannot occur three times after $B_1$.
			\item [(4)] $B_2$ cannot occur three times after $B_3$.
		\end{itemize} 
		\begin{table}[h!]
		\begin{center}
			\caption{Exceptional pairs.}
			\label{tab:table11}
			\begin{adjustbox}{max width=\textwidth}
				\begin{tabular}{c|c|c|c|c} 
					& $B_0'$ & $B_1'$ & $B_2'$ & $B_3'$\\
					\hline
					$B_0$ &	$a'=a+1, a+2$ & \checkmark  & $a'=-1,0$ or $b=2,3$  & $a'=0,1$ or $d=2,3$  \\
					\hline
					$B_1$ &	$a'=b+1, b+2$ & $b'=b+1, b+2$ & $b=-1,0$ or $c=3,4$  & $b=0,1$ or $d=3,4$  \\
					\hline
					$B_2$ &	$a'=2,3$ or $c=-1,0$ & $b'=2,3$ or $c=3,4$ & $c'=c+1, c+2$ & \checkmark\\
					\hline
					$B_3$ & $a'=3,4$ or $d=-1,0$ & $b'=3,4$ or $d=0,1$ & $c'=d+1, d+2$ & $d'=d+1,d+2$\\
					\hline
				\end{tabular}
			\end{adjustbox}
			\vspace{3mm}
			
			where $B_0=aH+D, B_1=bH+2D, B_2=H+cD, B_3=2H+dD$.
		\end{center}
	\end{table}
	\end{proof}
\end{thm}

\begin{thm}
	Let $X$ be $\PP^2 \times \PP^2$. Then any exceptional collection of line bundles of length 9 on $X$ is full.
	\begin{proof}
		By Lemma \ref{nor}, we only need to show that any exceptional collection of line bundles of length $9$ in Theorem \ref{cl11} is full. By Lemma \ref{mt}, it suffices to show that the exceptional collections of type $(2), (22), (13)$, $(27)$, $(23)$, $(12)$, $(35)$, $(43)$, $(59)$ and $(54)$ in Theorem \ref{cl11} are full. 
		
		We first show that exceptional collection of type $(2)$ is full. Consider the projective bundle $f:X\cong \PP_{\PP^2}(\OO(-a)\oplus \OO(-a)\oplus \OO(-a))\rightarrow \PP^2$, where $\OO_X(1)=H+aD$. Then by the projective bundle formula, we have the following semiorthogonal decomposition of $\D(X)$.
		\begin{align}
		\D(X)&=\langle f^*\D(\PP^2),f^*\D(\PP^2)\otimes \OO_X(1), f^*\D(\PP^2)\otimes \OO_X(2)\rangle \nonumber\\
		&=\langle f^*\OO_{\PP^2},f^*\OO_{\PP^2}(1),f^*\OO_{\PP^2}(2),\nonumber\\
		&\qquad f^*\OO_{\PP^2}\otimes \OO_X(1),f^*\OO_{\PP^2}(1)\otimes\OO_X(1),f^*\OO_{\PP^2}(2)\otimes \OO_X(1),\nonumber\\
		&\qquad f^*\OO_{\PP^2}\otimes \OO_X(2),f^*\OO_{\PP^2}(1)\otimes\OO_X(2),f^*\OO_{\PP^2}(2)\otimes \OO_X(2)\rangle\nonumber\\
		&=\langle \OO_X,\OO_X(D),\OO_X(2D),\nonumber\\
		&\qquad \OO_X(H+aD),\OO_X(H+(a+1)D),\OO_X(H+(a+2)D),\nonumber\\
		&\qquad \OO_X(2H+2aD),\OO_X(2H+(2a+1)D),\OO_X(2H+(2a+2)D)\rangle.\nonumber
		\end{align}
		Based on this semiorthogonal decomposition, we show that the exceptional collection of type $(2)$ is full by mathematical induction. For a given $a\in \mathbb{Z}$, assume that for $b=k$ the exceptional collection
		\begin{align}\label{mi7} \tag{E10}
		\{D, 2D, H+aD, H+(a+1)D, H+(a+2)D, 2H+kD, 2H+(k+1)D, 2H+(k+2)D\}
		\end{align}
		is full. Then for $d=k-1$, we have the exceptional collection
		\begin{align}\label{mi8} \tag{E11}
		\{D, 2D, H+aD, H+(a+1)D, H+(a+2)D, 2H+(k-1)D, 2H+kD, 2H+(k+1)D\}
		\end{align}
		and for $d=k+1$, we have the exceptional collection
		\begin{align}\label{mi9} \tag{E12}
		\{D, 2D, H+aD, H+(a+1)D, H+(a+2)D, 2H+(k+1)D, 2H+(k+2)D, 2H+(k+3)D\}.
		\end{align}
		Comparing the exceptional collections (\ref{mi8}) and (\ref{mi9}) with (\ref{mi7}) and by Lemma \ref{mt1}, the exceptional collections (\ref{mi8}) and (\ref{mi9}) are full. Hence the exceptional collection of type $(2)$ is also full.
		
		Now we show that exceptional collection of type $(76)$ is full. Consider the projective bundle $g:X\cong \PP_{\PP^2}(\OO(-b)\oplus \OO(-b)\oplus \OO(-b))\rightarrow \PP^2$, where $\OO_X(1)=bH+D$. Then by the projective bundle formula, we obtain the following semiorthogonal decomposition of $\D(X)$.
		\begin{align}
		\D(X)&=\langle g^*\D(\PP^2),g^*\D(\PP^2)\otimes \OO_X(1),g^*\D(\PP^2)\otimes \OO_X(2)\rangle \nonumber\\
		&=\langle g^*\OO_{\PP^2},g^*\OO_{\PP^2}(1),g^*\OO_{\PP^2}(2),\nonumber\\
		&\qquad g^*\OO_{\PP^2}\otimes \OO_X(1),g^*\OO_{\PP^2}(1)\otimes \OO_X(1),g^*\OO_{\PP^2}(2)\otimes \OO_X(1),\nonumber\\
		&\qquad g^*\OO_{\PP^2}\otimes \OO_X(2),g^*\OO_{\PP^2}(1)\otimes \OO_X(2),g^*\OO_{\PP^2}(2)\otimes \OO_X(2)\rangle\nonumber\\
		&=\langle \OO_X,\OO_X(H),\OO_X(2H),\OO_X(bH+D),\OO_X((b+1)H+D),\OO_X((b+2)H+D)\nonumber\\
		&\qquad \OO_X(2H+2bD),\OO_X(2H+(2b+1)D),\OO_X(2H+(2b+2)D)\rangle.\nonumber
		\end{align}
		Based on this semiorthogonal decomposition, we show that the exceptional collection of type $(76)$ is full by mathematical induction. For a given $b,c\in \mathbb{Z}$, assume that for $d=k$ the exceptional collection
		\begin{align}\label{mi10} \tag{E13}
		\{H, 2H, cH+D, (c+1)H+D, (c+2)H+D, kH+2D, (k+1)H+2D,(k+2)H+2D\}
		\end{align}
		is full. Then for $d=k-1$, we have the exceptional collection
		\begin{align}\label{mi11} \tag{E14}
		\{H, 2H, cH+D, (c+1)H+D, (c+2)H+D, (k-1)H+2D, kH+2D,(k+1)H+2D\}
		\end{align}
		and for $d=k+1$, we have the exceptional collection
		\begin{align}\label{mi12} \tag{E15}
		\{H, 2H, cH+D, (c+1)H+D, (c+2)H+D, (k+1)H+2D, (k+2)H+2D,(k+3)H+2D\}.
		\end{align}
		Comparing the exceptional collections (\ref{mi11}) and (\ref{mi12}) with (\ref{mi10}) and by Lemma \ref{mt1}, the exceptional collections (\ref{mi11}) and (\ref{mi12}) are full. Hence the exceptional collection of type $(76)$ is also full.
		
		The other exceptional collections are summarized in Table \ref{tab:table13}.
		
		\begin{table}[h!]
			\begin{center}
				\caption{Showing the fullness.}
				\label{tab:table13}
				\begin{tabular}{c|c} 
					Type & How to show \\
					\hline
					$(4)$ & Set $b=0$ in $(1)$ and mutate $2H+2D$ and $3H+D$  \\
					\hline
					$(13)$ & Set $b=1$ in $(1)$ and mutate $2H+3D$ and $3H+D$  \\
					\hline
					$(22)$ & Mutate $2H+2D$ and $3H+D$ in $(13)$  \\
					\hline
					$(31)$ & Mutate $2H+3D$ and $3H+2D$ in $(22)$  \\
					\hline
					$(46)$ & Mutate $3H+D$ and $2H+2D$ in $(31)$  \\
					\hline
					$(64)$ & Set $d=c+1$ in $(76)$ and mutate $(c+1)H+2D$ and $(c+2)H+2D$  \\
					\hline
					$(79)$ & Set $d=c$ in $(76)$ and mutate $(c+2)H+D$ and $(c+1)H+2D$  \\
					\hline
					$(70)$ & Mutate $(c+2)H+D$ and $(c+1)H+2D$ in $(79)$  \\
					\hline
					$(88)$ & Mutate $(c+1)H+D$ and $cH+2D$ in $(79)$  \\
					\hline
					$(54)$ & Mutate $(c+2)H+D$ and $(c+1)H+2D$ in $(88)$  \\
					\hline
					$(97)$ & Set $a=0$ in $(7)$ and mutate $2H+3D$ and $3H+D$  \\
					\hline
					$(106)$ & Set $a=0$ in $(16)$ and mutate $2H+2D$ and $3H+D$  \\
					\hline
					$(115)$ & Set $a=1$ in $(13)$ and mutate $H+3D$ and $2H+D$  \\
					\hline
					$(124)$ & Mutate $2H+3D$ and $3H+2D$ in $(115)$  \\
					\hline
					$(133)$ & Set $a=0$ in $(22)$ and mutate $H+2D$ and $2H+D$  \\
					\hline
					$(142)$ & Set $a=1$ in $(22)$ and mutate $H+3D$ and $2H+D$  \\
					\hline
					$(151)$ & Mutate $2H+3D$ and $3H+2D$ in $(142)$  \\
					\hline
					$(160)$ & Mutate $H+2D$ and $2H+D$ in $(106)$  \\
					\hline
					$(216)$ & Set $a=1$ in $(40)$ and mutate $2H+3D$ and $3H+D$  \\
					\hline
					$(169)$ & Mutate $H+3D$ and $3H+D$ in $(216)$  \\
					\hline
					$(178)$ & Mutate $2H+3D$ and $3H+2D$ in $(169)$  \\
					\hline
					$(198)$ & Mutate $2H+2D$ and $3H+D$ in $(169)$  \\
					\hline
					$(189)$ & Mutate $2H+2D$ and $H+3D$ in $(198)$  \\
					\hline
					$(209)$ & Mutate $H+2D$ and $2H+D$ in $(189)$  \\
					\hline
					$(223)$ & Mutate $2H+3D$ and $3H+2D$ in $(189)$  \\
					\hline
					$(232)$ & Mutate $H+3D$ and $2H+2D$ in $(223)$  \\
					\hline
					$(243)$ & Mutate $3H+D$ and $H+3D$ in $(209)$  \\
					\hline
					$(252)$ & Mutate $3H+D$ and $2H+2D$ in $(243)$  \\
					\hline
					$(261)$ & Mutate $H+3D$ and $2H+2D$ in $(209)$  \\
					\hline
					$(270)$ & Set $a=1$ in $(34)$ and mutate $2H+3D$ and $3H+D$  \\
					\hline
					$(279)$ & Mutate $H+3D$ and $3H+D$ in $(270)$  \\
					\hline
					$(286)$ & Mutate $2H+3D$ and $3H+2D$ in $(216)$  \\
					\hline
					$(299)$ & Mutate $2H+3D$ and $3H+2D$ in $(252)$  \\
					\hline
					$(304)$ & Mutate $2H+2D$ and $3H+D$ in $(299)$  \\
					\hline
					$(313)$ & Mutate $2H+3D$ and $3H+2D$ in $(279)$  \\
					\hline
				\end{tabular}
			\end{center}
		\end{table}
	
			\begin{table}[h!]
		\begin{center}
			\begin{tabular}{c|c} 
			$(326)$ & Mutate $2H+3D$ and $3H+2D$ in $(209)$  \\
			\hline
			$(331)$ & Mutate $H+3D$ and $2H+2D$ in $(326)$  \\
			\hline
			$(341)$ & Mutate $3H+D$ and $H+3D$ in $(313)$  \\
			\hline
			\end{tabular}
		\end{center}
	\end{table}
	\newpage
		Therefore, all the exceptional collections of length $9$ on $X$ are full.
	\end{proof}
\end{thm}

As a result, the conjecture \ref{Kuzc} holds true for smooth toric Fano fourfolds of Picard rank two and hence our main result Theorem \ref{main} holds.

\section{Final Remark}

 As we can see, the list of cases that we need to consider become longer when the Picard rank gets higher and even longer for the product of projective spaces. We finalize this paper with a question:
 
 \begin{que}
 	Is there a way to prove the conjecture \ref{Kuzc} without checking all the exceptional collections of maximal length? 
 \end{que}

\end{document}